\documentclass[aap]{imsart}

\RequirePackage{amsthm,amsmath,amsfonts,amssymb}
\RequirePackage[numbers,sort&compress]{natbib}
\RequirePackage[colorlinks,citecolor=blue,urlcolor=blue]{hyperref}
\RequirePackage{graphicx}
\usepackage{bigints}
\usepackage{comment}
\usepackage{subfig}
\startlocaldefs
\theoremstyle{plain}

\newtheorem{theorem}{Theorem}[section]
\newtheorem{lemma}[theorem]{Lemma}
\theoremstyle{remark}
\newtheorem{definition}[theorem]{Definition}
\newtheorem{example}{Example}

\newtheorem{corollary}[theorem]{Corollary}
\newtheorem{proposition}[theorem]{Proposition}
\newtheorem{conjecture}[theorem]{Conjecture}

\newtheorem{remark}[theorem]{Remark}

\newtheorem{algorithm}{Algorithm}

\DeclareSymbolFont{AMSb}{U}{msb}{m}{n}
\DeclareMathSymbol{\N}{\mathbin}{AMSb}{"4E}
\DeclareMathSymbol{\Z}{\mathbin}{AMSb}{"5A}
\DeclareMathSymbol{\R}{\mathbin}{AMSb}{"52}
\DeclareMathSymbol{\Q}{\mathbin}{AMSb}{"51}
\DeclareMathSymbol{\I}{\mathbin}{AMSb}{"49}
\DeclareMathSymbol{\C}{\mathbin}{AMSb}{"43}

\endlocaldefs

\begin{document}

\begin{frontmatter}
\title{The Inverse Problem for Single Trajectories of Rough Differential Equations}
\runtitle{The Inverse Problem for Single Trajectories of RDEs}

\begin{aug}
\author[A]{\fnms{Thomas}~\snm{Morrish}\ead[label=e1]{thomas.morrish@warwick.ac.uk}},
\author[B, C]{\fnms{Theodore}~\snm{Papamarkou}\ead[label=e2]{theodore@polyshape.com}},
\author[A]{\fnms{Anastasia}~\snm{Papavasiliou}\ead[label=e3]{a.papavasiliou@warwick.ac.uk}}
\and
\author[A]{\fnms{Yang}~\snm{Zhao}\ead[label=e4]{yz919@hotmail.com}}
\address[A]{Department of Statistics,
University of Warwick\printead[presep={,\ }]{e1,e3,e4}}

\address[B]{PolyShape\printead[presep={,\ }]{e2}}

\address[C]{School of Applied Mathematical and Physical Sciences,
National Technical University of Athens}

\end{aug}

\begin{abstract}
Motivated by the need to develop a general framework for performing statistical inference for discretely observed random rough differential equations, our aim is to construct a geometric $p$-rough path ${\bf X}$ whose response $Y$, when driving a rough differential equation, matches the observed trajectory $y$. We call this the \textit{continuous inverse problem} and start by rigorously defining its solution. We then develop a framework where the solution can be constructed as a limit of solutions to appropriately designed \textit{discrete inverse problems}, so that convergence holds in $p$-variation. Our approach is based on calibrating the bounded variation paths whose limit defines the rough path `lift' of path $X$ to rough path ${\bf X}$ to the observed trajectory $y$. Moreover, we develop a general numerical algorithm for constructing the solution to the discrete inverse problem. The core idea of the algorithm is to use the signature representation of the path, iterating between the response and the control, each time correcting according to the required properties. 

We apply our framework to the case where the geometric $p$-rough path ${\bf X}$ is defined as the limit of piecewise linear paths in the $p$-variation topology. This is a model assumption that includes some of the most common random rough path families, such as fractional Brownian motion with Hurst parameter $h>\frac{1}{4}$, as the sequence of its nested dyadic piecewise linear interpolations converge almost surely in $p$-variation, for $p>\frac{1}{h}$. We express the discrete inverse problem for a fixed observation rate as a solution to a system of equations driven by piecewise linear paths and prove convergence to the solution of the continuous inverse problem for observation time $\delta\to 0$. 
Finally, we show that, in this context, the numerical algorithm for solving the discrete inverse problem simplifies to an iterative simultaneous update of the local gradients and we prove that it converges in $p$-variation uniformly with respect to $\delta$. We discuss how we can improve the efficiency of the algorithm by considering splitting methods and reducibility conditions, and we demonstrate the results through a series of numerical examples that highlight the advantages of our algorithm in terms of robustness and efficiency.
\end{abstract}



\end{frontmatter}


\section{Introduction}

Consider the Rough Differential Equation (RDE)
\begin{equation}
\label{eq:main}
dY_t = f(Y_t) \cdot d{\bf X}_t,\ \ Y_0 = y_0,\ 0\leq t \leq T,
\end{equation}
where the control ${\bf X}\in G\Omega_p(\R^m)$ is a geometric $p$-rough path above $X$ and the vector field $f: \R^d\to L\left(\R^m,\R^d\right)$ is ${\rm Lip}(\gamma)$, for $\gamma>p\geq 1$, in the sense of E. Stein \cite[Definition 10.2, p. 213]{FrizVictoir} -- that is, $f$ can be locally approximated by a polynomial of degree $\lfloor \gamma\rfloor$ and the remainder is H\"older continuous with exponent $(\gamma-\lfloor\gamma\rfloor)$. Then, solution $Y$ exists and 
is of finite $p$-variation. The RDE \eqref{eq:main} defines an It\^o-Lyons map 
$I: \R^d\times G\Omega_p(\R^m)\to C^{p-var}([0,T],\R^d)$, 
mapping an initial value $y_0$ and $p$-rough path ${\bf X}$ to path $Y$, which is continuous in $p$-variation topology (see \cite{FrizVictoir, TerryBook} for a classical exposition or \cite{Chevyrev} for a recent review). We will write $Y = I(y_0, \bf X)$ or $Y=I({ \bf X})$, when initial conditions are implied by the context. Note that for $p<2$, $I$ is the classic It\^o map depending only on $X$.

Equation \eqref{eq:main} is very general and covers a large number of models. The initial motivating problem for this paper is that of extracting information about the vector field $f$, given a trajectory $y = \{{y}_t, t\in[0,T]\}$, of the solution (response) $Y$. This is a fundamental modeling question that comes up in numerous applications and takes many different forms, depending on the assumptions. One common assumption is that $X$ is a realization of a Brownian path and possibly time, and $f$ belongs to a parametric family, which corresponds to the well-studied problem of parametric statistical inference for diffusion processes (see \cite{Kutoyants} or \cite{Craigmile} for a review). The Brownian assumption can be relaxed, leading to the more general problem of inference about the vector field, when $X$ is a realization of a random path with known distribution \cite{HuNualart,ZhangXiao, Darrick}. Moreover, the RDE framework has been used to model recurrent neural networks, where the aim is to learn the vector field when both $X$ and $Y$ are known \cite{Kidger}. 

The dual problem to making inference about the vector field $f$ is reconstructing the control $X$ from observations of $Y$, assuming that the vector field is known -- we call this the {\it inverse problem}. It can be of independent interest, allowing one to make inference about the distribution of the random control $X$ driving the system \cite{KubiliusSkorniakov}. It also provides an indirect way of learning the vector field: in the case where the control is a realization of a random path with known distribution, reconstructing the specific realization of the control $X$ that is consistent with the observations, conditioned on the vector field $f$ can lead to the construction of the likelihood. 

The inverse problem for RDEs has also been studied in \cite{BailleulDiehl}, but in a different context: the authors reconstruct the truncated signature of the control ${\bf X}$ on a fixed window from observations of the increments of the solution to the RDE on that window, for a number of different initial conditions. Instead, we focus on reconstructing the rough path ${\bf X}$ from a trajectory of $Y$, which is a more realistic assumption for most applications. Note that, in the rough path framework, observing a trajectory of $Y$ is equivalent to continuously observing the increments of rough path ${\bf Y}$ but not the rough path `lift'. 
Thus, for $p\geq 2$, the trajectory of $Y$ on its own does not contain sufficient information to reconstruct the rough path ${\bf X}$. However, the missing information can be extracted from the model and, in particular, from the definition of the geometric $p$-rough path ${\bf X}$ as the limit of bounded variation paths in $p$-variation topology.

In section \ref{sec: problem}, we argue that in order to be able to reconstruct the rough path ${\bf X}$ from continuous observations of the increments of $Y$, we can start by considering finite-dimensional families of bounded variation paths that we know converge in $p$-variation to ${\bf X}$ when fitted to its trajectory and then fit them instead to the observations. This leads to the definition of the {\it discrete inverse problem} that involves 
fitting the finite-dimensional family of bounded variation paths approximating $X$ to the trajectory of $Y$ on a finite partition of size $\delta$ (the {\it calibration step}). We then give general conditions that the model needs to satisfy, so that solutions to the discrete problem are good approximations of solutions to the continuous problem in $p$-variation, in the sense that the distance between any solution to the discrete inverse problem and the space of solutions to the continuous inverse problem goes to $0$, as $\delta\to 0$. Note that we do not assume the uniqueness of solutions. Instead, we assume that the It\^o-Lyons map has a weakly continuous inverse, in a sense that we make precise, which is connected to the local invertibility property of the It\^o-Lyons map, discussed in \cite{Bailleul}.

In section \ref{sec: piecewise}, we focus on the case where piecewise linear interpolations of $X$ converge in $p$-variation to ${\bf X}$. Although in this case piecewise linear interpolations of $Y$ also converge, making it possible to define the lift ${\bf Y}$ first and then solve the inverse problem directly, working directly with $X$ still has an advantage as it allows us to combine the discretely observed response $Y$ with information about the distribution of $X$, which can improve results.
For example, if $X$ is a Gaussian process, the distribution of piecewise linear interpolations reduces to a finite-dimensional Gaussian distribution, which allows us to either express its likelihood precisely or perform exact simulation of the interpolation at a finer scale.

When considering piecewise linear approximations to $X$, the corresponding discrete inverse problem becomes equivalent to finding a piecewise linear path $\hat{X}^\delta$ that, when driving \eqref{eq:main}, is consistent with the observations. This leads to a localized formulation of the calibration step: Each linear segment of the piecewise linear path $\hat{X}^\delta$ is a solution to a system of equations derived by requiring that the terminal value of the solution to an Ordinary Differential Equation (ODE) matches the observations. We derive easy-to-check conditions for the solutions of the discrete inverse problem to converge to the solution of the continuous inverse problem in $p$-variation, under the assumption that the It\^o-Lyons map has a weakly continuous inverse.


An exact solution to the discrete inverse problem will rarely be available analytically. In section \ref{sec: numerical}, we discuss different numerical methods for constructing the solution. The first is an application of the Newton-Raphson method on the solution of the ODE as a function of the gradient of the control. The second method is based on a novel approach that relies on the signature representation of the path; by reformulating the discrete inverse problem in terms of the pair $({\bf X},{\bf Y})$, the main idea is to iterate between the driving rough path ${\bf X}$ and the rough path lift ${\bf Y}$ of the solution $Y$ through the corresponding It\^o-Lyons and inverse It\^o-Lyons maps, each time correcting the paths ${X}$ and ${Y}$ so that they satisfy the requirements to be a solution to the discrete inverse problem, in a way that changes the corresponding pair of signatures $({\bf X},{\bf Y})$ as little as possible. The main advantage of this approach is that it does not rely on the ODE solution but rather on the It\^o-Lyons map. As such, it can be extended to other types of approximations. Moreover, it sheds light on when and why the inverse problem can be solved exactly by quantifying the error in terms of the information contained in the discrete observations and that contained in the underlying continuous path. Under the assumption of piecewise linear approximations to $X$, we are able to reformulate the algorithm in terms of the gradients of each linear segment, and we show that the algorithm converges in $p$-variation, uniformly with respect to $\delta$.
In section \ref{sec: extensions}, we discuss how to further improve the efficiency of the algorithm by using reducibility conditions and splitting methods, building on ideas in \cite{ait_sahalia,Foster,splitting} and in section \ref{sec: examples}, we present numerical experiments in which we compare the two numerical algorithms in terms of robustness and efficiency.

\subsection{Notation}
Throughout the paper, we will use the following notation:
\begin{itemize}
    \item Bold font refers to the signature of a trajectory, a tensor object. For example, $\bf X$ refers to the signature, or the rough path lift of the trajectory $X$. Sometimes, when we are dealing with a path of bounded variation (and the signature is defined uniquely by the trajectory), we may use the trajectory and signature (and the It\^o and It\^o-Lyons maps) interchangeably.  
    \item Lower case $y$ refers to a sample trajectory of $Y$ or ${\bf Y}$, over a default time interval $[0,T]$. Thus, if we had a sample signature over $[0,T]$, ${\bf y}$, then we would have that $y={\bf y}^1+y_0$, as ${\bf y}^1$ records the increments of the path.
    \item For rough path ${\bf Z}$, $Z_t$ refers to the path $Z$ evaluated at the time $t$, while ${\bf Z}^k_{s,t}$ refers to the $k^{\rm th}$ level of the signature of the path $Z$ between times $s$ and $t$. For $k=1$, we also use the notation ${Z}_{s,t}=Z_t-Z_s$. 
\end{itemize}

\section{General framework}
\label{sec: problem}

We start by defining the continuous inverse problem for single trajectories of rough differential equations.

\begin{definition} [Continuous Inverse Problem] 
\label{def:cont inverse problem}
Suppose that we are given a trajectory ${y} = \{{y}_t, t\in[0,T]\}$ in $\R^d$ that is of finite $p$-variation. We will say that a $p$-rough path ${\bf X}\in G\Omega_p(\R^m)$ is a solution to the {\it continuous inverse problem} if the trajectory of the solution $Y$ to \eqref{eq:main} driven by ${\bf X}$ is the same as $y$, i.e., $Y_t = {y}_t, \forall t\in [0,T]$ for some fixed $T>0$.
\end{definition}

Let us assume for a moment that $m=d$, and $g:\R^{d}\to L(\R^d,\R^d)$ defined as $g(y) = f(y)^{-1}$ exists and is ${\rm Lip}(\gamma)$ for $\gamma>p\geq 1$. Then, given a trajectory $y = \{{y}_t, t\in[0,T]\}$, we would expect that 
\begin{equation} 
\label{eq: naive inverse}
\int_0^tg({y}_u) \cdot d{y}_u 
\end{equation}
is a natural candidate for $X$. Indeed, if $p<2$, integral \eqref{eq: naive inverse} is well-defined as a Young integral and the trajectory of its response through \eqref{eq:main} will match $y$. However, if $p\geq 2$, integral \eqref{eq: naive inverse} is not well-defined in general. Instead, we need to construct a geometric $p$-rough path ${\bf Y}$ whose trajectory is the same as $y$ (i.e. $\int_s^t dY_u = {y}_t - {y}_s$ for all $0\leq s<t\leq T$). Then,
\begin{equation}
\label{eq: exact inverse}
\int_0^t g(Y_u) \cdot d{\bf Y}_u, 
\end{equation}
interpreted as a rough path integral, is well-defined and would be a solution to the continuous inverse problem.

Constructing a geometric $p$-rough path ${\bf Y}$ that corresponds to the trajectory ${y}$, for $p\geq 2$, requires the construction of the rough path `lift' of $Y$, i.e. all iterated integrals 
\begin{equation}
\label{eq: lift}
{\bf Y}^{k}_{s,t} = \int_s^t {\bf Y}^{k-1}_{s,u}\otimes d{\bf Y}_u
\end{equation}
for $k=1, 2,\dots,\lfloor p \rfloor$, for all $0\leq s < t \leq T$.

This leads to a more general question: given a trajectory $z:[0,T]\to\R^d$ of finite $p$-variation, how can we construct a geometric $p$-rough path ${\bf Z}$ that corresponds to $z$? If we can find (i) a finite-dimensional family of bounded variation paths $Z^{\delta}:\R^{M_\delta}\to {\rm BV}(\R^d)$, parameterized, for each $\delta$, by a vector of free parameters $\theta^{\delta}\in\R^{M_\delta}$; and (ii) a family of maps ${\ \Theta}^\delta:C^{p-var}\left( [0,T], \R^d \right)\to \R^{M_\delta}$, such that the $Z^{\delta}({\Theta}^{\delta}(z))$ converges in the $p$-variation topology as $\delta \to 0$ to a $p$-rough path ${\bf Z}$ and ${\bf Z}^1_{s,t} = z_t - z_s$ for every $0\leq s < t \leq T$, then ${\bf Z}$ is the $p$-rough lift corresponding to $z$. We will then say that $\Theta^\delta$ is the calibration function and $\theta^\delta_z = \Theta^\delta(z)$ is the parameter calibrated to $z$. Note that the calibration function does not need to be unique for a given family $Z^\delta(\theta^\delta)$ and that our choice will depend on the context.

As the concept of calibration to the data is fundamental for the rest of the paper, we give below a couple of illustrative examples:


\begin{example}[Piecewise linear approximation]
\label{ex: pwl}
Throughout this paper, the finite-dimensional family of bounded variation paths of our choice for constructing the rough path lift will be the family of piecewise linear paths over the grid ${\mathcal D}_\delta = \{t_k = k\delta; k=0,\dots,N_\delta\}$, with $N_\delta \delta = T$ fixed, where the free parameters are the gradients of the path on intervals of length $\delta$, which have dimension $M_\delta = d \cdot N_\delta$. That is
\[
    Z^\delta(\theta^\delta)_u
    = Z^\delta(\theta^\delta)_{t_i} + \theta^\delta_i\,(u-t_i).
\]
Given a trajectory $z$, an obvious choice for the calibration function is 
\[
\theta^\delta_z = \Theta^\delta (z) = \bigg\{ \frac{z_{t_{i+1}}-z_{t_i}}{\delta}; i=0,\dots,N_\delta-1\bigg\}.
\]
Then, for $u\in[t_i,t_{i+1}]$, and a trajectory $z$, we have
\[
    Z^\delta(\theta^\delta_z)_u
    = Z^\delta(\theta^\delta_z)_{t_i} + \frac{z_{t_{i+1}}-z_{t_i}}{\delta}\,(u-t_i)=z_{t_i} + \theta^\delta_{z,i}\,(u-t_i).
\]
Since $Z^\delta(\theta^\delta)$ are bounded-variation paths, constructing the $p$-rough path lift \eqref{eq: lift} for $Z^\delta(\theta^\delta_z)$ amounts to computing Riemann--Stieltjes integrals. If the limit exists in $p$-variation topology for $\delta\to0$, then that limit is a geometric $p$-rough path corresponding to $z$ and we denote it by ${\bf Z}$. While we will focus on the case where the finite-dimensional family of bounded variation paths is that of piecewise linear paths, in this paper we also consider alternative calibration functions. 
\end{example}

\begin{example}[Karhunen-Lo\`eve decomposition]
Below, we give an example of an alternative finite-dimensional family of bounded variation paths to the piecewise linear case, that can be used to construct the rough path lift. While we will not develop this further, they also fit the general framework described in this section.

Suppose that $z\in L^2([0,T])$ and $\{e_i\}_{i\ge 1}$ is an orthonormal basis of $L^2([0,T])$, such that 
\[
    z_t = \lim_{M\to\infty}\sum_{i=1}^M \phi_i e_i(t),
\]
where the limit is to be interpreted as the $L_2$ limit. An example of such a $z$ would be a trajectory of the zero-mean stochastic process that admits the Karhunen--Lo\`eve decomposition (see, e.g., \cite{loeve_book}). Then, if the truncated sums
\[
Z^M(\theta^M)_t = \sum_{i=1}^M \theta^M_i e_i(t).
\]
are of bounded variation, they provide a finite-dimensional family of bounded-variation paths that can be used to construct the lift. Moreover, we can choose
\[
    \theta^M_{z,i} = \Theta^M(z)_i = \int_0^T z_u e_i(u) du
\]
as our calibration function. Then, if the $p$-variation limit of $Z^M(\theta^M_z)$ exists, then it defines the $p$-rough path lift ${\bf Z}$ that corresponds to $z$.

In this example, we used the $Z^M$ and $\theta^M$ rather than $Z^\delta$ and $\theta^\delta$, as it is natural to take $M\to\infty$ rather than $\delta\to 0$. For the rest of the paper, we will stick to the $Z^\delta,\theta^\delta$ notation, with $M_\delta\to\infty$ as $\delta \to 0$.
\end{example}

Let us now go back to the original problem of constructing a solution to the continuous inverse problem. The approach described earlier, aiming to directly compute \eqref{eq: exact inverse}, has several shortcomings: an obvious one is the assumption regarding the existence of $g=f^{-1}$. Another, less obvious but more important one is that we do not always have prior information about how to choose the finite-dimensional family of bounded variation paths ${Y}^\delta({\bf \theta}^\delta)$ and calibration function $\Theta^\delta$, such that ${Y}^\delta(\Theta^\delta(y))$ converges in the $p$-variation topology, as $Y$ is defined as a response to ${\bf X}$ through \eqref{eq:main}. We do, however, expect to know this for ${\bf X}$ which is defined directly as a geometric $p$-rough path. Thus, rather than trying to solve the inverse problem going `backwards' with $y$ as our starting point by aiming to construct the integral \eqref{eq: exact inverse}, we start with the finite-dimensional family and bounded-variation paths ${X}^\delta(\theta^\delta)$ and we go `forward', by now calibrating the parameter $\theta^\delta$ to the trajectory of the response $Y$ rather than that of $X$. This addresses both problems of constructing the rough path lift and assuming the existence of $g$, leading to the following definition of the {\it discrete inverse problem}.

\begin{definition} [Discrete Inverse Problem] 
\label{def:discrete inverse problem}
Let $y = \{y_t, t\in[0,T]\}$ be a continuous path in $\R^d$ of finite $p$-variation. Let $X^\delta(\theta^\delta)$ be a finite-dimensional family of bounded variation paths parameterized by $\theta^\delta$, and $\Theta^\delta_X:C^{p-var}([0,T],\R^m)\to\R^{M_\delta}$ a calibration function, such that if we know the trajectory ${x}$ of $X$ then for $\theta^\delta_x = \Theta^\delta_X(x)$, ${\bf X}^\delta(\theta_x^\delta)\to {\bf X}$ as $\delta\to 0$ in the $p$-variation topology. 

Suppose that there exists a calibration function $\Theta^\delta_Y:C^{p-var}([0,T],\R^d)\to\R^{M_\delta}$ on the space of trajectories $y$, defined as the solution with respect to $\theta^\delta$ of equation 
\begin{equation} 
\label{eq: y calibration}
I({X}^\delta(\theta^\delta))_{t_i} = y_{t_i},\ \forall t_i\in{\mathcal D}_\delta,
\end{equation}
where ${\mathcal D}_\delta = \{k\delta; k=0,\dots,N_{\delta}, N_{\delta}\delta = T\}$ and $I$ is the It\^o map defined by \eqref{eq:main}. That is, for $\theta^\delta_y=\Theta_Y^\delta(y)$ and $\hat{X}^\delta = X^\delta(\theta^\delta_y)$, $\hat{Y}^\delta = I(\hat{X}^\delta)$ is the solution to
\begin{equation}
\label{eq:approx}
d\hat{Y}^\delta_t = f(\hat{Y}^\delta_t) \cdot d{X}^\delta(\theta^\delta_y)_t,\ 0\leq t \leq T, \ \ \hat{Y}^\delta_0 = Y_0.
\end{equation}
Then, we will say that $\hat{X}^\delta = X^\delta(\theta_y^\delta)$ solves the discrete inverse problem corresponding to the finite-dimensional family of bounded variation paths $X^\delta(\theta^\delta)$. Note that, for fixed $\delta>0$, the solution to the discrete inverse problem requires only discrete observations of ${y}$ on ${\mathcal D}_\delta$.
\end{definition}

\begin{remark}
For different discrete information about $y$, it is possible to extend this framework by replacing \eqref{eq: y calibration} to fit the given information. Moreover, in principle it could be possible to extend this framework to non-geometric rough paths ${\bf X}$, but it would require an appropriate interpretation of $X^\delta(\theta^\delta)$ and $I$ in \eqref{eq: y calibration}, which this is beyond the scope of this paper.
\end{remark}

\begin{remark}
It is useful to note that solutions to the discrete inverse problem are trajectories, while solutions to the continuous inverse problem are rough path lifts.  
\end{remark}

The key step in constructing the solution $X^{\delta}(\theta_y^\delta)$ to the discrete inverse problem is the computation of $\theta_y^\delta = \Theta^\delta_Y(y)$, which we call the `calibration step'. In practice, this requires solving the system of equations defined in \eqref{eq: y calibration}, which also involve solutions to the ODE defined in \eqref{eq:approx}. 

A desirable property for the solutions to the discrete inverse problem is that they converge to the solution to the continuous inverse problem in $p$-variation as $\delta\to 0$. We discuss below conditions for this to happen.

\subsection{Convergence of solutions to the discrete inverse problem}
\label{subsec: discrete to continuous}

Let $\hat{X}^\delta = X^\delta(\theta_y^\delta)$ be a solution to the discrete inverse problem for some continuous path ${y}$ and $\delta>0$. We would like to show that for $\delta$ small enough, its rough path lift $\hat{\bf X}^\delta = {\bf X}^\delta(\theta_y^\delta)$ will be within $\epsilon$ $p$-variation distance from a solution to the continuous inverse problem ${\bf X}$. Note that we do not assume uniqueness of solution to either the discrete or continuous inverse problem. Thus, the previous statement is equivalent to showing that the distance between $\hat{\bf X}^\delta$ and the space of solutions to the continuous inverse problem, which we will denote by ${\mathcal I}^{-1}({y})$, goes to $0$ with $\delta \to 0$, provided that ${\mathcal I}^{-1}({y})\neq \emptyset$ (there exists at least one solution to the continuous inverse problem), i.e.
\begin{equation}
\label{eq: def of convergence}
\lim_{\delta\to 0}d_p({\bf X}^\delta(\theta_y^\delta),{\mathcal I}^{-1}({y}))=0,
\end{equation}
where $d_p({\bf X},A) = \inf_{\tilde{\bf X}\in A}\|{\bf X}-\tilde{\bf X}\|_{p-{\rm var},[0,T]}$ for $A\subset G\Omega_p(\R^d)$. We will use $\|\cdot\|_{p-{\rm var},[0,T]}$ to denote the $p$-variation distance of either $p$-rough path paths (i.e. including the rough path lift) or paths (without the lift) on the corresponding Euclidean space, depending on the context. We denote by $|\cdot |$ the corresponding Euclidean metric.

Proving convergence of ${\bf X}^\delta(\theta^\delta_y)$ to a solution of the continuous inverse problem, in the sense of \eqref{eq: def of convergence}, requires two steps: (i) showing that $I({X}^\delta(\theta^\delta_y))$ and $I({X}^\delta(\theta^\delta_x))$ both converge to $Y$ in $p$-variation and (ii) showing that if $I({X}^\delta(\theta^\delta_y))$ and $I({X}^\delta(\theta^\delta_x))$ are close in $p$-variation, the controls ${X}^\delta(\theta^\delta_y)$ and ${X}^\delta(\theta^\delta_x)$ will also be close, in the sense that there exists a solution ${\bf X}$ to the continuous problem such that $X^\delta(\theta^\delta_y)$ calibrated to the observations and $X^\delta(\theta^\delta_x)$ calibrated to the trajectory of the solution ${\bf X}$ will be close. The second step is equivalent to a type of inverse function theorem for the It\^o-Lyons map. This will depend on specific model assumptions and, in particular, the family of parameterized paths $X^\delta(\theta^\delta)$ that converges in $p$-variation almost surely, under the assumptions on the distribution of $X$. Thus, at this stage, we will assume it holds by assuming that $I$ has a weakly continuous inverse, in the following sense:

\begin{definition}
\label{def:inv Ito continuity}
Let $I$ be the It\^o-Lyons map $I:G\Omega_p(\R^m)\to {C}^{p-var}([0,T],\R^d)$ corresponding to \eqref{eq:main}. We will say that $I$ has a weakly continuous inverse if there exists an $\epsilon>0$ and a constant $C>0$ depending only on $f$ and $p$, such that for every ${\bf X}\in G\Omega_p(\R^m)$ and $Y\in C^{p-var}([0,T],\R^d)$ with $I({\bf X})$ and $Y$ within $\epsilon$ $p$-variation distance, for $\epsilon$ sufficiently small,
\begin{equation}
\label{eq: Ito lower bound}
\inf_{\tilde{\bf X}: I(\tilde{\bf X}) = Y} \|{\bf X} - \tilde{\bf X}\|_{p-{\rm var},[0,T]} \leq C \| I({\bf X}) - Y \|_{p-{\rm var},[0,T]}.
\end{equation}
\end{definition}

\begin{remark} 
Property \eqref{eq: Ito lower bound} is connected to the local invertibility of the It\^o-Lyons map. Conditions for that to hold can be found in \cite{Bailleul}. However, a constructive proof of such a result is still missing.
\end{remark}


Going back to the first step required for the proof, it follows from the continuity of the It\^o-Lyons map (Universal Limit Theorem) and the convergence of ${X}^\delta(\theta^\delta_x)$ to ${\bf X}$ in $p$-variation, by definition, that $I({X}^\delta(\theta^\delta_x))$ will converge to $Y$. By construction, $I(X^\delta(\theta^\delta_y))$ will converge to ${y}$ pointwise. The following lemma says that as long as we can control the distance of $X^\delta(\theta^\delta_x)$ and $X^\delta(\theta^\delta_y)$ in-between observations, we will be able to control their responses on the whole time horizon in $p$-variation. In other words, $I(X^\delta(\theta^\delta_y))$ will also converge to $Y$ in $p$-variation.

\begin{lemma}
\label{lemma: discrete solution Y convergence}
Let $X^{\delta}(\theta^\delta_y)$ be a solution to the discrete inverse problem and ${\bf X}$ be a solution to the continuous inverse problem, with $\theta^\delta_x$ such that $X^{\delta}(\theta^\delta_x) \to {\bf X}$ in $p$-variation. Let $\omega:\Delta_T = \{0\leq s < t\leq T\}\to \R$ be a control, such that 
\[
|{\bf X}^{\delta}(\theta^\delta_y)_{s,t}^{k}|,\ |{\bf X}^{\delta}(\theta^\delta_x)_{s,t}^{k}|,\ |{\bf X}_{s,t}^{k}| <\omega(s,t)^\frac{k}{p},
\]
for all $0\leq s<t \leq T$ and $k=1,\dots,\lfloor p \rfloor$. Moreover, we assume that
\begin{equation}
\label{eq: in-segment smoothness}
\| {\bf X}^\delta(\theta^\delta_y) - {\bf X}^\delta(\theta^\delta_x)\|_{p-{\rm var},[t_i,t_{i+1}]} \leq K \hat{\eta}(\delta)\omega(t_i,t_{i+1})^\frac{1}{p},
\end{equation}
for every $[t_i,t_{i+1}]$, where $K>0$ is a constant independent of $i$ and $\hat{\eta}(\delta)$ is such that $\hat{\eta}(\delta)\to 0$ as $\delta \to 0$. 
Then, there exists a constant $C>0$ and a function ${\eta}(\delta)$ with ${\eta}(\delta)\to 0$ as $\delta \to 0$, such that 
\begin{equation}
\label{eq: xy distance}
  | I({\bf X}^{\delta}(\theta^\delta_y))_{s,t} - I({\bf X}^{\delta}(\theta^\delta_x))_{s,t} | \leq C \eta(\delta)\omega(s,t)^\frac{1}{p},
\end{equation}
for all $0\leq s<t \leq T$, where all constants depend only on $f,p,\gamma$ and $\omega(0,T)$.
\end{lemma}

\begin{proof}
For ease of notation, we set $\hat{Y}^\delta = \hat{Y}^\delta(y) = I(X^{\delta}(\theta^\delta_y))$ and $\hat{Y}^\delta(x) = I(X^{\delta}(\theta^\delta_x))$. Moreover, $Y=I({\bf X})$ as before. We prove \eqref{eq: xy distance} by considering the following different cases:
    
    First, we consider the case where $s=t_i,t = t_j \in {\mathcal D}_\delta$. By construction, $I(X^{\delta}(\theta^\delta_y))_{t_i} = \hat{Y}^\delta(y)_{t_i} = Y_{t_i}$ for every $t_i\in{\mathcal D}_\delta$. It follows that $\hat{Y}^{\delta}(y)_{t_i,t_j} = Y_{t_i,t_j}$. Moreover, since $X^{\delta}(\theta^\delta_x)\to {\bf X}$ in the $p$-variation topology, it follows from the Universal Limit Theorem that 
\[
|\hat{Y}^{\delta}(x)_{s,t}-Y_{s,t}| < C_1\eta_1(\delta) \omega(s,t)^\frac{1}{p},
\]
for some $C_1>0$ and $\eta_1(\delta) = \|{\bf X}^{\delta}(\theta^\delta_x)-{\bf X}\|_{p-{\rm var},[0,T]}$ (see \cite[Thm 6.3.1, p. 179]{TerryBook}). By assumption, $\eta_1(\delta)\to 0$ as $\delta\to 0$.
Thus,
\[
|\hat{Y}^{\delta}(y)_{t_i,t_j} - \hat{Y}^{\delta}(x)_{t_i,t_j}| = |Y_{t_i,t_j}-\hat{Y}^{\delta}(x)_{t_i,t_j}| < C_1 \eta_1(\delta) \omega(t_i,t_j)^\frac{1}{p},
\]
as required.

Next, we consider any $0\leq s<t\leq T$. Let $t_i, t_j \in {\mathcal D}_\delta$ be such that $t_{i-1}<s\leq t_i$ and $t_j \leq t < t_{j+1}$. It is sufficient to consider the case where $t_i\leq s<t\leq t_{i+1}$, as the general case follows from the previous result and 
\begin{eqnarray*}
    |\hat{Y}^{\delta}(y)_{s,t} - \hat{Y}^{\delta}(x)_{s,t}| &\leq& |\hat{Y}^{\delta}(y)_{s,t_i} - \hat{Y}^{\delta}(x)_{s,t_i}| +|\hat{Y}^{\delta}(y)_{t_i,t_j} - \hat{Y}^{\delta}(x)_{t_i,t_j}|\\
    && +|\hat{Y}^{\delta}(y)_{t_j,t} - \hat{Y}^{\delta}(x)_{t_j,t}|.
\end{eqnarray*}
Let $u\in[0,T]$. We denote by $\tilde{Y}^{u,\delta}$ the response to $X^\delta(\theta^\delta_y)$ through \eqref{eq:main} with initial condition $\tilde{Y}^{u,\delta}_{u} = \hat{Y}^\delta(x)_{u}$. Then,
\begin{equation}
\label{eq: in-segment bound}
   |\hat{Y}^{\delta}(y)_{s,t} - \hat{Y}^{\delta}(x)_{s,t}| \leq |\hat{Y}^{\delta}(y)_{s,t} - \tilde{Y}^{s,\delta}_{s,t}| + |\tilde{Y}^{s,\delta}_{s,t} - \hat{Y}^{\delta}(x)_{s,t}| .
\end{equation}
The first term on the right-hand side of the inequality above is the difference of two processes driven by the same path but with different initial conditions and it is bounded by
\begin{eqnarray}
\label{eq: same initial conditions}
\nonumber   |\hat{Y}^{\delta}(y)_{s,t} - \tilde{Y}^{s,\delta}_{s,t}| &\leq& C_2 |\hat{Y}^{\delta}(y)_{s} - \hat{Y}^{\delta}(x)_{s}| \cdot |\hat{X}^\delta(\theta^\delta_y)_{s,t}| \\
   &\leq& C_2 |\hat{Y}^{\delta}(y)_{s} - \hat{Y}^{\delta}(x)_{s}|\omega(s,t)^\frac{1}{p}.
\end{eqnarray}
See \cite[Lemma 10.6, p. 216]{FrizVictoir}. 
Moreover, we write
\begin{eqnarray*}
|\hat{Y}^{\delta}(y)_{s} - \hat{Y}^{\delta}(x)_{s}| &\leq& |\hat{Y}^{\delta}(y)_{s} - \tilde{Y}^{t_i,\delta}_{s}|+|\tilde{Y}^{t_i,\delta}_{s}-\hat{Y}^{\delta}(x)_{s}| \\
& \leq & C_3|\hat{Y}^{\delta}(y)_{t_i} - \hat{Y}^{\delta}(x)_{t_i}|\cdot |X^\delta({\theta}^\delta_y)_{t_i,s}| + |\tilde{Y}^{t_i,\delta}_{t_i,s}-\hat{Y}^{\delta}(x)_{t_i,s}| \\
& \leq & C_3|\hat{Y}^{\delta}(y)_{0,t_i} - \hat{Y}^{\delta}(x)_{0,t_i}| \omega(t_i,s)^\frac{1}{p} \\
&&+ C_4 \|{\bf X}^\delta(\theta^\delta_y)-{\bf X}^\delta({\theta}^\delta_x)\|_{p-{\rm var,[t_i,t_{i+1}]}}\omega(t_i,s)^\frac{1}{p}\\
& \leq & C_3 \eta_1(\delta) \omega(0,t_i)^\frac{1}{p}
\omega(t_i,s)^\frac{1}{p} + C_4 K \hat{\eta}(\delta)\omega(t_i,t_{i+1})^\frac{1}{p}\omega(t_i,s)^\frac{1}{p} \leq C_5 \eta_2(\delta),
\end{eqnarray*}
for appropriately defined $\eta_2(\delta)$, such that $\eta_2(\delta)\to 0$ as $\delta \to 0$. Thus, the first term of the upper bound in \eqref{eq: in-segment bound} is bounded by
\[ 
|\hat{Y}^{\delta}(y)_{s,t} - \tilde{Y}^{s,\delta}_{s,t}| \leq C_6 \eta_2(\delta) \omega(s,t)^\frac{1}{p}.
\]
Applying once again the Universal Limit Theorem on interval $[t_i,t_{i+1}]$ (see \cite[Theorem 2.9, p. 73]{Baudoin}) and assumption \eqref{eq: in-segment smoothness}, we get the following bound for the second term of \eqref{eq: in-segment bound},
\begin{equation*}
       |\tilde{Y}^{s,\delta}_{s,t} - \hat{Y}^{\delta}(x)_{s,t}| \leq C_7 \|{\bf X}^\delta({\theta}^\delta_y)-{\bf X}^\delta({\theta}^\delta_x)\|_{p-{\rm var,[t_i,t_{i+1}]}}\omega(s,t)^\frac{1}{p} 
    \leq C_8 \hat{\eta}(\delta)\omega(s,t)^\frac{1}{p}.  
\end{equation*}
Putting everything together, we can find $C$ and $\eta(\delta)$ with $\eta(\delta)\to 0$ as $\delta\to 0$, such that \eqref{eq: xy distance} holds for any $0\leq s < t \leq T$.

\end{proof}


We can now prove the convergence in $p$-variation of the solutions to the discrete inverse problem $X^{\delta}({\theta}^\delta_y)$ to a solution of the continuous inverse problem ${\bf X}$.

\begin{theorem}
\label{thm: general convergence}
Let $Y\in C^{p-var}([0,T],\R^d)$ be the solution to \eqref{eq:main} with trajectory ${y}$ and let $X^{\delta}({\theta}^\delta_y)$ be a corresponding solution to the discrete inverse problem. Suppose that there exists an ${\bf X}\in {\mathcal I}^{-1}({y})$, such that the corresponding $X^{\delta}({\theta}^\delta_x)$ satisfies \eqref{eq: in-segment smoothness}. Then, if $I$ has a weakly continuous inverse, 
\begin{equation}
\label{eq: convergence of X(cy)}
\lim_{\delta\to 0}d_p\left({\bf X}^\delta({\theta}_y^\delta),{\mathcal I}^{-1}({y})\right)=0,
\end{equation}
where ${\mathcal I}^{-1}({y}) = \{ {\bf X}\in G\Omega_p(\R^d)\ |\ I({\bf X})_t = y_t,\ \forall t\in[0,T] \}$.
\end{theorem}

\begin{proof}
It follows from lemma \ref{lemma: discrete solution Y convergence} that we can take $\delta$ small enough, so that $I({\bf X}^\delta(\theta^\delta_y))$ and $I({\bf X})=Y$ are sufficiently close. Thus, we can use the assumption that $I$ has a weakly continuous inverse to write
\begin{eqnarray*}
 & d_p({\bf X}^\delta({\theta}_y^\delta),{\mathcal I}^{-1}(y)) =   \inf_{\tilde{\bf X}: I(\tilde{\bf X}) = Y} \|{\bf X}^\delta({\theta}_y^\delta) - \tilde{\bf X}\|_{p-{\rm var},[0,T]}  \\
&  \leq  C \cdot \|I(X^\delta({\theta}_y^\delta)) - Y\|_{p-{\rm var},[0,T]} = C \cdot \|I(X^\delta({\theta}_y^\delta)) - I({\bf X})\|_{p-{\rm var},[0,T]} \\
 & \leq  C \left( \|I(X^\delta({\theta}_y^\delta)) - I(X^\delta({\theta}_x^\delta))\|_{p-{\rm var},[0,T]} + \|I(X^\delta({\theta}_x^\delta)) - I({\bf X})\|_{p-{\rm var},[0,T]}\right),
\end{eqnarray*}
also using the existence of ${\bf X}$ such that $I({\bf X})=Y$ and that the corresponding $X({\theta}^\delta_x)$ satisfies \eqref{eq: in-segment smoothness}. It follows from lemma \ref{lemma: discrete solution Y convergence} that the first term goes to $0$ as $\delta \to 0$. The second term also goes to $0$ by construction and the Universal Limit Theorem, which proves the result.
\end{proof}


Having developed a general framework for addressing the inverse problem for single (discretely observed) trajectories of solutions to rough differential equations, we now focus on the case where ${X}^\delta$ is the family of piecewise linear paths on $\mathcal{D}_\delta$, parameterized by the gradients (i.e. we now write $\theta^\delta=c^\delta$ for gradients over intervals of width $\delta$), under the assumption that piecewise linear interpolations of $X$ converge in $p$-variation. However, we note that this framework allows one to also address cases where the piecewise linear interpolations of $X$ do not converge and one needs to consider different families of approximations (e.g., \cite{Unterberger}).


\section{Piecewise Linear Approximations}
\label{sec: piecewise}

Suppose that piecewise linear interpolations of the control ${\bf X}\in G\Omega_p(\R^m)$ driving \eqref{eq:main} on a sequence of nested dyadic partitions (i.e. $\delta = 2^{-n}$ for $n\geq 1$) converge in $p$-variation. That is, if $X^\delta({c}^\delta)$ is the family of piecewise linear paths starting at $X_0$ (we will assume that $X_0=0$) and parameterized by the gradients of its linear segments denoted by ${c}^\delta_x = \{ c^\delta_{x,k}\}_{k=1}^{N_\delta}$, with $c^\delta_{x,k} = \frac{1}{\delta}\left(X_{k\delta}-X_{(k-1)\delta}\right)$ so that $X^\delta({c}_x^\delta)$ interpolates $X$ on partition ${\mathcal D}_\delta$, then $X^\delta({c}^\delta_x)\to {\bf X}$ in $p$-variation. For Random Rough Differential Equations (RRDEs), we want convergence to hold almost surely, so this is a distributional assumption and it holds in many cases of interest -- in particular, it holds when $X$ is a fractional Brownian motion for $h>\frac{1}{4}$ (see \cite{TerryBook}). Note that in order to keep notation simple and consistent, we will keep working with $\delta\to 0$ rather than $n\to\infty$.

Having secured the existence of at least one sequence ${c}^\delta = {c}_x^\delta$, such that $X^\delta({c}^\delta_x)$ converges in $p$-variation to ${\bf X}$, we apply the framework for constructing the solution to the discrete inverse problem to this parameterization. Under this choice of parameterization for $X$, the discrete inverse problem (definition \ref{def:discrete inverse problem}) takes the following form:

\begin{definition} [Discrete Inverse Problem -- piecewise linear paths] 
\label{def: pwl discrete inverse problem}
Consider the family of piecewise linear paths $X^\delta({c}^\delta)$ on partition ${\mathcal D}_\delta$, parameterized by the gradients ${c}^\delta$. Given a vector ${y}_{{\mathcal D}_\delta} = \left\{y_{k\delta}\right\}_{k=0}^{N_\delta}$ and a vector field $f\in{\rm Lip}(\gamma)$, for $\gamma>p\geq 1$, we say that $\hat{X}^\delta = X^\delta({c}_y^\delta)$ solves the discrete inverse problem for this parameterization if there exists a ${c}_y^\delta$ such that its response $\hat{Y}^\delta$ through \eqref{eq:main} interpolates the observations, i.e. $\hat{Y}^\delta_{k\delta} = y_{k\delta}, \forall k=0,\dots,N_\delta$. Equivalently, we say that $\left( \hat{X}^\delta,\hat{Y}^\delta\right)$ solves the discrete inverse problem for piecewise linear paths if the following conditions are satisfied:
\begin{itemize}
\item[(a)] $\hat{X}^\delta = X^\delta({c}_y^\delta)$ is a piecewise linear path on ${\mathcal D}_\delta$ with gradients ${c}_y^\delta$, starting at $\hat{X}_0 = 0$.
\item[(b)] $\hat{Y}^\delta = I(\hat{X}^\delta)$, where $I$ is the It\^o map defined by \eqref{eq:main}, i.e. $\hat{Y}^\delta$ is the solution to
\begin{equation}
\label{eq:approx2}
d\hat{Y}^\delta_t = f(\hat{Y}^\delta_t) \cdot d\hat{X}^\delta_t,\ 0\leq t \leq T, \ \ \hat{Y}^\delta_0 = y_0.
\end{equation}
\item[(c)] $\hat{Y}^\delta_{k\delta} = y_{k\delta},\ \forall k=0,1,\dots,N_\delta$.
\end{itemize}
\end{definition}

\begin{remark}
Note that in definition~\ref{def: pwl discrete inverse problem}, we only assume discrete observations of the response rather than a continuous trajectory, since they are sufficient for calibrating the path. 
\end{remark}



Next, we study how to solve the calibration problem and the corresponding properties of the solution space.



\subsection{Calibration}
\label{subsec: calibration}

Suppose that $\hat{X}^\delta = X^\delta({c}^\delta_y)$ is a piecewise linear path on ${\mathcal D}_\delta$, with ${c}^\delta_y = \{ c^\delta_{y,k}\}_{k=1}^{N_\delta}$ corresponding to the gradients of each linear segment, i.e.,
\begin{equation}
\label{X piecewise linear}
\hat{X}^\delta_t = \hat{X}^\delta_{(k-1)\delta} + c^\delta_{y,k}(t-(k-1)\delta),\ \ {\rm for}\ t\in[(k-1)\delta,k\delta],\ k=1,\dots,N_\delta.
\end{equation}
For $\hat{X}^\delta$ to be the solution to the discrete inverse problem for piecewise linear paths, it also needs to satisfy conditions (b) and (c) of definition \ref{def: pwl discrete inverse problem}. Using \eqref{X piecewise linear}, this becomes equivalent to studying the existence of ${c}_y^\delta$, such that $\hat{Y}^\delta_{k\delta} = y_{k\delta}$ for all $k=0,1,\dots,N_{\delta}$, where $\hat{Y}^\delta$ is the solution to 
\begin{equation}
\label{eq:local-all}
d\hat{Y}^\delta_t = \left( f(\hat{Y}^\delta_t) \cdot c^\delta_{y,k} \right) dt,\ t\in [(k-1)\delta, k\delta].
\end{equation}
This defines a system of $N_\delta$ independent equations of the form
\begin{equation}
\label{eq:local-generic}
d\tilde{Y}_t = \left( f(\tilde{Y}_t) \cdot c \right) dt, t\in [0, \delta],
\end{equation}
for $c=c^\delta_{y,k}$ and initial and terminal conditions $\tilde{Y}_0 = y_{(k-1)\delta}$ and $\tilde{Y}_\delta = y_{k\delta}$ respectively, for each $k=1,\dots,N_\delta$. Since $f\in{\rm Lip}(\gamma)$ with $\gamma>1$, the solution to \eqref{eq:local-generic} exists and is unique. We denote it by $F\left( t; y,c\right)$, where $y$ is the initial condition and $c$ the constant now incorporated into the new vector field $f\cdot c$. Then, for $t\in[(k-1)\delta, k\delta]$, 
\begin{equation}
\label{Y solution}
\hat{Y}^\delta_t = F(t-(k-1)\delta; y_{(k-1)\delta}, c^\delta_{y,k}).
\end{equation}
With $\hat{X}^\delta$ and $\hat{Y}^\delta$ expressed as \eqref{X piecewise linear} and \eqref{Y solution} respectively, conditions (a) and (b) in the definition \ref{def: pwl discrete inverse problem} are satisfied, while condition (c) is equivalent to 
\begin{equation}
\label{eq:c}
F(\delta; y_{(k-1)\delta}, c^\delta_{y,k}) = y_{k\delta}.
\end{equation}
Thus, constructing $\hat{X}^\delta=X^\delta({c}^\delta_y)$ is equivalent to solving the system of equations \eqref{eq:c} with respect to $c^\delta_{y,k}$ for $k=1,\dots,N_\delta$. 

Moreover, in order to understand the limiting behavior of the sequence $X^\delta({c}^\delta_y)$, we need to understand the space of solutions to this system. The following lemma links the rank of the Jacobian of the system to that of the vector field $f$, allowing us to infer the dimension of the solution space.

\begin{lemma}
\label{lemma: rank Z}
Consider the initial value problem
\begin{equation}
\label{eq: general ODE}
d\tilde{Y}_t = \left(f(\tilde{Y}_t) \cdot c \right) dt, \tilde{Y}_0 = y,
\end{equation}
where $f:\R^d\to L(\R^m,\R^d)$ is ${\rm Lip}(2)$, $c\in\R^m$ and $y\in\R^d$. Let $\tilde{Y}_t = \tilde{Y}(c)_t = F(t;y,c)$ be the solution to \eqref{eq: general ODE} and $\tilde{Z}(c)_t = G(t;y,c)$ be the derivative of the solution with respect to $c$, i.e.
\begin{equation}
\label{Z definition} 
G(t;y,c) = \bigtriangledown_c F(t;y,c) \in \R^{d\times m}.
\end{equation}
Then, as $t\downarrow 0$,
\begin{equation}
\label{eq: Jacobian rank}
\frac{1}{t}G(t;y,c) = f(y) + {\mathcal O}(t).
\end{equation}
In particular, if ${\rm rank}\left(f(y)\right) = \min(d,m)$ for all $y$, it follows that for every $y$ and $c$ there exists $t_0=t_0(y,c)>0$ such that for every $t\in(0,t_0]$,
\[ {\rm rank}\left( G(t;y,c)\right) = \min(d,m).\]
\end{lemma}

\begin{proof}
The vector field $f\cdot c$ of \eqref{eq: general ODE} is linear with respect to $c$ and ${\rm Lip}(2)$ with respect to $y$, and thus it is continuously differentiable. In particular, the solution map $(t,y,c)\mapsto F(t;y,c)$ exists and is continuously differentiable with respect to $c$ for every $y$ and $t$ (\cite{CoutinLejay} or \cite{kelley2010}, Theorem 8.49).

Moreover, the process $\tilde{Z}_t=\tilde{Z}(c)_t$ is the Jacobian of the solution with respect to $c$, that is $\tilde{Z}(c)_t=\bigtriangledown_cF(t;y,c)\in\R^{d\times m}$, and it satisfies
\begin{equation}
\label{Z equation} 
\frac{d}{dt} \tilde{Z}(c)_t  = A(\tilde{Y}(c)_t;c)\cdot \tilde{Z}(c)_t + f(\tilde{Y}(c)_t),
\end{equation}
where 
\begin{equation}
\label{def: A}
A(y;c) = \left(\bigtriangledown_y f(y)\right)\cdot c ,
\end{equation}
and $\bigtriangledown_y f \in\R^{d\times d\times m}$. The initial condition is
\[
\tilde{Z}(c)_0 = G(0;y,c) = \bigtriangledown_cF(0;y,c)=\bigtriangledown_c y \equiv {\bf 0}_{d\times m},
\]
where ${\bf 0}_{d\times m}$ is the $d\times m$-matrix with all entries equal to $0$.

Let $\Phi(t,s)$ denote the fundamental matrix of the homogeneous part of \eqref{Z equation}, that is the (invertible) solution of
\[
\frac{\partial}{\partial t}\Phi(t,s)=A(\tilde{Y}(c)_t;c)\Phi(t,s),~\Phi(s,s)=I_d.
\]
By the variation-of-constants formula applied columnwise to \eqref{Z equation}, it follows that
\begin{equation*}
\tilde{Z}(c)_t=\int_0^t \Phi(t,s)\,f(\tilde{Y}(c)_s)\,ds.
\end{equation*}

Fix $y$ and $c$. Since $f\in{\rm Lip}(2)$, $f$ is locally Lipschitz and $\tilde{Y}(c)_t$ is continuous in $t$, so there exists $t_1>0$ such that $\tilde{Y}(c)_{[0,t_1]}$ remains in a bounded neighborhood on which $f$ and $\bigtriangledown_y f$ are bounded. In particular, there exists $M<\infty$ such that
$\sup_{s,t\in[0,t_1]}\|\Phi(t,s)\|\leq e^{M t_1}$
and
$\sup_{s\in[0,t_1]}\|f(\tilde{Y}(c)_s)\|<\infty.$
Moreover, for $t\in(0,t_1]$ and $s\in[0,t]$ one has the bound
\[
\|\Phi(t,s)-I_d\|\leq |e^{M(t-s)}-1| \leq C\, (t-s)
\]
for a constant $C<\infty$ depending only on $M$ and $t_1$.

Using $\tilde{Z}(c)_t=\int_0^t \Phi(t,s)f(\tilde{Y}(c)_s) ds$, we write
\[
\tilde{Z}(c)_t - t f(y)
= \int_0^t \left(\Phi(t,s)-I_d\right) f(\tilde{Y}(c)_s) ds
+\int_0^t \left(f(\tilde{Y}(c)_s)-f(y)\right) ds.
\]
The first integral is ${\mathcal O}(t^2)$ by the estimate on $\Phi(t,s)-I_d$ and boundedness of $f(\tilde{Y}(c)_s)$ on $[0,t_1]$.
For the second integral, since $f$ is locally Lipschitz on the same neighborhood, there exists $L<\infty$ such that
\[
\|f(\tilde{Y}(c)_s)-f(y)\|\leq L \|\tilde{Y}(c)_s-y\|.
\]
Finally, from \eqref{eq: general ODE}, $\tilde{Y}(c)_s-y=\int_0^s f(\tilde{Y}(c)_u)\cdot c du$, hence $\|\tilde{Y}(c)_s-y\|={\mathcal O}(s)$ as $s\downarrow 0$, and therefore the second integral is also ${\mathcal O}(t^2)$.
This proves \eqref{eq: Jacobian rank}.

If ${\rm rank}(f(y))=\min(d,m)$, then there exists an $r\times r$ minor of $f(y)$, with $r=\min(d,m)$, whose determinant is non-zero. By \eqref{eq: Jacobian rank}, $\frac{1}{t}G(t;y,c)$ converges to $f(y)$ as $t\downarrow 0$, so by continuity of that determinant there exists $t_0=t_0(y,c)>0$ such that the same minor of $\frac{1}{t}G(t;y,c)$ remains non-zero for all $t\in(0,t_0]$. Hence ${\rm rank}(G(t;y,c))\geq r$ for $t\in(0,t_0]$, and since ${\rm rank}(G(t;y,c))\leq r$ always, it follows that ${\rm rank}(G(t;y,c))=r$ for $t\in(0,t_0]$.
\end{proof}



\begin{corollary}
\label{cor:from_lemma_rank_Z}
Consider the non-linear system of equations defined by \eqref{eq:c}, for fixed $k$, where $k\in\{1,2,\dots,N_\delta\}$ and $F$ is the solution to \eqref{eq:local-generic}. Fix any solution $c^\ast\in\R^m$ of \eqref{eq:c}, that is, $F(\delta;y_{(k-1)\delta},c^\ast)=y_{k\delta}$. Moreover, suppose that ${\rm rank}(f(y))=\min(d,m)$ for all $y$. Then the following local statements hold for $\delta\in(0,\delta_0]$ for some $\delta_0>0$ depending on $y_{(k-1)\delta}$ and $c^\ast$.  

\begin{itemize}
\item If $d<m$, then the system is under-determined in a neighborhood of $c^\ast$ and the space of solutions is locally a $C^1$ submanifold of $\R^m$ of dimension $m-d$. 
\item If $d=m$, then $c^\ast$ is locally unique and depends continuously on the terminal value $y_{k\delta}$, equivalently $c\mapsto F(\delta;y_{(k-1)\delta},c)$ is locally a $C^1$ diffeomorphism at $c^\ast$. 
\item If $d>m$, then $c^\ast$ is locally isolated. Moreover, for fixed $\delta>0$ and $y_{(k-1)\delta}$, the set of terminal values $y_{k\delta}\in\R^d$ for which \eqref{eq:c} admits a solution has $d$-dimensional Lebesgue measure zero, so for terminal data $y_{k\delta}$ there is no solution.
\end{itemize}
\end{corollary}

\begin{proof}
Fix $k$ and set $y_-:=y_{(k-1)\delta}$ and $y_+:=y_{k\delta}$. Define $\Phi_\delta:\R^m\to\R^d$ by $\Phi_\delta(c):=F(\delta;y_-,c)$. Then \eqref{eq:c} is equivalent to $\Phi_\delta(c)=y_+$. By assumption there exists $c^\ast\in\R^m$ with $\Phi_\delta(c^\ast)=y_+$.

By Lemma~\ref{lemma: rank Z}, for each $c^\ast$ there exists $\delta_0=\delta_0(y_-,c^\ast)>0$ such that for all $\delta\in(0,\delta_0]$,
\[
D\Phi_\delta(c^\ast)=\nabla_c F(\delta;y_-,c^\ast)=G(\delta;y_-,c^\ast)
\]
satisfies
\[
\frac{1}{\delta}D\Phi_\delta(c^\ast)=f(y_-)+{\mathcal O}(\delta),
\]
hence ${\rm rank}(D\Phi_\delta(c^\ast))=\min(d,m)$ for all $\delta\in(0,\delta_0]$, since ${\rm rank}(f(y_-))=\min(d,m)$.

Assume first that $d<m$. Then ${\rm rank}(D\Phi_\delta(c^\ast))=d$, so $D\Phi_\delta(c^\ast)$ is surjective. By the submersion theorem, the level set $\Phi_\delta^{-1}(y_+)$ is, in a neighborhood of $c^\ast$, a $C^1$ embedded submanifold of $\R^m$ of codimension $d$, hence of dimension $m-d$.

Assume next that $d=m$. Then ${\rm rank}(D\Phi_\delta(c^\ast))=d$, so $D\Phi_\delta(c^\ast)$ is invertible. By the inverse function theorem, $\Phi_\delta$ is a local $C^1$ diffeomorphism at $c^\ast$. In particular, in a neighborhood of $c^\ast$ the solution to $\Phi_\delta(c)=y_+$ is unique, and the inverse map $y_+\mapsto c$ is continuous.

Assume finally that $d>m$. Then ${\rm rank}(D\Phi_\delta(c^\ast))=m$, so $D\Phi_\delta(c^\ast)$ has full column rank. By the constant rank theorem, there exist neighborhoods $V$ of $c^\ast$ and $U$ of $y_+$, and $C^1$ coordinate changes on $V$ and $U$, such that in these coordinates $\Phi_\delta$ is given by $x\mapsto (x,0)\in\R^m\times\R^{d-m}$. It follows that $\Phi_\delta^{-1}(y_+)\cap V=\{c^\ast\}$, so $c^\ast$ is locally isolated. Moreover, the same normal form shows that $\Phi_\delta(V)$ is an $m$-dimensional embedded $C^1$ submanifold of $\R^d$ in a neighborhood of $y_+$. Since $d>m$, $\Phi_\delta(V)$ has $d$-dimensional Lebesgue measure zero, and therefore the set of terminal values $y_+\in\R^d$ for which $\Phi_\delta(c)=y_+$ has a solution is contained in $\Phi_\delta(\R^m)$ and is Lebesgue-null. This proves the non-existence claim.
\end{proof}

In the following remark, we discuss why both the over-determined ($m<d$) and under-determined ($m>d$) cases are problematic and we suggest ways of addressing the issues.

\begin{remark}
\label{rmk: uneven dimensions}
Let us assume that ${\rm rank}\left(f(y)\right) = \min(d,m)$ for all $y$. The case where ${\rm rank}\left(f(y)\right) < \min(d,m)$ can be dealt with in a similar manner.

\begin{itemize}
\item[(i)] Suppose $m>d$, i.e. there are more unknowns parameterizing the driving process $X$ than observations (under-determined case).  
We address this by reducing the dimension of the solution space, starting with the continuous inverse problem. Let $X = (X^s,X^u)$, where $X^s$ are the first $(m-d)$ coordinates (without loss of generality) and $X^u$ the remaining $d$ coordinates of $X$. Under the assumption that $X$ is a random rough path with known distribution, we sample $X^s$ from the corresponding marginal distribution of $X$ and set the gradients equal to the corresponding gradients of the sampled paths, denoted by ${c}^\delta_{s,x^s}$. Similarly to $X$, we write ${c}^\delta = ({c}^\delta_s,{c}^\delta_u)$. Then, system \eqref{eq:c} becomes
\[ F\left(\delta;y_{t_{i-1}},(c^\delta_{s,x^s},c^\delta_{u,y^u})\right) = y_{t_i}
\]
and we need to solve only for ${c}^\delta_{u,y^u}$. As the number of unknowns for each segment is now $d$, the system is no longer under-determined. In other words, we are now addressing the inverse problem for system
\[
\left( \begin{array}{c} dY_t \\ dX^u_t \end{array} \right) = \left( \begin{array}{cc} f_{1:d,1:m-d}(Y_t) &  f_{1:d,m-d+1:m}(Y_t)\\ 0 & 1 \end{array} \right) \left( \begin{array}{c} dX^s_t \\ dX^u_t \end{array} \right)
\]
where $f_{1:d,i:k}(y)$ is the projection of $f$ to a map acting only on the projection of $X$ to the corresponding $k-i$ dimensions. 
\item[(ii)] Suppose that $m<d$, i.e. there are fewer unknowns than equations, leading to non-existence of solutions (over-determined case). This issue can emerge even under the assumption that a solution $X$ to the continuous inverse problem exists and the core reason is that we have restricted the solution to the discrete inverse problem to too few degrees of freedom compared to the information available, when in fact the additional information can be used to achieve better accuracy in the approximation. Indeed, one possible way of addressing the issue is by considering piecewise linear paths on a (partially) finer partition than ${\mathcal D}_\delta$, so that degrees of freedom match the number of equations. 

By considering piecewise linear approximations on different time parameterizations of $X$, the inverse problem can be rephrased as an inverse problem of a system where $m=d$. To demonstrate this point, we consider the case where $m=1$ and $d=2$. Let $X^\delta({c}^\delta_x)$ and $X^\delta(\tilde{c}^\delta_x)$ be piecewise linear interpolations of $X$ on a partition of size $\frac{\delta}{2}$, 
defined as
\[
X^\delta({c}^\delta_x)_t = \left\{ \begin{array}{cc}
 X_{(k-1)\delta}+2\frac{X_{k\delta}-X_{(k-1)\delta}}{\delta}(t-(k-1)\delta), & t \in [(k-1)\delta, (k-\frac{1}{2})\delta]\\
   X_{k\delta},  &t \in [(k-\frac{1}{2})\delta,k\delta]
\end{array}
\right.
\]
and 
\[
X^\delta(\tilde{c}^\delta_x) = \left\{ \begin{array}{cc}
 X_{(k-1)\delta} ,   & t \in [(k-1)\delta, (k-\frac{1}{2})\delta]\\
   X_{(k-1)\delta} + 2\frac{X_{k\delta}-X_{(k-1)\delta}}{\delta}(t-(k-\frac{1}{2})\delta),  &t \in [(k-\frac{1}{2})\delta,k\delta]
\end{array}
\right.,
\]
for all $k=1,\dots,N_\delta$.
We can define the solution to the discrete inverse problem as the solution to the following system:
\[
\left( \begin{array}{c} d\hat{Y}^{(1,\delta)}_t \\ d\hat{Y}^{(2,\delta)}_t \end{array} \right) = \left( \begin{array}{cc} f(Y_t) &  0\\ 0 & f(Y_t) \end{array} \right) \left( \begin{array}{c} dX^\delta({c}^\delta_x)_t \\ dX^\delta(\tilde{c}^\delta_x)_t \end{array} \right).
\]
Note that this is different than \eqref{eq:approx2} but they both converge to \eqref{eq:main} in $p$-variation, as it is possible to show that $(X^\delta({c}^\delta_x),X^\delta(\tilde{c}^\delta_x))\to ({\bf X},{\bf X})$ in $p$-variation.

An alternative approach to increasing the dimension of the driving path approximation is to consider approximation by Abelian rough paths \cite{Flint}.
\end{itemize}
\end{remark}

For the remaining paper, we will assume that $m=d$ and $f(y)$ has rank $d$ for all $y\in\R^m$, for simplicity of the exposition. Note that, in that case, neither uniqueness nor existence is guaranteed since the system is non-linear, but we can have at most countably many solutions. 


\subsection{Convergence of solutions to the discrete inverse problem for piecewise linear paths}
\label{subsec: convergence piecewise linear}

We now consider the limiting behavior of the solution $\hat{X}^\delta = X^\delta({c}^\delta_y)$ to the discrete inverse problem as $\delta\to 0$, under the assumption that $X({c}^\delta)$ is the family of piecewise linear paths, with parameter ${c}^\delta$ corresponding to the gradients of the linear segments. It follows from theorem \ref{thm: general convergence} that \eqref{eq: convergence of X(cy)} will hold provided that the assumptions are satisfied. Below, we give simple-to-check conditions for assumptions of lemma \ref{lemma: discrete solution Y convergence} to hold, involving only properties of the vector field $f$. 
We build upon the following lemma, which is a simple corollary of the classical inverse function theorem. 

\begin{lemma}
\label{lemma: inverse funciton theorem}
    Let $F:\R^d\to \R^d$ be a twice differentiable function and choose $c,c_0\in \R^d$ such that $|\bigtriangledown F(c_0)|\neq 0$ and $F(c)$ is sufficiently close to $F(c_0)$, in the sense that for a neighborhood $V = {\mathcal N}_{\rho}(c_0) = \{x\in\R^d:\ |x-c_0|<\rho\}$ for some $\rho>0$ with ${\inf_{x\in{\mathcal N}_\rho(c_0)}|\bigtriangledown F(x)| = \kappa > 0}$, we have $F(c)\in F(V)$. Then, there exists a $\tilde{c}\in V$ such that 
    \begin{equation}
        | \tilde{c}-c_0 | \leq C | F(c) - F(c_0) |,
    \end{equation}
    for some $C \leq \frac{1}{\kappa}<\infty$.
\end{lemma}
\begin{proof}
    Let $W=F(V)$. It follows from the multivariate inverse function theorem that there exists a $\rho>0$ such that $F$ has a local continuous inverse $F^{-1}_{c_0,\rho}:W\to V$ and $\inf_{x\in{\mathcal N}_\rho(c_0)}|\bigtriangledown F(x)| = \kappa > 0$. Moreover, the derivative of the inverse is equal to 
    \[ \bigtriangledown F^{-1}_{c_0,\rho}(y) = \left( \bigtriangledown F (F^{-1}_{c_0,\rho}(y) \right)^{-1},\]
    for every $y\in W$. Since $c$ is such that $F(c) \in W$, there exists a $\tilde{c}\in V$, such that $F(\tilde{c}) = F(c)$. Using the properties of $F^{-1}_{c_0,\rho}$ and applying Taylor's theorem, we write
    \begin{eqnarray*}
        \tilde{c} - c_0 = F^{-1}_{c_0,\rho}\left( F(c) \right) - F^{-1}_{c_0,\rho}\left( F(c_0) \right) &=& \bigtriangledown F^{-1}_{c_0,\rho}(\xi) \cdot \left( F(\tilde{c}) - F(c_0) \right) \\
        &=& \left( \bigtriangledown F (\zeta) \right)^{-1}\cdot \left( F(\tilde{c}) - F(c_0) \right),
    \end{eqnarray*}
    for some $\xi\in W$ and corresponding $\zeta\in V$. It follows that 
    \begin{eqnarray*}
    |\tilde{c} - c_0 | &\leq & |\left( \bigtriangledown F (\zeta) \right)^{-1}|\cdot |\left( F(\tilde{c}) - F(c_0) \right)| \leq  \sup_{x\in V}|\left( \bigtriangledown F (x) \right)^{-1}|\cdot |\left( F(\tilde{c}) - F(c_0) \right)| \\
    &\leq & \frac{1}{\inf_{x\in V}|\left( \bigtriangledown F (x \right)|}|\left( F({c}) - F(c_0) \right)| = \frac{1}{\kappa}|\left( F({c}) - F(c_0) \right)|.
    \end{eqnarray*}
\end{proof}


In the following proposition, we give conditions on $f$ so that assumption \eqref{eq: in-segment smoothness} of lemma \ref{lemma: discrete solution Y convergence} holds.

\begin{proposition}
\label{prop: pws in-segment smoothness}
Let ${y}$ be the trajectory of $Y = I({\bf X})$, where $I$ is the It\^o-Lyons map corresponding to \eqref{eq:main}, $m=d$ and ${\bf X}\in G\Omega_p(\R^d)$ is the limit of nested dyadic piecewise linear interpolations in $p$-variation. 
Let ${c}^\delta_y$ be a solution to the corresponding discrete inverse problem. Then, if $f$ satisfies 
\begin{equation}
    \label{eq: f assumptions}
    f\in{\rm Lip}(2),\ \ \inf_{t\in[0,T]}|f(y_t)|=\kappa >0,\ {\rm and}\ \sup_{t\in[0,T]}|\bigtriangledown f(y_t)|=K<+\infty,
\end{equation}
and for $\delta$ sufficiently small, $\exists~\tilde{\bf X}\in {\mathcal I}^{-1}({y})$, such that $X^\delta({c}^\delta_y)$ and $X^\delta({c}^\delta_{\tilde{x}})$ satisfy \eqref{eq: in-segment smoothness}.
\end{proposition}

\begin{proof}
By construction, ${\bf X}\in {\mathcal I}^{-1}({y})$ and let $X^\delta({c}^\delta_x)$ be the corresponding piecewise linear interpolation. Let ${c}^\delta_y$ be a solution to the discrete inverse problem. We assume that such a solution exists. We want to show that there exists an $\tilde{\bf X}\in{\mathcal I}^{-1}({y})$ (possibly different than ${\bf X}$) such that \eqref{eq: in-segment smoothness} holds. 

Consider a segment $[t_{i-1},t_{i}]$, where $t_{i}-t_{i-1}=\delta,\ \forall i\geq 1$. By definition, $X^\delta({c}^\delta)$ is linear in that segment, with gradient $c^\delta_{i}$. It follows that for any ${c}^\delta$ and $\delta$ small enough,
\begin{equation}
\label{eq: in-segment p-var bound}
    \| {\bf X}^\delta({c}^\delta_{y})- {\bf X}^\delta({c}^\delta) \|_{p-{\rm var},[t_{i-1},t_{i}]} = 
    |c^\delta_{y,i} - c^\delta_{i}| \delta.
\end{equation}
Let $F_{\delta,i}:\R^d\to \R^d$ be defined as $F_{\delta,i}(c) = F(\delta; y_{t_{i-1}},c)$ where $F(t;y,c)$ is the solution to \eqref{eq:local-generic}, as before. Let ${c}^\delta_{y,i}$ and ${c}^\delta_{x,i}$ play the role of $c_0$ and $c$ in lemma \ref{lemma: inverse funciton theorem}, respectively. Then, it follows from lemma \ref{lemma: inverse funciton theorem} that for ${c}^\delta_{x,i}$ such that $F_{\delta,i}({c}^\delta_{x,i})$ and $F_{\delta,i}({c}^\delta_{y,i})$ are sufficiently close, there exists a $\tilde{c}^\delta_i$ such that $F_{\delta,i}(\tilde{c}^\delta_i) = F_{\delta,i}({c}^\delta_{x,i})$ and
\begin{equation}
\label{eq: Dc bounded by DF}
|\tilde{c}^\delta_i - {c}^\delta_{y,i}| \leq C |F_{\delta,i}({c}^\delta_{x,i})-F_{\delta,i}({c}^\delta_{y,i})|,
\end{equation}
for some $C>0$. By definition, $F_{\delta,i}({c}^\delta_{y,i})=y_{t_i}$. Also, $F_{\delta,i}({c}^\delta_{x,i})=\tilde{y}^{t_{i-1},\delta}_{t_i}$, where 
$\tilde{Y}^{u,\delta}$ is the response to $X^\delta({c}^\delta_{x})$ with initial conditions $\tilde{Y}^{u,\delta}_{u} = y_{u}$ and trajectory $\tilde{y}^{u,\delta} = \{\tilde{y}^{u,\delta}_t;
 t\geq u\}$. Following similar arguments as in lemma \ref{lemma: discrete solution Y convergence}, we can show that there exists a $C_1>0$ and $\eta_1(\delta)$ with $\eta_1(\delta)\to 0$, such that
\begin{equation}
    \label{eq: difference of ys bound}
    |y_{t_i} - \tilde{y}^{t_{i-1},\delta}_{t_i}| \leq C_1 \eta_1(\delta) \omega(t_{i-1},t_i)^\frac{1}{p}.
\end{equation}
Thus, we can choose $\delta>0$ sufficiently small so that \eqref{eq: Dc bounded by DF} holds. Lemma \ref{lemma: inverse funciton theorem} also states that constant $C$ in \eqref{eq: Dc bounded by DF} is bounded by 
$C\leq \sup_{c\in V}\left( |\bigtriangledown F_{\delta,i}(c)^{-1}|\right)$, 
where we can choose $V$ to be a neighborhood of $c^\delta_{y,i}$ with radius $\epsilon(\delta)$, such that its image contains $F_{\delta,i}(c^\delta_{x,i})$. Since $F_{\delta,i}(c^\delta_{x,i})\to F_{\delta,i}(c^\delta_{y,i})=y_{t_i}$ as $\delta\to 0$, we can choose $\epsilon(\delta)\to 0$, as $\delta\to 0$.
Since $f\in{\rm Lip}(2)$, the solution map $(t,y,c)\mapsto F(t;y,c)$ is $C^1$ in $c$, and
\[
\bigtriangledown F_{\delta, i}(c)=
\bigtriangledown_c F(\delta; y_{t_{i-1}},c) = G(\delta;y_{t_{i-1}}, c).
\]
Moreover, for $c\in V$ it holds that $|c|\delta \le |c^\delta_{y,i}| \delta + \delta \epsilon(\delta)$, hence the trajectories $u \mapsto F(u; y_{t_{i-1}}, c)$, $u \in [0, \delta]$, remain in a small neighborhood of $y_{t_{i-1}}$ uniformly over $c \in V$ when $\delta$ is small. Thus, the estimate in lemma~\ref{lemma: rank Z} applies uniformly on $V$ and yields
\[
\frac{1}{\delta} \bigtriangledown F_{\delta, i}(c)
= \frac{1}{\delta} G(\delta; y_{t_{i-1}}, c)
= f(y_{t_{i-1}}) + {\mathcal O}(\delta),~c\in V.
\]
Using \eqref{eq: f assumptions}, we may choose $\delta>0$ small enough so that
\[
\inf_{c\in V} \left|\frac{1}{\delta} \bigtriangledown F_{\delta, i}(c)\right| \ge \frac{\kappa}{2},
\]
hence
\[
\sup_{c\in V}\left|\bigtriangledown F_{\delta, i}(c)^{-1}\right|
\le \frac{2}{\kappa} \delta^{-1},
\]
and, consequently,
\begin{equation}
\label{eq: segment p-var bound}
\begin{split}
\| {\bf X}^\delta({c}^\delta_{y})- {\bf X}^\delta(\tilde{c}^\delta) \|_{p-{\rm var},[t_{i-1},t_{i}]}
& \leq \frac{2C_1}{\kappa} \eta_1(\delta) \omega(t_{i-1}, t_i)^\frac{1}{p} \\
 &\leq C_2 \eta_2(\delta)\omega(t_{i-1},t_i)^\frac{1}{p},
 \end{split}
\end{equation}
for $\eta_2(\delta)\to 0$ as $\delta\to 0$, where $X^\delta(\tilde{c}^\delta)$ is the piecewise linear path corresponding to local gradients $\tilde{c}^\delta_i$ for $i=1,\dots,N_\delta$. Note that the bound in \eqref{eq: segment p-var bound} is independent of the choice of ${c}^\delta_y$ and only depends on the interval through $\omega(t_{i-1},t_i)^\frac{1}{p}$. 

It remains to show that there exists a $\tilde{\bf X}\in{\mathcal I}^{-1}({y})$ such that 
\begin{equation}
    \label{eq: from tilde c to tilde x}
    \| {\bf X}^\delta({c}^\delta_{\tilde{x}})- {\bf X}^\delta(\tilde{c}^\delta) \|_{p-{\rm var},{[t_{i-1},t_i]}} \leq C \eta (\delta)\omega(t_{i-1},t_i)^\frac{1}{p},
\end{equation}
for some constant $C<\infty$ and $\eta(\delta)\to 0$ as $\delta\to 0$. Without loss of generality, we take $[t_{i-1},t_i] = [0,\delta]$. We will construct a sequence $\{\tilde{c}^{2^{-n} \delta}_i\}_{i=1}^{2^n}$ for $n\geq 0$ and $\delta$ fixed, such that ${\bf X}^\delta(\tilde{c}^{2^{-n} \delta})$ is Cauchy on $[0,\delta]$ with respect to the $p$-variation topology and thus converges in $p$-variation as $n\to\infty$ to a rough path $\tilde{\bf X}$ on $[0,\delta]$. Moreover, we will show that \eqref{eq: from tilde c to tilde x} holds and $\tilde{\bf X}$ is a solution to the continuous inverse problem on $[0,\delta]$. The concatenation of such limits on all intervals $[t_{i-1},t_i]$ will be a solution to the continuous inverse problem on $[0,T]$.

For $n=0$, $\tilde{c}^\delta_1$ is constructed as above, i.e. it is the solution of 
\[
F(\delta; y_0,\tilde{c}^\delta_1) = F(\delta; y_0,c^\delta_{x,1}) = y^\delta_\delta
\]
that lies within the neighborhood $V$ of $c^\delta_{y}$ where $F(\delta; y_0,\cdot)$ is invertible. 
For $n=1$, $\tilde{c}^{\delta/2}_1$ and $\tilde{c}^{\delta/2}_2$ are defined as solutions to 
\begin{eqnarray*}
 F\left(\frac{\delta}{2}; y_0, \tilde{c}^{\delta/2}_1\right) &=& F\left(\frac{\delta}{2}; y_0, c^{\delta/2}_{x,1}\right) = y^{\delta/2}_{\frac{\delta}{2}},   \\
 F\left(\frac{\delta}{2}; y^{\delta/2}_{\frac{\delta}{2}}, \tilde{c}^{\delta/2}_2\right) &=& F\left(\frac{\delta}{2}; y^{\delta/2}_{\frac{\delta}{2}}, c^{\delta/2}_{x,2}\right) = y^{\delta/2}_{\delta},
\end{eqnarray*}
respectively. We have chosen $\tilde{c}^\delta_1$ so that $F(\delta; y_0,\tilde{c}^\delta_1) = F(\delta; y_0, c^\delta_{x,1})$. When viewed as a function of $t$, however, $F$ is invariant with respect to time parameterization since it corresponds to the increment of the response and thus, it follows that $F(t; y_0,\tilde{c}^\delta_1) = F(t; y_0, c^\delta_{x,1})$, $\forall t\geq 0$. Thus,
\[
F\left(\frac{\delta}{2}; y_0,\tilde{c}^\delta_1\right) = F\left(\frac{\delta}{2}; y_0,{c}^\delta_{x,1}\right) \approx F\left(\frac{\delta}{2}; y_0, c^{\delta/2}_{x,1}\right) = F\left(\frac{\delta}{2}; y_0, \tilde{c}^{\delta/2}_1\right),
\]
or, more precisely, using the Universal Limit Theorem
we get
\begin{equation*}
\begin{split}
\left|F\left(\frac{\delta}{2}; y_0, \tilde{c}^{\delta/2}_1\right) - F\left(\frac{\delta}{2}; y_0,\tilde{c}^\delta_1\right)\right|
&=  \left| F\left(\frac{\delta}{2}; y_0, c^{\delta/2}_{x,1}\right) - F\left(\frac{\delta}{2}; y_0,{c}^\delta_{x,1}\right)\right|\\
&\leq C_3 \|{\bf X}^{\delta/2}({c}^{\delta/2}_x) - {\bf X}^\delta({c}^\delta_x)\|_{p-{\rm var},[0,\frac{\delta}{2}]}.
\end{split}
\end{equation*}
Since ${\bf X}^\delta({c}^\delta_x)$ converges in $p$-variation, it follows that there exist an $\eta_3(\delta)\to 0$ and a constant $C_4$, such that $\|{\bf X}^{\delta/2}({c}^{\delta/2}_x) - {\bf X}^\delta({c}^\delta_x)\|_{p-{\rm var},[0,\frac{\delta}{2}]} \leq C_4 \eta_3(\delta) \omega(0,\frac{\delta}{2})^\frac{1}{p}$. It follows from lemma \ref{lemma: inverse funciton theorem} applied to the function $F_{\frac{\delta}{2},1}(\cdot) = F(\frac{\delta}{2}; y_0, \cdot)$ that, for $\delta$ sufficiently small, there exists a neighborhood $\tilde{V}$ of $\tilde{c}^\delta$ that contains a solution $\tilde{c}^{\delta/2}_1$, such that
\[
\tilde{c}^\delta_1 - \tilde{c}^{\delta/2}_1 = \bigtriangledown F_{\frac{\delta}{2},1}(\zeta_1)^{-1}\cdot \left( F\left(\frac{\delta}{2}; y_0,\tilde{c}^\delta_1\right) - F\left(\frac{\delta}{2}; y_0, \tilde{c}^{\delta/2}_1\right)\right),
\]
for some $\zeta_1\in \tilde{V}$. Similarly,
\begin{equation*}
\begin{split}
&F\left(\frac{\delta}{2}; y^{\delta/2}_{\frac{\delta}{2}},\tilde{c}^\delta_1\right) \approx F\left(\frac{\delta}{2}; y^{\delta}_{\frac{\delta}{2}},\tilde{c}^\delta_1\right) = F\left(\frac{\delta}{2}; y^\delta_{\frac{\delta}{2}},{c}^\delta_{x,1}\right) = y^\delta_\delta\\
&\approx y^{\delta/2}_{\delta}=F\left(\frac{\delta}{2}; y^{\delta/2}_{\frac{\delta}{2}}, c^{\delta/2}_{x,2}\right) = F\left(\frac{\delta}{2}; y^{\delta/2}_{\frac{\delta}{2}}, \tilde{c}^{\delta/2}_2\right),
\end{split}
\end{equation*}
where the first approximation comes from the difference in initial conditions and the second from the change in drivers. Using \cite[Lemma 2.7 and Theorem 2.9]{Baudoin}, we get
\begin{equation*}
\begin{split}
    \left|F\left(\frac{\delta}{2}; y^{\delta/2}_{\frac{\delta}{2}}, \tilde{c}^{\delta/2}_2\right) - F\left(\frac{\delta}{2}; y^{\delta/2}_{\frac{\delta}{2}},\tilde{c}^\delta_1\right) \right| &\leq 
     \left|F\left(\frac{\delta}{2}; y^{\delta/2}_{\frac{\delta}{2}}, \tilde{c}^{\delta/2}_2\right) - F\left(\frac{\delta}{2}; y^{\delta}_{\frac{\delta}{2}},\tilde{c}^\delta_1\right) \right| \\
     & + \left|F\left(\frac{\delta}{2}; y^{\delta}_{\frac{\delta}{2}}, \tilde{c}^{\delta}_1\right) - F\left(\frac{\delta}{2}; y^{\delta/2}_{\frac{\delta}{2}},\tilde{c}^\delta_1\right) \right| \\
      &\leq |y^{\delta/2}_\delta - y^\delta_\delta| + C_5 \cdot|y^\delta_{\frac{\delta}{2}} - y^{\delta/2}_{\frac{\delta}{2}}| \cdot \omega(\frac{\delta}{2},\delta)^\frac{1}{p} \\
      &\leq C_6 \|{\bf X}^{\delta/2}({c}^{\delta/2}_x) - {\bf X}^\delta({c}^\delta_x)\|_{p-{\rm var},[0,\delta]}.
\end{split}
\end{equation*}
So, for $\delta$ sufficiently small, there also exists a solution $\tilde{c}^{\delta/2}_2\in\tilde{V}$ satisfying
\[
\tilde{c}^\delta_1 - \tilde{c}^{\delta/2}_2 = \bigtriangledown F_{\frac{\delta}{2},y^\delta_{\frac{\delta}{2}}}(\zeta_2)^{-1}\cdot \left( F\left(\frac{\delta}{2}; y^{\delta/2}_{\frac{\delta}{2}},\tilde{c}^\delta_1\right) - F\left(\frac{\delta}{2}; y^{\delta/2}_{\frac{\delta}{2}}, \tilde{c}^{\delta/2}_2\right)\right),
\]
for some $\zeta_2\in\tilde{V}$. It follows that
\[
\| {\bf X}^\delta(\tilde{c}^\delta) - {\bf X}^{\delta/2}(\tilde{c}^{\delta/2}) \|_{p-{\rm var},[0,\delta]} \leq C_7  \| {\bf X}^\delta({c}^\delta_x) - {\bf X}^{\delta/2}({c}^{\delta/2}_x) \|_{p-{\rm var},[0,\delta]}.
\]
Following a similar construction for $n\geq 1$, this generalises to 
\begin{equation*}
\begin{split}
&\| {\bf X}^{2^{-n} \delta}(\tilde{c}^{2^{-n} \delta}) - {\bf X}^{2^{-(n+1)} \delta}(\tilde{c}^{2^{-(n+1)}\delta}) \|_{p-{\rm var},[0,\delta]} \\
&\leq C_8 \|{\bf X}^{2^{-n} \delta}({c}^{2^{-n} \delta}_x) - {\bf X}^{2^{-(n+1)}\delta}({c}^{2^{-(n+1)}\delta}_x)\|_{p-{\rm var},[0,\delta]},
\end{split}
\end{equation*}
which goes to $0$ with $n\to\infty$, since we know that ${\bf X}^{2^{-n}\delta}({c}^{2^{-n} \delta}_x)$ converges in $p$-variation as $n\to\infty$ and for $\delta$ fixed. Thus, ${\bf X}^{2^{-n}\delta}(\tilde{c}^{2^{-n} \delta})$ also converges in $p$-variation on $[0,\delta]$. Moreover, its limit $\tilde{\bf X}$ will be a solution to the continuous inverse problem since, by construction, $I(X(\tilde{c}^{2^{-n} \delta}))_{2^{-n}k\delta} = I(X({c}^{2^{-n} \delta}_x))_{2^{-n}k\delta}$, for every $k=0,1,\dots,2^n$. 

Finally, to show \eqref{eq: from tilde c to tilde x}, we note that 
\[
\|{\bf X}^\delta(\tilde{c}^\delta) - {\bf X}^\delta({c}^\delta_{\tilde{x}})\|_{p-{\rm var},[0,\delta]} = |\tilde{c}^\delta_1 - c^\delta_{\tilde{x},1}|\delta = 
\lim_{n\to\infty} |\tilde{c}^\delta_1 - 2^{-n}\sum_{k=1}^{2^n} \tilde{c}^{2^{-n}\delta}_{k}|\delta
\]
and following the same construction as above, we get 
\begin{eqnarray*}
   \|{\bf X}^\delta(\tilde{c}^\delta) - {\bf X}^\delta({c}^\delta_{\tilde{x}})\|_{p-{\rm var},[0,\delta]} &\leq& C_9 \lim_{n\to\infty}\|\hat{Y}^\delta(x) - \hat{Y}^{2^{-n}\delta}(x)\|_{p-{\rm var},[0,\delta]} \\
   &=& C_9 \|\hat{Y}^\delta(x) - Y\|_{p-{\rm var},[0,\delta]} \leq C \eta(\delta)\omega(0,\delta)^\frac{1}{p},
\end{eqnarray*}
where $\hat{Y}^\delta(x) = I(X^\delta({c}^\delta_x))$, with trajectory $\{y^\delta_t; t\geq 0\}$.

\end{proof}

We can now state the main theorem for this section. It is a direct consequence of theorem \ref{thm: general convergence} and proposition \ref{prop: pws in-segment smoothness}. 

\begin{theorem}
\label{thm: pwl convergece}
Let $Y=I({\bf X})$ be a solution to \eqref{eq:main} with trajectory ${y}:[0,T]\to \R^d$, with $I$ defined on the closure of the space of nested dyadic piecewise linear paths in $p$-variation. Let $X^{\delta}({c}^\delta)$ be the family of piecewise linear paths in dyadic partition ${\mathcal D}_\delta$ and assume that there exists a solution $X^{\delta}({c}^\delta_y)$ to the discrete inverse problem. Suppose that $m=d$ and $f$ satisfies assumptions \eqref{eq: f assumptions}. Then, \eqref{eq: convergence of X(cy)} holds, provided that $I$ has a weakly continuous inverse.   
\end{theorem}

\begin{remark}
In theorem \ref{thm: pwl convergece}, we assume that the $I$ has a weakly continuous inverse. Ideally, this is something that we would like to prove, under appropriate conditions on the system. While corollary~\ref{cor:from_lemma_rank_Z} can be interpreted as a result on the local stability of the inverse of the It\^o map restricted to piecewise linear drivers on one time step, extending to the whole path space equipped with the $p$-variation topology is highly non-trivial. 
\end{remark}

\color{black}


\section{Numerical Algorithms}
\label{sec: numerical}

For the remainder of the paper, we focus on the case where the nested dyadic piecewise linear interpolations of $X$ driving \eqref{eq:main} converge and we assume that $T$ is fixed, $m=d$ and $f$ satisfies assumptions \eqref{eq: f assumptions}. Then, we know from theorem \ref{thm: pwl convergece} that the solution to the discrete inverse problem for $\delta$ sufficiently small will be close to a solution of the continuous inverse problem. 

When solutions ${c}^\delta$ to \eqref{eq:c} have an exact closed-form expression, we are able to construct solutions to the discrete inverse problem exactly. However, in most cases, an analytic solution does not exist, and we need to resort to numerical approximations. Note that, while connected to solving \eqref{eq:c}, our goal is to construct a piecewise linear path $\hat{X}^\delta$ corresponding to observations ${y}_{{\mathcal D}_{\delta}}$, for $N_{\delta}\delta = T$, that solves the discrete inverse problem of definition \ref{def: pwl discrete inverse problem}. In this section, we discuss iterative algorithms for constructing numerical approximations $\hat{X}^\delta(n)$ of $\hat{X}^\delta$. We denote by ${c}^\delta(n) = \left(c^\delta(n)_k \right)_{k=0}^{N_\delta}$ the gradients of $\hat{X}^\delta(n)$. We define the total error of the numerical approximation as the $p$-variation distance between $\hat{X}^\delta(n)$ and the set of all possible solutions to the discrete inverse problem, denoted by ${\mathcal I}^{-1}_\delta({y})$:
\begin{equation}
\label{eq:error}
e_p(n;\delta) = d_p\left( \hat{X}^\delta(n), {\mathcal I}^{-1}_\delta({y}) \right).
\end{equation}
As the ultimate goal is to construct an approximation to a solution of the continuous inverse problem, numerical algorithms need to converge uniformly with respect to $\delta$, i.e.
\[ \lim_{n\to\infty}\sup_{\delta>0} e_p(n;\delta) = 0,
\]
for $\delta$ in a sufficiently small neighborhood of $0$. Then, if assumptions of theorem \ref{thm: pwl convergece} are satisfied and the space of solutions ${\mathcal I}^{-1}_\delta({y})$ is compact, we can control the total error $d_p\left(\hat{X}^\delta(n), {\mathcal I}^{-1}({y})\right)$.

\begin{remark}
\label{rem: L1 error bound}
Note that, as is standard in the context of rough paths, we define the $p$-variation norm on the rough path, i.e. for any $p$-rough path ${\bf X}$,
\[
\|{\bf X} \|_{p-{\rm var},[0,T]} = \sup_{k=1,\dots,\lfloor p \rfloor} \sup_{{\mathcal D}}\left( \sum_{t_i\in {\mathcal D}} |{\bf X}^k_{t_i,t_{i+1}}|^\frac{p}{k}\right)^\frac{k}{p},
\]
where $|\cdot |$ are the Euclidean norms in the corresponding Euclidean space and ${\mathcal D}$ are finite partitions of $[0,T]$. In practice, this is very expensive to compute. Moreover, ${\mathcal I}^{-1}_\delta({y})$ will not be known in advance. When testing the performance of different algorithms in section \ref{sec: examples}, we use the $L_1$ norm between observed data and corresponding response to the system, when driven by the numerical solution to the discrete inverse problem,
i.e. we use the distance
\begin{equation}
\label{eq:error2}
    \|y_{{\mathcal D}_\delta} - \hat{Y}^\delta(n)_{{\mathcal D}_\delta}\|_1 = \sum_{i=1}^{N_\delta} |y_{t_i} - \hat{Y}^\delta(n)_{t_i}|,
\end{equation}
where $\hat{Y}^\delta(n) = I\left(\hat{X}^\delta(n) \right)$. We claim that, for fixed $\delta$, this provides an upper bound to the $p$-variation norm of the difference of the piecewise linear paths times a constant that depends on $\delta$. While the constant will explode as $\delta\to0$ when the underlying path $X$ is not of bounded variation, we will show that the convergence of the difference is fast enough to compensate for this, leading to convergence in $p$-variation with respect to $n$, uniformly in $\delta$. 

The claim is supported by the following inequalities that can be derived by taking advantage of the piecewise linear assumption for $X^\delta({c}^\delta)$ and $X^\delta({c}^\delta(n))$, which guarantees that the supremum over all partitions will be achieved at a subset of ${\mathcal D}_\delta$ (see \cite{Driver}). First, it is straight-forward to show that
\[ 
\| \hat{X}^\delta - \hat{X}^\delta(n) \|_{1-{\rm var},[0,T]} \leq \delta \sqrt{d} \sum_{i=1}^{N_\delta} |c^\delta_i-c^\delta(n)_i|,
\]
using standard inequalities between $L_1$ and $L_2$ norms applied to $|\cdot |$.
Using the piecewise linear assumption for $X^\delta({c}^\delta)$ and $X^\delta({c}^\delta(n))$ and Chen's identity, it is possible to get the following bound
\begin{equation}
\label{eq: p-to-1 bound}
    \| \hat{\bf X}^\delta - \hat{\bf X}^\delta(n) \|_{p-{\rm var},[0,T]} \leq C_1 \cdot C_2(\hat{X}^\delta,\hat{X}^\delta(n))\cdot 
\delta \sum_{i=1}^{N_\delta} |c^\delta_i-c^\delta(n)_i|,
\end{equation}
for constant $C_1$ depending only on $p$ and dimension $d$ and $C_2(\hat{X}^\delta,\hat{X}^\delta(n))$ bounded by
\begin{equation}
\label{eq: p-to-1 bound constant}
    C_2(\hat{X}^\delta,\hat{X}^\delta(n)) \leq \max_{k=1,\dots,\lfloor p \rfloor} \max_{\ell=0,\dots,k-1} \left( \|\hat{X}^\delta\|_{1-{\rm var},[0,T]}^\ell\|\hat{X}^\delta(n)\|^{k-1-\ell}_{1-{\rm var},[0,T]} \right).
\end{equation}
Moreover, under the assumptions of proposition \ref{prop: pws in-segment smoothness} and using similar arguments based on the inverse function theorem for $F(\delta;Y_{t_i},\cdot)$ and the initial value correction, we get
\[
|c^\delta_i-c^\delta(n)_i| \leq C_3 |y_{t_i}-\hat{Y}^\delta(n)_{t_i}| + C_4 |y_{t_{i-1}}-\hat{Y}^\delta(n)_{t_{i-1}}|.
\]

\end{remark}

Below, we discuss two different numerical approaches for constructing the solution to the discrete inverse problem. The first is based on classical methods, solving for one gradient at a time, while the second takes a global approach allowing solutions for each segment in the partition ${\mathcal D}_\delta$ to also inform those in neighboring segments. We will compare the performance of the two approaches in more detail in section \ref{sec: examples}.


\subsection{Local approach: the Newton-Raphson scheme}
\label{sec: classic}

The first approach is based on using classical numerical methods to solve \eqref{eq:c}. In fact, constructing a numerical solution for \eqref{eq:c} for fixed $k$ is a standard problem that involves numerically solving a system of equations and an ODE. For example, applying the Newton-Raphson algorithm in the context of \eqref{eq:c} gives
\begin{equation}
\label{eq:NR}
\begin{split}
c^\delta(n+1)_k &= c^\delta(n)_k + D_c F(\delta; y_{(k-1)\delta},c^\delta(n)_k)^{-1} \left( y_{k\delta} - F(\delta; y_{(k-1)\delta},c^\delta(n)_k) \right) \\
&= c^\delta(n)_k + G(\delta; y_{(k-1)\delta},c^\delta(n)_k)^{-1} \left( y_{k\delta} - F(\delta; Y_{(k-1)\delta},c^\delta(n)_k) \right),
\end{split}
\end{equation}
Both $F(\delta; y_{(k-1)\delta},c^\delta(n)_k)$ and $G(\delta; y_{(k-1)\delta},c^\delta(n)_k)$ are solutions to the system of ODEs described in \eqref{eq: general ODE} and \eqref{Z equation} respectively, which can be solved using standard numerical techniques, such as the Euler method with step $\delta^\prime << \delta$. Note that the solutions we construct in this way are local as they use only local information, namely the initial and terminal value of $Y$ in each segment, and are independent of each other. 

Let us denote by $\eta^\delta(n)_k = | c^\delta_k - c^\delta(n)_k|$ the local errors corresponding to each linear segment. It follows from the compactness assumption for ${\mathcal I}^{-1}_\delta({y})$, 
that there exists an $\hat{X}^\delta\in{\mathcal I}^{-1}_\delta({y})$ such that
\[
e_p(n;\delta) = \| \hat{\bf X}^\delta(n) - \hat{\bf X}^\delta\|_{p-var,[0,T]}.
\]
It follows from \eqref{eq: p-to-1 bound} that
\begin{eqnarray*}
   \| \hat{\bf X}^\delta(n) - \hat{\bf X}^\delta\|_{p-var,[0,T]} \leq  C_1 C_2(\hat{X}^\delta,\hat{X}^\delta(n))\frac{T}{N_\delta}\sum_{i=1}^{N_\delta} \eta^\delta(n)_i.
\end{eqnarray*}
Thus, since the local errors $\eta^\delta(n)_k$ are independent of each other, the $p$-variation error will converge uniformly with respect to $\delta$ if each of the errors also converges uniformly with respect to $\delta$. Without loss of generality, we consider $\eta^\delta(n)_1$. Then, ignoring any numerical errors in the computation of $F$ and $G$, $\eta^\delta(n)_1$ will be the standard error of the Newton-Raphson algorithm which, under appropriate conditions, satisfies
\[
\eta^\delta(n+1)_{1} \leq C_\delta  \eta^\delta(n)_{1}^2 \leq C^{-1}_\delta \left(C_\delta\eta^\delta(0)\right)^{2^n},
\]
where 
\[
C_\delta = \frac{1}{2}\sup_{c} |\partial^2_c F_{\delta,y_0}(c)|\cdot \sup_c|\partial_c F_{\delta,y_0}(c)|^{-1}.
\]
Assuming they are both bounded for fixed $\delta$, we have already seen that $|\partial_c F_{\delta,y_0}(c)|$ is of order $\delta$ and, in a similar way, we can show that the second derivative $|\partial^2_c F_{\delta,y_0}(c)|$ will be of order $\delta^2$. Thus, $C_\delta$ is of order $\delta$. Taking $c^\delta(0)_1 = 0$ implies $\eta^\delta(0)_1 = |c^\delta_1|$, which is of order $\delta^{\frac{1}{p}-1}$. Thus, $|C_\delta\eta^\delta(0)|$ will be of order $\delta^{\frac{1}{p}}$. Moreover, using \eqref{eq: p-to-1 bound constant} and the fact that $\|\hat{X}^\delta\|_{1-{\rm var},[0,T]}$ and $\|\hat{X}^\delta(n)\|_{1-{\rm var},[0,T]}$ will be of order $\delta^{\frac{1}{p}-1}$, we finally get that for $\delta$ sufficiently small,
\begin{equation}
\label{eq: NR error bound}
   \| \hat{\bf X}^\delta(n) - \hat{\bf X}^\delta\|_{p-var,[0,T]}  \lesssim  \delta^{\frac{1-p}{p}\lfloor p \rfloor}\delta^{-1}\delta^\frac{2^n}{p}.
\end{equation}
It follows that the upper bound converges to $0$ with $n$, uniformly in $\delta$, as for $n$ large enough, the error of the difference $\delta^\frac{2^n}{p}$ becomes very small very quickly and will compensate for the constant depending on $\delta$ that explodes as $\delta^{\frac{1-p}{p}\lfloor p \rfloor-1}\lesssim \delta^{-p}$.

\subsection{Signature approach}
\label{sec: signature}

We now present a novel approach that aims to numerically solve the discrete inverse problem directly, rather than through equation \eqref{eq:c}. We start by reformulating the discrete inverse problem of definition \ref{def: pwl discrete inverse problem} in terms of the signature of $X$. Although we still assume that nested dyadic piecewise linear approximations to $X$ converge, one of the advantages of this approach is that it has the potential to be extended to different families of approximations to $X$. 

We remind the reader that the signature of a path $X:[0,T]\to\R^m$ is defined as the collection of iterated integrals
\begin{equation}
\label{signature}
{\bf X}_{0,T} := \left( 1, \int_0^T dX_u, \int_0^T\int_0^{u_2} dX_{u_1}\otimes dX_{u_2}, \dots, \int_0^T\cdots\int_0^{u_2} dX_{u_1}\otimes\cdots\otimes dX_{u_n},\dots  \right)
\end{equation}
belonging to the tensor algebra on $\R^m$, $T(\R^m)$, with the logarithm of the signature (log-signature) belonging to the corresponding free Lie algebra. We call
\[ {\bf X}_{0,T}^n =  \underset{0<u_1<\cdots<u_n<T}{\int\cdots\int} dX_{u_1}\otimes\cdots\otimes dX_{u_n}\in(\R^m)^n\]
the $n^{\rm th}$ level of the signature and we denote by ${\bf X}_{0,T}^{[n]}$ the corresponding $n^{\rm th}$ level of the log-signature \cite{TerryBook}. The signature characterizes the path, up to tree-like equivalence \cite{Horatio}, thus providing an alternative representation for it. Moreover, for ${\bf Y}$ the signature of the response $Y:[0,T]\to\R^d$, the It\^o map corresponding to \eqref{eq:main} can be expressed as a linear map ${\bf I}$ from ${\bf X}$ to ${\bf Y}$ (we will also use ${\bf I}(y_0)$ when we want to emphasize initial condition $y_0$). We can reformulate definition \ref{def: pwl discrete inverse problem} of the discrete inverse problem in terms of signatures of paths, as follows:
\begin{definition} [Discrete Inverse Problem on Signatures - piecewise linear paths] 
\label{def:inverse problem signature}
Consider the family of piecewise linear paths $X^\delta({c}^\delta)$ on partition ${\mathcal D}_\delta$, parameterized by the gradients ${c}^\delta$. Given a vector  ${y}_{{\mathcal D}_\delta} = \left\{y_{k\delta}\right\}_{k=0}^{N_\delta}$ and a vector field $f\in{\rm Lip}(\gamma)$, for $\gamma>p\geq 1$, we say that $\hat{X}^\delta = X^\delta({c}_y^\delta)$ solves the discrete inverse problem for this parameterization if the following conditions are satisfied:
\begin{itemize}
\item[(a)] All but the first coordinate of the log-signature of the path $\hat{X}^\delta$ on $[(k-1)\delta,k\delta]$ are $0$, which is equivalent to the path being linear on that segment, i.e.
\[ ({\bf \hat{X}^\delta})_{(k-1)\delta,k\delta}^{[r]} = 0, \forall k = 1,\dots,N_{\delta}\ {\rm and}\ \forall r\geq 2. \]
\item[(b)] If ${\bf I}$ is the It\^o map defined on signatures, corresponding to \eqref{eq:main}, then
\[ {\bf \hat{Y}}^\delta = {\bf I} \cdot {\bf \hat{X}^\delta}.\]
\item[(c)] The first level of the log-signature of $\hat{Y}^\delta$ corresponds to the increments given by the observations, i.e.
\[ ({\bf \hat{Y}^\delta})^{[1]}_{(k-1)\delta,k\delta} = y_{k\delta} - y_{(k-1)\delta},\ \forall k=1,\dots,N.\]
\end{itemize}
\end{definition}


\subsubsection{Description of algorithm}
\label{algorithm}

The main idea of the algorithm is to solve a series of continuous inverse problems, replacing ${y}$ by a series of paths $Y(\delta, n)$ interpolating ${y}$ on ${\mathcal D}_\delta$ that are simple enough so that the problem can be solved, with the sequence constructed so that $Y(\delta, n)$ converges to $\hat{Y}^\delta$ uniformly in $\delta$. Then, for small enough $\delta$ and large enough $n$, a solution $X(\delta,n)$ to the continuous inverse problem for $Y(\delta, n)$ is a good approximation of the solution of the continuous inverse problem for ${y}$. To simplify the notation, we fix $\delta$ and simply write $X(n)$ and $Y(n)$ instead of $X(\delta, n)$ and $Y(\delta, n)$. We also write $\hat{X}$ and $\hat{Y}$ instead of $\hat{X}^\delta$ and $\hat{Y}^\delta$, unless dependence on $\delta$ needs to be emphasized.

\begin{figure}
  \includegraphics[scale=.75]{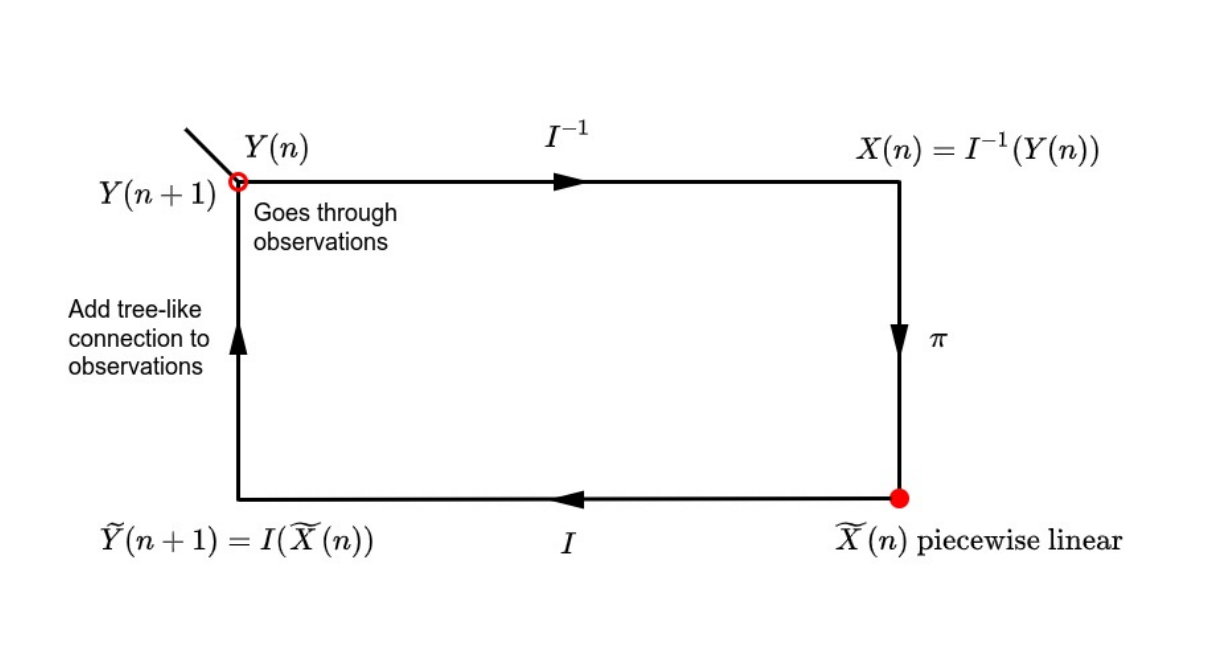}
  \caption{Schematic description of algorithm \ref{alg: sig}.}
\label{fig: algorithm}
\end{figure}

To make the continuous inverse problem for $Y(n)$ possible to solve, first we assume that $Y(n)$ are bounded variation paths. Then, we construct a sequence of bounded variation paths 
\[ \cdots \to Y(n) \to X(n) \to \tilde{X}(n) \to \tilde{Y}(n+1) \to Y(n+1) \to \cdots\]
each time trying to correct for conditions (a), (b) and (c) of definition \ref{def:inverse problem signature} by applying the smallest possible change to the signature of the corresponding path (see figure \ref{fig: algorithm}). While no pair $(X(n), Y(n))$ is expected to satisfy all three conditions, the goal is to construct a contraction map so that $Y(n)$ converges in $p$-variation to some solution of the discrete inverse problem $\hat{Y}^\delta\in{\mathcal I}^{-1}_\delta({y})$. Then, since $X(n)$ are the corresponding solutions to the continuous inverse problem for ${y}$ replaced by the trajectory of $Y(n)$ interpolating ${y}$ on ${\mathcal D}_\delta$, their limit will be a solution to the discrete inverse problem $\hat{X}^\delta$. Note that by the `sandwich principle', this is equivalent to $\tilde{X}(n)$ and $\tilde{Y}(n)$ converging in $p$-variation to $\hat{X}^\delta$ and $\hat{Y}^\delta$ respectively.

A key step in the algorithm is the solution of the continuous inverse problem for $Y(n)$, i.e. mapping $Y(n)$ to a bounded variation path $X(n)$, such that $I({X(n)}) = Y(n)$. We call this map the {\it inverse It\^o} map denoted by $I^{-1}$. 
In the case $m=d$, our original intuition was to try to solve the continuous inverse problem using the integral \eqref{eq: naive inverse}, with the shortcoming being that the integral was not well defined for path ${y}$. However, if we restrict ourselves to bounded variation paths $Y(n)$ that interpolate ${\mathcal D}_\delta$ while their range remains within a subset of $\R^d$ where the inverse of $f$, $g$ exists and is locally Lipschitz, then we can define $I^{-1}$ as mapping $Y(n)$ to an $X(n)$ defined as
\begin{equation}
\label{eq:inverse Ito}
{X}(n)_{0,t} = \int_0^t g({Y}(n)_u)\cdot d{Y}(n)_u,\ t\in[0,T].
\end{equation}

\begin{figure}
\centering
\includegraphics[scale=.6]{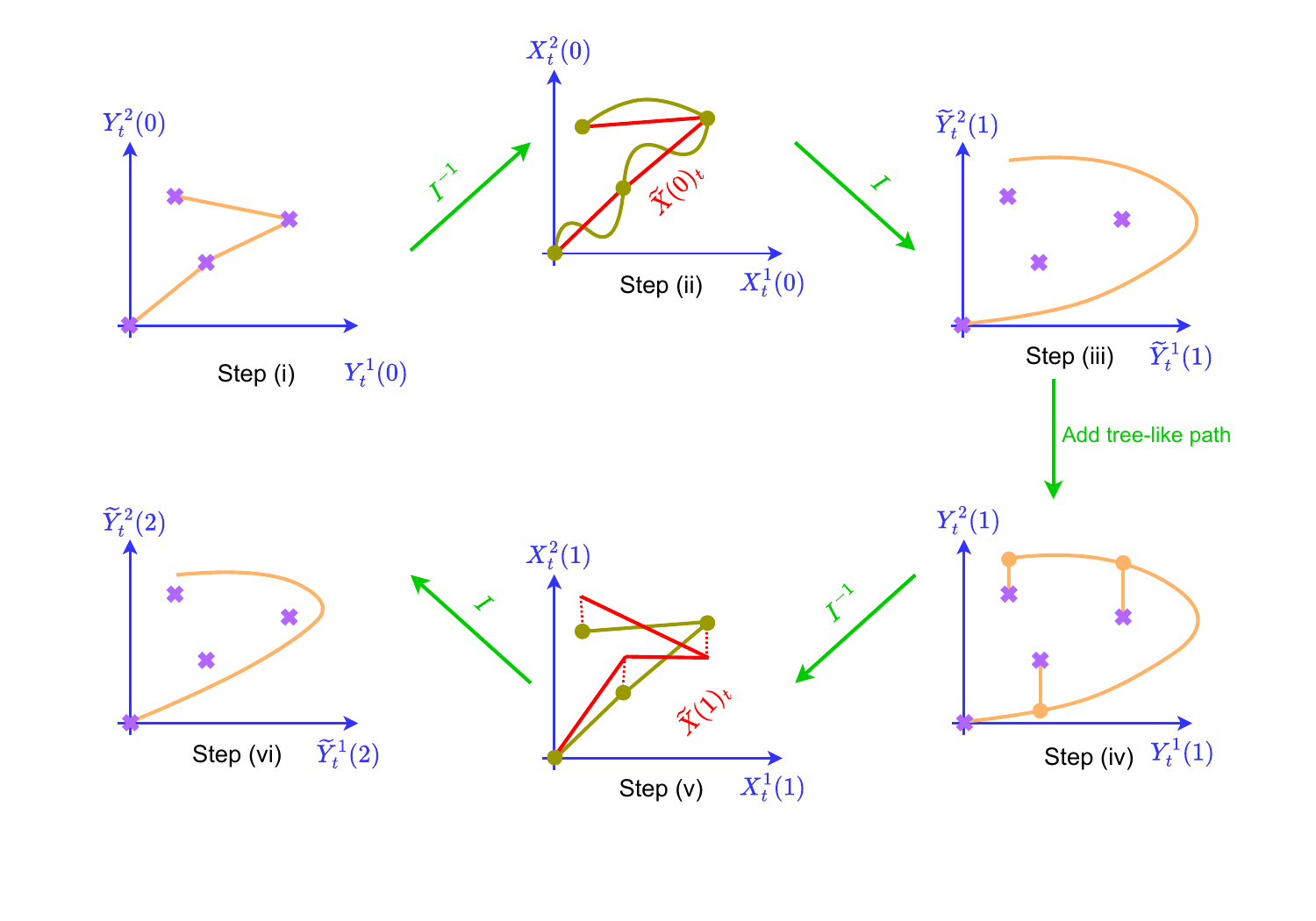}
\caption{Step-by-step description of algorithm \ref{alg: sig}.}
\label{fig: steps}
\end{figure}

By construction, the response to $X(n)$ interpolates the observations ${y}_{{\mathcal D}_\delta}$, thus satisfying conditions (b) and (c) of definition \ref{def:inverse problem signature}. However, for ${X}(n)$ to be the solution to the discrete inverse problem of definition \ref{def: pwl discrete inverse problem}, it also needs to satisfy condition (a), i.e. be piecewise linear on the partition ${\mathcal D}_\delta$. Thus, the next step in the algorithm is to correct for this by mapping $X(n)$ to a piecewise linear path $\tilde{X}(n)$ that approximates $X(n)$, with a natural choice being the piecewise linear interpolation. Note that ${X}(n)$ will normally depend on the whole continuous path ${Y}(n)$ on $[0,T]$, and not only its values on the partition ${\mathcal D}_\delta$ (note that this will not always be the case -- see section \ref{sec: reducibility}). Thus, different paths ${Y}(n)$ interpolating the observations ${y}_{{\mathcal D}_\delta}$ will lead to different paths ${X}(n)$ with possibly different increments on the partition, thus different $\tilde{X}(n)$. Now, the response $\tilde{Y}(n+1) = I(\tilde{X}(n))$ will no longer interpolate the observations, ${y}_{{\mathcal D}_\delta}$, failing to satisfy condition (c), and a new correction is needed (see figure \ref{fig: steps}).
This is where the intuition from using signatures becomes crucial, leading to the following:

\begin{algorithm}[on signatures]
\label{alg: sig}
Given a vector  ${y}_{{\mathcal D}_\delta} = \left\{y_{k\delta}\right\}_{k=0}^{N_\delta}$ and a vector field $f\in{\rm Lip}(\gamma)$, for $\gamma>p\geq 1$, defining the It\^o map for \eqref{eq:main}, with $I^{-1}$ a corresponding inverse It\^o map, we construct a sequence of bounded variation paths $(X(n),Y(n))$ as follows:

\noindent {\bf Step 0 (initialization)}:
\begin{itemize}
\item Initialize $Y(0)$ as the piecewise linear path going through  ${y}_{{\mathcal D}_\delta}$ (see figure \ref{fig: steps} (i)). Then, condition (c) of definition \ref{def:inverse problem signature} is satisfied. Note that alternative interpolations of the observations can also be considered (see section \ref{sec: extensions_1}).
\end{itemize}

\noindent {\bf Step $n\to n+1$}
\begin{itemize}
\item[1.] Given $Y(n)$ satisfying condition (c) by construction, define $X(n) = I^{-1}(Y(n))$, so that $(X(n), Y(n))$ satisfy conditions (b) and (c) of definition \ref{def:inverse problem signature} (see figure \ref{fig: steps} (ii) and (v)), but not condition (a).
\item[2.] Define $\tilde{X}(n)$ as the projection of $X(n)$ onto the space of piecewise linear paths on ${\mathcal D}_\delta$  (see figure \ref{fig: steps} (ii) and (v)). In terms of signatures, this means that $\forall k=1,\dots,N_\delta$
\[ ({\bf \tilde{X}}(n))^{[1]}_{(k-1)\delta,k\delta} = ({\bf X}(n))^{[1]}_{(k-1)\delta,k\delta}\ \ {\rm and}\ \ ({\bf \tilde{X}}(n))^{[r]}_{(k-1)\delta,k\delta}=0,\ r\geq 2.
\]
This is the minimum change to the signature on each segment, so that it corresponds to a linear path and  $(\tilde{X}(n), Y(n))$ satisfy conditions (a) and (c) but not (b).
\item[3.] Define $\tilde{Y}(n+1) = I(\tilde{X}(n))$, so that $(\tilde{X}(n), \tilde{Y}(n+1))$ satisfies conditions (a) and (b), but it now fails to satisfy (c) (see figure \ref{fig: steps} (iii) and (vi)).
\item[4.] Define $Y(n+1)$ by adding tree-like paths to connect to the observations, i.e. $Y(n+1)$ can be expressed as a concatenation of paths, with each segment corresponding to 
\[ 
L (\tilde{Y}(n+1)_{(k-1)\delta},Y_{(k-1)\delta})_{1,0}\otimes\tilde{Y}(n+1)_{(k-1)\delta,k\delta} \otimes L (\tilde{Y}(n+1)_{k\delta},Y_{k\delta})_{0,1},
\]
where for any two points $y,\tilde{y}\in\R^d$, $L(y,\tilde{y})_u = y + (\tilde{y}-y)u$ is the linear path connecting $y$ to $\tilde{y}$ for $u\in[0,1]$ (see figure \ref{fig: steps} (iv)). Then, $(\tilde{X}(n),Y(n+1))$ satisfy (a) and (c), but not (b). Note that by adding a tree-like path connecting to the observations, we are not changing the signature of the path over the whole interval $[0,T]$ \cite{Horatio}.
\end{itemize}
\end{algorithm}

\begin{remark}
\label{rem: partially known}
Algorithm \ref{alg: sig} assumes that the solution to the discrete inverse problem $\hat{X}^\delta$ is completely unknown. However, there are many examples where some coordinates of $\hat{X}^\delta$ might be known (e.g. when \eqref{eq:main} has a drift or when some coordinates are sampled from the distribution -- see remark \ref{rmk: uneven dimensions}). In that case, the algorithm remains the same except that we can initialize the known coordinates correctly. 
\end{remark}

Although the above description explains the intuition behind the algorithm, it is not very easy to implement. However, it can be simplified by using the piecewise linear assumption on $\tilde{X}(n)$. We can write
\begin{equation} 
\label{eq: alg piecewise linear}
\tilde{X}(n)_t = \tilde{X}(n)_{k\delta} + \tilde{c}(n)_k (t-k\delta),\ {\rm for}\ t\in[k\delta,(k+1)\delta],\ k=1,\dots,N_{\delta}-1.
\end{equation}
The following proposition gives an equivalent way of expressing algorithm \ref{alg: sig} as an update on vector ${\bf \tilde{c}}(n)= (\tilde{c}(n)_k)_{k=1}^{N_{\delta}}$.

\begin{proposition}
\label{prop: alg}
Let $I$ be the It\^o map for \eqref{eq:main} with vector field $f\in{\rm Lip}(\gamma)$, for $\gamma>p\geq 1$ and let $I^{-1}$ be a corresponding inverse It\^o map. Let ${\bf \tilde{c}}(n)= (\tilde{c}(n)_k)_{k=1}^{N_{\delta}}$ be defined as the gradients of the linear segments of the piecewise linear path $\tilde{X}(n)$ constructed through algorithm \ref{alg: sig}, given ${y}_{{\mathcal D}_\delta} = \left\{y_{k\delta}\right\}_{k=0}^{N_\delta}$. Then, ${\bf \tilde{c}}(n)= (\tilde{c}(n)_k)_{k=1}^{N_\delta}$ can be constructed directly as follows:

\noindent {\bf Step 0 (initialization)}:
\begin{itemize}
\item Let $Y(0)$ be the linear interpolation of the observations, as described in step 0 of algorithm \ref{alg: sig}. Then, for $k=1,\dots,N_\delta$, we define $\tilde{c}(0)_k$ as
\begin{equation}
\label{tilde c initialization}
\tilde{c}(0)_k = \frac{1}{\delta}I^{-1}\left(Y(0)\right)_{(k-1)\delta,k\delta} ,
\end{equation}
where we denote by $I^{-1}(Z)^1_{s,t}$ the increment of $I^{-1}(Z)$ on $[s,t]$ -- note that $[s,t]$ will be omitted when implied by the context. If some coordinates of the gradients are already known, then we initialize them to their actual value (see remark \ref{rem: algorith imprementation comments}).
\end{itemize}

\noindent {\bf Step $n\to n+1$}
\begin{itemize}
\item Given ${\bf \tilde{c}}(n)$, we define ${\bf \tilde{c}}(n+1)$ by
\begin{equation}
\label{tilde c update}
\tilde{c}(n+1)_k = \tilde{c}(n)_k + \left(r(n+1)_k - r(n+1)_{k-1}\right),
\end{equation}
where ${r}(n+1)_0 = 0$ and for $k=1,\dots,N_\delta$, 
\begin{equation}
\label{r}
r(n)_k =  \frac{1}{\delta}I^{-1}\left( L(\tilde{Y}(n)_{k\delta},y_{k\delta}) \right),
\end{equation}
where, as defined earlier, $L(y,\tilde{y})$ is the linear path connecting $y$ with $\tilde{y}$ and $\tilde{Y}(n)$ is the solution of \eqref{eq:main}, driven by $\tilde{X}(n)$ corresponding to ${\bf \tilde{c}}(n)$, i.e. $\tilde{Y}(n) = I(y_0,\tilde{X}(n)) = I(y_0,X^\delta(\tilde{c}(n)))$. 
\end{itemize}
\end{proposition}


\begin{proof}
Using the definition of $\tilde{X}(n)$, $X(n)$, $\tilde{Y}(n)$ and $Y(n)$, we write
\begin{eqnarray*}
  &&{X}(n+1)_{(k-1)\delta,k\delta} = I^{-1}\left(Y(n+1)\right)_{(k-1)\delta,k\delta} \\
  &=& I^{-1}\left( L(y_{(k-1)\delta},\tilde{Y}(n+1)_{(k-1)\delta},)\otimes\tilde{Y}(n+1)\otimes L(\tilde{Y}(n+1)_{k\delta},y_{k\delta})\right)
  \\
  &=& I^{-1}\left(
  L(y_{(k-1)\delta}, \tilde{Y}(n+1)_{(k-1)\delta})\right)\otimes 
  I^{-1}\left(\tilde{Y}(n+1)\right)\otimes I^{-1}\left(L(\tilde{Y}(n+1)_{k\delta},y_{k\delta})\right)  \\
  &=& I^{-1}\left(L(\tilde{Y}(n+1)_{(k-1)\delta},y_{(k-1)\delta})^{-1}\right)\otimes X(n)_{(k-1)\delta,k\delta} \otimes I^{-1}\left( L(\tilde{Y}(n+1)_{k\delta},y_{k\delta})\right).
\end{eqnarray*}
Moreover, by definition    
$$\tilde{X}(n+1)_{(k-1)\delta,k\delta} = X(n+1)_{(k-1)\delta,k\delta}$$ since $\tilde{X}(n+1)$ is the linearization of $X(n+1)$. Thus, 
\begin{eqnarray*}
    \tilde{X}(n+1)_{(k-1)\delta,k\delta} &=& X(n)_{(k-1)\delta,k\delta} - I^{-1}\left(L(\tilde{Y}(n+1)_{(k-1)\delta},y_{(k-1)\delta})\right) \\
    &&+ I^{-1}\left( L(\tilde{Y}(n+1)_{k\delta},y_{k\delta})\right).
\end{eqnarray*}
The results follows by dividing both sides by $\delta$.
\end{proof}


\begin{remark} 
\label{rem: algorith imprementation comments}
In practice, there are a few more things to consider when implementing the algorithm of proposition \ref{prop: alg}.
\begin{itemize}
\item Similarly to the Newton-Raphson scheme, this algorithm is designed to solve the system of equations required to compute the gradients of a solution to the discrete inverse problem. Both algorithms require the solution $F$ of the ODE corresponding to the It\^o map -- in the case of Newton-Raphson, this is explicit, while for the signature-based algorithm, this is hidden in the computation of $\tilde{Y}(n)$. Newton-Raphson algorithm also requires the computation of the derivative $G$ of the solution to the ODE with respect to constant $c$, which is not required by the signature-based algorithm. Instead, it requires the computation of the integral corresponding to the solution of the inverse It\^o map. In most cases, numerical methods are needed to compute $F$, $G$ and the integral.
\item When $\delta$ is small, it might be necessary to scale constants $\tilde{c}(n)_k$ and $r(n)_k$ by multiplying by $\delta$. Then, an alternative initialization is $\tilde{c}(0) = 0$. 
\item In order to be able to apply the inverse It\^o map in \eqref{r} to the tree-like corrections, we need to make sure that they are within its range. For example, when using \eqref{eq:inverse Ito}, the tree-like corrections should be such that their range of values is within a set where $g$ is locally Lipschitz. When this fails, we reinitialize the gradient to the corresponding segment. As will become clear in the proof of Lemma \ref{lemma:contraction}, the initialization error must be sufficiently small for the algorithm to converge uniformly. This requirement is implicitly linked to the time horizon, $T$; as $T$ grows, more frequent reinitializations will be needed.  
\end{itemize}
\end{remark}

Note that, as the solution to the inverse problem given a continuous path $X(n)$ is uniquely defined through \eqref{eq:inverse Ito}, it is the initialization of the algorithm that determines which solution of the discrete inverse problem we will converge to. Thus, for fixed initialization, we have a unique solution. In what follows, we assume that the solution to the discrete inverse problem $\hat{\bf X}$ is unique.
\color{black}

\subsubsection{Proof of Convergence}
\label{sec: convergence}

As before, we define the error of the signature-based numerical method as
\[
e_p(n;\delta) = d_p\left( \tilde{\bf X}^\delta(n), {\mathcal I}^{-1}_\delta({y}) \right),
\]
where the numerical solution is now denoted by $\tilde{X}^\delta(n)$ to emphasize dependence on $\delta$, but this is omitted when implied by the context. We want to show that 
\[ 
\lim_{n\to\infty}\sup_{\delta>0} e_p(n;\delta) = 0.
\]
Similar to the proof of convergence of the Newton-Raphson algorithm, we will first show convergence of 
\begin{equation}
\label{eq: p-var upper bound}
    \delta \sum_{i=1}^{N_\delta} |c^\delta_i-\tilde{c}^\delta(n)_i| = \sum_{i=1}^{N_\delta} |\hat{X}^\delta_{(i-1)\delta,i\delta}-\tilde{X}^\delta(n)_{(i-1)\delta,i\delta}|
\end{equation}
with respect to $n$. Note that, unlike the Newton-Raphson case, individual errors $|c^\delta_i-\tilde{c}^\delta(n)_i|$ are no longer independent and the interaction between segment errors needs to be taken into consideration when proving convergence of the sum. The proof is based on showing that \eqref{eq: p-var upper bound} is a contraction with respect to $n$. 
Moreover, the contraction coefficient goes to zero with $\delta\to 0$,
which will allow us to compensate for the exploding constant in the upper bound of the $p$-variation norm given by \eqref{eq: p-to-1 bound} and thus prove that $e_p(n;\delta)$ converges to $0$ with $n\to\infty$ uniformly with respect to $\delta$.
\color{black}

The proof of contraction is based on the following observation: the way the signature-based algorithm works is by correcting the paths at each step to match the assumptions, either forcing an assumption on the first level of the log-signature (the path has to agree with observations) or all other levels (the path needs to be piecewise linear). As such, the key idea to understand the evolution of the error in each segment and prove contraction is to express the It\^o and inverse It\^o maps as linear maps acting on monomials of the log-signature. This allows us to separate the effect of the linear and nonlinear parts of the input path on the linear and nonlinear parts of the output path. We remind the reader that the signature can be represented as linear combinations of monomials on the log-signature and that the It\^o and inverse It\^o maps can be viewed as linear maps on the signature.

To make this more concrete, we introduce some additional notation: for any rough path ${\bf Z}$ on a fixed interval $[s,t]$, we denote by ${\bf P}_{1}(Z,s,t)$ the vector of all monomials of first level of the log-signature on $[s,t]$ (i.e. the increments) and by ${\bf P}_{R}(Z,s,t)$ the vector of all remaining monomials -- to simplify notation, we will drop dependence on $s,t$ when this is implied by the context. Moreover, we use the decomposition of the linear It\^o map ${\bf I}$ acting on the signature into ${\bf I}^f_{ij}$, with $i,j\in\{ 1,R\}$, to separate the part contributing to ($i$) or acting on ($j$) monomials of the increment ($i,j=1$) and the part acting on remaining monomials ($i,j=R$), and similarly for ${\bf I}_{ij}^b$ corresponding to the decomposition of the inverse It\^o map.

Expressing the It\^o map as a linear map on the signature will always hold under the Lipschitz assumption on $f$ (see \cite{Gluckstad}). Thus, using the above notation, for any rough path ${\bf X}$ and $Y=I(y,{\bf X})$ we can write
\begin{equation}
\label{eq: forward Ito decomposition}
    \left\{ 
    \begin{array}{ccccc}
    {\bf P}_{1}(Y,s,t) &=& {\bf I}^f_{11}(y_s) {\bf P}_{1}(X,s,t) &+& {\bf I}^f_{1R}(y_s) {\bf P}_{R}(X,s,t) \\
    {\bf P}_{R}(Y,s,t) &=& {\bf I}^f_{R1}(y_s) {\bf P}_{1}(X,s,t) &+& {\bf I}^f_{RR}(y_s) {\bf P}_R(X,s,t).
    \end{array}
    \right.   
\end{equation}
The inverse It\^o map corresponds to integration, so under the assumption that $g$ exists and is (locally) Lipschitz, we can also write 
\begin{equation}
\label{eq: backward Ito decomposition}
    \left\{ 
    \begin{array}{ccccc}
    {\bf P}_1(X,s,t) &=& {\bf I}^b_{11}(y_s) {\bf P}_1(Y,s,t) &+& {\bf I}^b_{1R}(y_s) {\bf P}_R(Y,s,t) \\
    {\bf P}_R(X,s,t) &=& {\bf I}^b_{R1}(y_s) {\bf P}_1(Y,s,t) &+& {\bf I}^b_{RR}(y_s) {\bf P}_R(Y,s,t).
    \end{array}
    \right..   
\end{equation} 

\begin{remark}
    When the inverse It\^o map is acting on a path $Y$, we use the notation $X= I^{-1}(Y)$ -- we do not include initial conditions for $X$ as we are typically interested in increments of $X$. However, when we use the corresponding linear map ${\bf I}^b$, this will also depend on the initial value of $Y$, as it is relevant and is lost when replacing the path $Y$ by its signature ${\bf Y}$.
\end{remark}
The following lemma gives the contraction result. 
\begin{lemma}
\label{lemma:contraction}
    Suppose that $\hat{X}^\delta$ is the solution to the discrete inverse problem defined in \ref{def: pwl discrete inverse problem} and the sequence $(\tilde{X}^\delta(n))_{n=0}^\infty$ is constructed by applying algorithm \ref{alg: sig}. Assume that $f$ is such that \eqref{eq: forward Ito decomposition} and \eqref{eq: backward Ito decomposition} hold on any segment $[s,t]$, with coefficients that are uniformly bounded on ${y}$ and $\omega((k-1)\delta,k\delta)\sim{\mathcal O}(\delta)$ for each $k=1,\dots,N_\delta$. 
    Then \eqref{eq: p-var upper bound} is a contraction with respect to $n$, i.e. 
    \begin{equation}
        \sum_{i=1}^{N_\delta} |\hat{X}^\delta_{(i-1)\delta,i\delta}-\tilde{X}^\delta(n+1)_{(i-1)\delta,i\delta}| \leq C_\delta \sum_{i=1}^{N_\delta} |\hat{X}^\delta_{(i-1)\delta,i\delta}-\tilde{X}^\delta(n)_{(i-1)\delta,i\delta}|
    \end{equation}
    where $C_\delta$ is of order ${\mathcal O}(\delta^\frac{1}{p})$, provided that the initialization error $\|\tilde{X}(0)-\hat{X}\|_{p-var,[0,T]}$ is sufficiently small.
\end{lemma}

\begin{proof}
We start by deriving a bound for $|\tilde{X}^\delta(n+1)_{(i-1)\delta,i\delta}-\hat{X}^\delta_{(i-1)\delta,i\delta}|$ in terms of past errors $|\tilde{X}^\delta(n)_{(j-1)\delta,j\delta}-\hat{X}^\delta_{(j-1)\delta,j\delta}|$, for $j=1,\dots i$, using \eqref{eq: forward Ito decomposition} and \eqref{eq: backward Ito decomposition}. Note that the dependence on past errors is a fundamental difference from the Newton-Raphson algorithm and originates from the fact that we are forcing $Y(n)$ to be continuous. 

First, we consider the case $i=1$, which corresponds to the segment $[0,\delta]$. In this case, the corresponding $\tilde{Y}(n)$ will not need a tree-like correction at the starting point, since the initial conditions are known. It follows from proposition \ref{prop: alg} that 
\begin{equation}
\label{eq: signature error 1}
    \tilde{X}(n+1)_{0,\delta} - \hat{X}_{0,\delta}= \left(\tilde{X}(n)_{0,\delta} - \hat{X}_{0,\delta}\right)+ I^{-1}(L(\tilde{Y}(n+1)_\delta,y_\delta)).
\end{equation}
It follows from the linearity assumption on $L(\tilde{Y}(n+1)_\delta,y_\delta)$ and the definition of ${\bf P}_R$ as all monomials of the log-signature that depend on at least one level of the log-signature higher than 1 that 
\[
{\bf P}_R(L(\tilde{Y}(n+1)_\delta,y_\delta))) \equiv  {\bf 0}.
\]
Thus, using \eqref{eq: backward Ito decomposition} to express the correction term $I^{-1}(L(\tilde{Y}(n+1)_\delta,y_\delta))$, we write
\begin{eqnarray}
\label{eq: correction term 1}
\nonumber {\bf P}_1\left(I^{-1}(L(\tilde{Y}(n+1)_\delta,y_\delta))\right) &=& {\bf I}^b_{11}(\tilde{Y}(n+1)_\delta)\cdot {\bf P}_1\left( L(\tilde{Y}(n+1)_\delta,y_\delta)\right)
\\
\nonumber &=&{\bf I}^b_{11}(\tilde{Y}(n+1)_\delta)\cdot {\bf P}_1\left(\tilde{Y}(n+1)^{-1}\otimes\hat{Y}\right)\\
\nonumber ={\bf P}_1\left(I^{-1}(\tilde{Y}(n+1)^{-1}\otimes\hat{Y})\right) &-&{\bf I}^b_{1R}(\tilde{Y}(n+1)_\delta)\cdot {\bf P}_R\left(\tilde{Y}(n+1)^{-1}\otimes\hat{Y}\right)\\
= {\bf P}_1\left(\tilde{X}(n)^{-1}\otimes\hat{X}\right) &-& {\bf I}^b_{1R}(\tilde{Y}(n+1)_\delta) \cdot {\bf P}_R\left(\tilde{Y}(n+1)^{-1}\otimes\hat{Y}\right),
\end{eqnarray}
using the fact that the increment of ${L}(\tilde{Y}(n+1)_\delta,y_\delta)$, equal to $y_\delta-\tilde{Y}(n+1)_\delta$ is the same as the increment of $\tilde{Y}(n+1)^{-1}\otimes\hat{Y}$ given by $(Y(n+1)_0-Y(n+1)_\delta)+(y_\delta-y_0)$, since $Y(n+1)_0 = y_0$ (which is always true for the first segment), and thus monomials of the increments will also be the same. Note that, using standard notation, for any rough path ${\bf Z}$ we denote by $Z^{-1}$ the time-reversed path.
Moreover, by definition,
\[
I^{-1}\left(Y(n+1)^{-1}\otimes \hat{Y}\right) = I^{-1}\left(Y(n+1)\right)^{-1}\otimes I^{-1}(\hat{Y}) = \tilde{X}(n)^{-1}\otimes \hat{X}.
\]
Since \eqref{eq: correction term 1} holds for all monomials of the increment, it will hold for the increment itself (monomial of order 1), i.e. 
\[
I^{-1}(L(\tilde{Y}(n+1)_\delta,y_\delta)) = \hat{X}_{0,\delta} - \tilde{X}(n)_{0,\delta} - \left({\bf I}^b_{1R}(\tilde{Y}(n+1)_\delta) \cdot {\bf P}_R(\tilde{Y}(n+1)^{-1}\otimes\hat{Y})\right)^1
\]
Plugging this into \eqref{eq: signature error 1}, we get
\begin{equation}
\label{eq: correction term 2}
 \tilde{X}(n+1)_{0,\delta} - \hat{X}_{0,\delta}=- \left({\bf I}^b_{1R}(\tilde{Y}(n+1)_\delta) \cdot {\bf P}_R(\tilde{Y}(n+1)^{-1}\otimes\hat{Y})\right)^1.
\end{equation}

The next step is to expand the operator ${\bf I}^b_{1R}$. To do this, we need to construct a basis for the monomials of the log-signature. Let ${\mathcal H}$ be the set of Hall words on $\{1,\dots,d\}$ and let ${w} = (w_1,\dots,w_n)\in{\mathcal H}^n$. Then, ${\mathcal H}^\infty = \cup_{n\geq 1}{\mathcal H}^n$ provides such a basis. 
For an arbitrary rough path $Z$ taking values in ${\mathbb R}^d$, let $P_{w_i}(Z)$ be the log-signature of $Z$ on Hall word $w_i\in{\mathcal H}$. We define the {\it total order} of the order $n$ monomial $P_{w}(Z)=\prod_{i=1}^n P_{w_i}(Z)$ as $|{w}| = |w_1|+\cdots+|w_n|$ where $|w_i|$ is the level of the log-signature corresponding to Hall word $w_i$. Moreover, suppose that ${\mathcal H}^{[1]},{\mathcal H}^{[R]}\subset {\mathcal H}^\infty$ are a partition of the space of sequences of Hall words, where ${w}= (w_1,\dots,w_n)\in{\mathcal H}^{[1]}$ if $|{w}|= n$ and ${w}\in{\mathcal H}^{[R]}$ otherwise. 


By definition, ${\bf I}^b_{1R}$ is acting exactly on monomials $P_{w}$ where ${w}\in{\mathcal H}^{[R]}$. In particular, the lowest order polynomials it acts on will be the order 1 monomials of the area, i.e. $P_{[i,j]}$, whose total order is 2. For clarity of the exposition, we only consider the case where $d=2$, as the proof for $d>2$ is the conceptually the same. Then, by truncating the sum at total order 2, we get
\begin{equation}
\label{eq: inverse Ito Taylor expansion}
  \left({\bf I}^b_{1R}(z) {\bf P}_R(Z)\right)^{(i)} = \psi_{[2,1]}^i(\zeta_z)P_{[2,1]}(Z)
\end{equation}
for some $\zeta_z$ in a neighborhood of $z$, where $\psi^{w^\prime}_{w}(\cdot) $ are the coefficients in the expansion of ${\bf I}^b$ depending on initial conditions. This follows from the Taylor expansion of $g$, under the locally Lipschitz assumption and the interpretation of the inverse It\^o map as an integral given by \eqref{eq:inverse Ito}. 

Similarly, for $Z = I(z,\tilde{Z})$, we can write
\[
    P_{[2,1]}(Z) = \left({\bf I^f_{R1}}(z){\bf P}_1(\tilde{Z})+{\bf I^f_{RR}}(z){\bf P}_R(\tilde{Z})\right)^{[2,1]}
    = \sum_{{w}\in{\mathcal H}^\infty}\phi^{[2,1]}_{w}(z) P_{w}(\tilde{Z})
\]
(see \cite{Gluckstad}). Since the left-hand-side cannot contain terms of lower order than the right-hand-side, it follows that $\phi^{[2,1]}_{w}(y) \equiv 0$, for $|{w}|=1$. As before, assuming $\delta$ sufficiently small so that $\|\tilde{Z}\|_{p-var,[0,\delta]}$ is sufficiently small, we can truncate at total order 2 (see \cite{Baudoin}, Prop. 2.5), leading to
\[
    P_{[2,1]}(Z) = \sum_{|{w}|=2}\phi^{[2,1]}_{w}(z) P_{w}(\tilde{Z}) + {\mathcal O}(\|{\tilde{Z}}\|_{p-var,[0,\delta]}^2),
\]
where $\phi^{w^\prime}_{w}(\cdot) $ are now the coefficients in the expansion of ${\bf I}^f$ depending on initial conditions. 

By plugging in the expansion
\[
Z^{i} = \sum_{{w}\in{\mathcal H}^\infty}\phi^{i}_{w}(z) P_{w}(\tilde{Z})
\]
for $i=1,2$, into the expression
\[
P_{[2,1]}(Z,s,t) = \frac{1}{2}\left(Z^{(2,1)}_{s,t})-Z^{(1,2)}_{s,t}\right) = \frac{1}{2}\left(\int_s^t\int_s^u dZ^{(2)}_{u_2}dZ^{(1)}_u - \int_s^t\int_s^u dZ^{(1)}_{u_2}dZ^{(2)}_u\right)
\]
we find that
\[
\phi^{[2,1]}_w \equiv 0, \forall w\in{\mathcal H}^{[1]},
\]
that is, the area term of the response $P_{[2,1]}(Z)$ only depends on the area term of the control $\tilde{Z}$ and higher-order terms. Thus, for $Z =\tilde{Y}(n+1)^{-1}\otimes\hat{Y}$ and $\delta$ sufficiently small so that $\|\tilde{X}(n)^{-1}\otimes\hat{X}\|_{p-var,[0,\delta]}$ we get that 
\begin{eqnarray}
\label{eq: one segment expansion}
\nonumber \tilde{X}(n+1)_{0,\delta} - \hat{X}_{0,\delta} &=&-\psi_{[2,1]}^1(y_{\xi_{n+1}})\phi^{[2,1]}_{[2,1]}(\tilde{Y}(n+1)_\delta)
P_{[2,1]}\left(\tilde{{ X}}(n)^{-1}\otimes{\hat{X}}\right) \\
&&+ {\mathcal O}(\|\tilde{{X}}(n)^{-1}\otimes{\hat{X}}\|_{p-var,[0,\delta]}^2),
\end{eqnarray}
for $y_{\xi_{n+1}}$ in a neighborhood of $\tilde{Y}(n+1)_\delta$.
By expanding the tensor product $\tilde{{ X}}(n)^{-1}\otimes{\hat{X}}$, we get 
\[
P_{[2,1]}(\tilde{{ X}}(n)^{-1}\otimes{\hat{X}}) = P_{[2,1]}(\tilde{{ X}}(n)^{-1})+P_{[2,1]}({\hat{X}}) + P_1(\tilde{X}(n))P_2(\hat{X})- P_2(\tilde{X}(n))P_1(\hat{X}).
\]
Since $\tilde{X}(n)$ and $\hat{X}$ are linear on one segment, this simplifies to
\begin{eqnarray*}
P_{[2,1]}(\tilde{{ X}}(n)^{-1}\otimes{\hat{X}}) &=& \tilde{X}(n)^{(1)}_{0,\delta}\hat{X}^{(2)}_{0,\delta}-\tilde{X}(n)^{(2)}_{0,\delta}\hat{X}^{(1)}_{0,\delta}\\
&=& (\tilde{X}(n)^{(1)}_{0,\delta}-\hat{X}^{(1)}_{0,\delta})\hat{X}^{(2)}_{0,\delta}-(\tilde{X}(n)^{(2)}_{0,\delta}-\hat{X}^{(2)}_{0,\delta})\hat{X}^{(1)}_{0,\delta}.
\end{eqnarray*}
It follows that
\[
|P_{[2,1]}(\tilde{{ X}}(n)^{-1}\otimes{\hat{X}})| \leq |\tilde{X}(n)_{0,\delta}-\hat{X}_{0,\delta}|\omega(0,\delta)^\frac{1}{p}.
\]
Following a similar approach to the one above, we can show that the $\frac{p}{k}$-variation of path $\left(\tilde{{ X}}(n)^{-1}\otimes{\hat{X}}\right)^k_{0,\delta}$ (the $k^{\rm th}$ level of the signature) is of order ${\mathcal O}\left(|\tilde{X}(n)_{0,\delta}-\hat{X}_{0,\delta}|\omega(0,\delta)^\frac{k-1}{p} \right)$, which implies that
\[
\|\tilde{{ X}}(n)^{-1}\otimes{\hat{X}}\|_{p-var,[0,\delta]} \sim {\mathcal O}\left(|\tilde{X}(n)_{0,\delta}-\hat{X}_{0,\delta}|\right).
\]
By the assumption that coefficients of the expansions \eqref{eq: forward Ito decomposition} and \eqref{eq: backward Ito decomposition} are uniformly bounded, it follows from \eqref{eq: one segment expansion} that
\[
|\tilde{X}(n+1)_{0,\delta} - \hat{X}_{0,\delta}| \leq C_1 |\tilde{X}(n)_{0,\delta}-\hat{X}_{0,\delta}|\omega(0,\delta)^\frac{1}{p}
\]
for some constant $C_1$ that depends on $f$ and $\delta$ and it is ${\mathcal O}(1)$ with respect to $\delta$. This gives us a contraction estimate for the first interval. 

We now consider an arbitrary interval $[k\delta,(k+1)\delta]$. The main difference with interval $[0,\delta]$ is that we now need to correct the value of $\tilde{Y}(n+1)$ at the start point as well. 
As a result, \eqref{eq: correction term 1} now becomes
\begin{eqnarray*}
&& {\bf P}_1\left(I^{-1}(L(\tilde{Y}(n+1)_{(k+1)\delta},y_{(k+1)\delta}))\right) 
\\ 
&=& {\bf I}^b_{11}(\tilde{Y}(n+1)_{(k+1)\delta})\cdot {\bf P}_1\left(\tilde{Y}(n+1)^{-1}\otimes L(\tilde{Y}(n+1)_{k\delta},y_{k\delta}))\otimes\hat{Y}\right)
\\
&=& {\bf P}_1\left(\tilde{X}(n)^{-1}\otimes I^{-1}(L(\tilde{Y}(n+1)_{k\delta},y_{k\delta})))\otimes\hat{X}\right)  
\\
&&-{\bf I}^b_{1R}(\tilde{Y}(n+1)_\delta) \cdot {\bf P}_R\left(\tilde{Y}(n+1)^{-1}\otimes L(\tilde{Y}(n+1)_{k\delta},y_{k\delta})\otimes\hat{Y}\right),
\end{eqnarray*}
Note that \eqref{tilde c update} of proposition \ref{prop: alg} now corresponds to 
\begin{eqnarray*}
    \tilde{X}(n+1)_{k\delta,(k+1)\delta} - \hat{X}_{k\delta,(k+1)\delta}&=& \left(\tilde{X}(n)_{k\delta,(k+1)\delta} - \hat{X}_{k\delta,(k+1)\delta}\right)\\
    + I^{-1}(L(\tilde{Y}(n+1)_{(k+1)\delta},y_{(k+1)\delta})) 
& - &I^{-1}(L(\tilde{Y}(n+1)_{k\delta},y_{k\delta})).
\end{eqnarray*}
As before, the terms $I^{-1}(L(\tilde{Y}(n+1)_{k\delta},y_{k\delta})) $ and $I^{-1}(L(\tilde{Y}(n+1)_{(k+1)\delta},y_{(k+1)\delta}))$ cancel, leading to
\begin{eqnarray*}
&& \tilde{X}(n+1)_{k\delta,(k+1)\delta} - \hat{X}_{k\delta,(k+1)\delta}= \\
&&- \left({\bf I}^b_{1R}(\tilde{Y}(n+1)_{(k+1)\delta}) \cdot {\bf P}_R(\tilde{Y}(n+1)^{-1}\otimes L(\tilde{Y}(n+1)_{k\delta},y_{k\delta})\otimes\hat{Y})\right)^1,
\end{eqnarray*}
which, when expanding the linear It\^o and inverse It\^o maps, becomes
\begin{eqnarray}
\label{eq: X step inner segment}
\nonumber && \tilde{X}(n+1)_{k\delta,(k+1)\delta} - \hat{X}_{k\delta,(k+1)\delta} =
 \\
 \nonumber && -\psi_{[2,1]}^1(y_{\xi_{k,n+1}})\phi^{[2,1]}_{[2,1]}(\tilde{Y}(n+1)_{(k+1)\delta})
P_{[2,1]}\left(\tilde{{ X}}(n)^{-1}\otimes I^{-1}(L(\tilde{Y}(n+1)_{k\delta},y_{k\delta}))\otimes{\hat{X}}\right) \\
&& \ \ \ \ \ + {\mathcal O}(\|\tilde{{ X}}(n)^{-1}\otimes I^{-1}(L(\tilde{Y}(n+1)_{k\delta},y_{k\delta}))\otimes{ \hat{X}}\|_{p-var,[k\delta,(k+1)\delta]}^2).
\end{eqnarray}
We assume that $\|\tilde{{ X}}(n)^{-1}\otimes I^{-1}(L(\tilde{Y}(n+1)_{k\delta},y_{k\delta}))\otimes{ \hat{X}}\|_{p-var,[k\delta,(k+1)\delta]}$ is small so that the expansion above holds, noting that it is now not sufficient that $\delta$ is small. We also need that $\|I^{-1}(L(\tilde{Y}(n+1)_{k\delta},y_{k\delta}))\|_{p-var,[k\delta,(k+1)\delta]}$ is small, which, as we will see below, requires that $|\tilde{Y}(n+1)_{k\delta}-y_{k\delta}|$ is small -- i.e., the error in the initial conditions, coming from the previous segment, also matters. We will get back to this later. 

Expanding $P_{[2,1]}\left(\tilde{{ X}}(n)^{-1}\otimes I^{-1}(L(\tilde{Y}(n+1)_{k\delta},y_{k\delta}))\otimes{\hat{X}}\right)$ gives
\begin{eqnarray*}
&   P_{[2,1]}\left(\tilde{{ X}}(n)^{-1}\otimes{\hat{X}}\right) + P_{[2,1]}\left(I^{-1}(L(\tilde{Y}(n+1)_{k\delta},y_{k\delta}))\right)
    \\
 &   +\frac{1}{2}\left( (\tilde{X}(n)^{(2)}_{k\delta,(k+1)\delta}+\hat{X}^{(2)}_{k\delta,(k+1)\delta})I^{-1}(L(\tilde{Y}(n+1)_{k\delta},y_{k\delta}))^{(1)}\right.
\\   
&\left.\ \ \ \ \ \ \ \ \ \ \ \ \ \ \ \ \ \ \ -  (\tilde{X}(n)^{(1)}_{k\delta,(k+1)\delta}+\hat{X}^{(1)}_{k\delta,(k+1)\delta})I^{-1}(L(\tilde{Y}(n+1)_{k\delta},y_{k\delta}))^{(2)}\right).
\end{eqnarray*}
By expanding the operator $I^{-1}$ as before and noting that $L(\tilde{Y}(n+1)_{k\delta},y_{k\delta})$ is linear, we get that
\[
\left|P_{[2,1]}\left(I^{-1}(L(\tilde{Y}(n+1)_{k\delta},y_{k\delta}))\right)\right| \lesssim |\tilde{Y}(n+1)_{k\delta}-y_{k\delta}|^2
\]
and
\[
\left|I^{-1}(L(\tilde{Y}(n+1)_{k\delta},y_{k\delta}))^{(i)} \right|\lesssim |\tilde{Y}(n+1)_{k\delta}-y_{k\delta}|.
\]
To simplify notation, we use $\lesssim$ to denote an inequality up to a constant, with constants depending on $f$ and are ${\mathcal O}(1)$ with respect to $\delta$. 
Finally, drawing all the estimates together, we get that
\begin{eqnarray}
\label{eq: X one step correction different initial conditions}
\nonumber |\tilde{X}(n+1)_{k\delta,(k+1)\delta} - \hat{X}_{k\delta,(k+1)\delta}| &\lesssim& |\tilde{X}(n)_{k\delta,(k+1)\delta} - \hat{X}_{k\delta,(k+1)\delta}|\omega(k\delta,(k+1)\delta)^\frac{1}{p} 
\\&& + |\tilde{Y}(n+1)_{k\delta}-y_{k\delta}|\omega(k\delta,(k+1)\delta)^\frac{1}{p} \nonumber
\\&& +|\tilde{Y}(n+1)_{k\delta}-y_{k\delta}|^2
\end{eqnarray}
In other words, the error for each segment will depend on the previous error on this segment as in Newton-Raphson but also on the initial condition error, which will depend on all previous errors. Noting that $\tilde{Y}(n+1)_{k\delta} = I(\tilde{Y}(n+1)_{(k-1)\delta},\tilde{X}(n))$ and $y_{k\delta} = I(y_{(k-1)\delta},\hat{X})$, we write
\begin{eqnarray}
\label{eq: Y one step correction}
\nonumber    |\tilde{Y}(n+1)_{k\delta}-y_{k\delta}| &\leq& |I(\tilde{Y}(n+1)_{(k-1)\delta},\tilde{X}(n))-I(\tilde{Y}(n+1)_{(k-1)\delta},\hat{X})| 
\\ && + |I(\tilde{Y}(n+1)_{(k-1)\delta},\hat{X})-I(y_{(k-1)\delta},\hat{X})| 
\nonumber    
\\ &\lesssim& |\tilde{X}(n)_{(k-1)\delta,k\delta}-\hat{X}_{(k-1)\delta,k\delta}| \nonumber
    \\ && + |\tilde{Y}(n+1)_{(k-1)\delta}-y_{(k-1)\delta}| \omega((k-1)\delta,k\delta)^\frac{1}{p}
\end{eqnarray}
where the first inequality is just the triangle inequality and the second is an application of Lyons' continuity theorem for linear processes $\tilde{X}(n)$ and $\hat{X}$ for the first term and the initial value correction (see \cite{FrizVictoir}, Lemma 10.6, p. 216) for the second. Plugging \eqref{eq: Y one step correction} into \eqref{eq: X one step correction different initial conditions} and assuming that $\omega(j\delta,(j+1)\delta)\sim {\mathcal O}(\delta)$ for all $j=0,\dots,N_{\delta}-1$, we get
\begin{eqnarray*}
    |\tilde{X}(n+1)_{k\delta,(k+1)\delta} - \hat{X}_{k\delta,(k+1)\delta}| \lesssim |\tilde{X}(n)_{k\delta,(k+1)\delta} - \hat{X}_{k\delta,(k+1)\delta}|\delta^\frac{1}{p}
\\  + |\tilde{X}(n)_{(k-1)\delta,k\delta}-\hat{X}_{(k-1)\delta,k\delta}|\delta^\frac{1}{p} + |\tilde{Y}(n+1)_{(k-1)\delta}-y_{(k-1)\delta}|\delta^\frac{2}{p}.
\end{eqnarray*}
Note that the square term $|\tilde{Y}(n+1)_{k\delta}-y_{k\delta}|^2$ only contributes terms of the same or higher order than $|\tilde{Y}(n+1)_{k\delta}-y_{k\delta}|\omega(k\delta,(k+1)\delta)^\frac{1}{p}$. Inductively, this leads to
\begin{eqnarray*}
    &|\tilde{X}(n+1)_{k\delta,(k+1)\delta} - \hat{X}_{k\delta,(k+1)\delta}| \lesssim \\
    & |\tilde{X}(n)_{k\delta,(k+1)\delta} - \hat{X}_{k\delta,(k+1)\delta}|\delta^\frac{1}{p}
  + \sum_{j=1}^{k}|\tilde{X}(n)_{(j-1)\delta,j\delta}-\hat{X}_{(j-1)\delta,j\delta}|\delta^\frac{k+1-j}{p}.
\end{eqnarray*}
Thus, summing over all terms, we get
\begin{eqnarray*}
        \sum_{i=1}^{N_\delta} |\hat{X}^\delta_{(i-1)\delta,i\delta}-\tilde{X}^\delta(n+1)_{(i-1)\delta,i\delta}| &\lesssim& \sum_{i=1}^{N_\delta} |\hat{X}^\delta_{(i-1)\delta,i\delta}-\tilde{X}^\delta(n)_{(i-1)\delta,i\delta}|\left(\delta^\frac{1}{p}+\sum_{j=i}^{N_{\delta}-1}\delta^\frac{j-i+1}{p} \right) 
        \\&\leq& \sum_{i=1}^{N_\delta} |\hat{X}^\delta_{(i-1)\delta,i\delta}-\tilde{X}^\delta(n)_{(i-1)\delta,i\delta}|\left(\delta^\frac{1}{p}+\sum_{k=0}^{{N_\delta}-1-i}\delta^\frac{k+1}{p} \right) 
        \\ &\leq& \delta^\frac{1}{p}\left(1+\sum_{k=0}^{{N_\delta}-1}\delta^\frac{k}{p} \right)\sum_{i=1}^{N_\delta} |\hat{X}^\delta_{(i-1)\delta,i\delta}-\tilde{X}^\delta(n)_{(i-1)\delta,i\delta}|
\end{eqnarray*}
The result follows, as $\left(\delta^\frac{1}{p}+\sum_{k=0}^{{N_\delta}-1}\delta^\frac{k+1}{p} \right)\sim{\mathcal O}(\delta^\frac{1}{p})$.

Going back to the assumption that expansion \eqref{eq: X step inner segment} holds, it is now clear that we need the initial error $\|\tilde{X}(0)-\hat{X}\|_{p-var,[0,T]}$ to be small on the whole time horizon $[0,T]$ so that the differences $|\tilde{Y}(n)_{k\delta} -y_{k\delta}|$ remain uniformly under control with respect to both $n$ and $k$.

\end{proof}

We can now prove the main result of this section, which is the uniform convergence of $\tilde{X}(n)$ in $p$-variation.

\begin{theorem}
    Suppose that $\hat{X}^\delta$ is the solution to the discrete inverse problem defined in \ref{def: pwl discrete inverse problem} and the sequence $(\tilde{X}^\delta(n))_{n=0}^\infty$ is constructed by applying algorithm \ref{alg: sig}. Then, under the assumptions of lemma \ref{lemma:contraction}, there exists a $\delta_0>0$ such that 
    \begin{equation}
        \lim_{n\to \infty}\sup_{0<\delta<\delta_0} \|\tilde{\bf X}(n)-\hat{\bf X}\|_{p-var} = 0.
    \end{equation}
\end{theorem}
\begin{proof}
    It follows from lemma \ref{lemma:contraction} that 
    \[
    \sum_{k=1}^{N_\delta} |\tilde{X}(n)_{(k-1)\delta,k\delta}-\hat{X}_{(k-1)\delta,k\delta}| \sim {\mathcal O}(\delta^\frac{n}{p}).
    \]
    Similar to \eqref{eq: NR error bound}, it follows from bounds \eqref{eq: p-to-1 bound} and \eqref{eq: p-to-1 bound constant} that 
    \[
    \|\tilde{X}(n)-\hat{X}\|_{p-var}\lesssim\delta^{\frac{1-p}{p}\lfloor p \rfloor}\delta^\frac{n}{p},
    \]
    which implies uniform convergence, as required.
    \color{black}
\end{proof}

\section{Variations to the algorithm}
\label{sec: extensions}
\subsection{Alternative tree-like corrections}
\label{sec: extensions_1}

One of the main challenges in the implementation of the algorithm in proposition \ref{prop: alg} is the computation of the inverse It\^o map $I^{-1}$ applied to the tree-like corrections \eqref{r}. When $m=d$ and $f$ satisfies \eqref{eq: f assumptions}, it is possible to use \eqref{eq:inverse Ito}, assuming that the range of values of the tree-like corrections remains within a set where $g$ is locally Lipschitz (see also remark \ref{rem: algorith imprementation comments}). Even then, however, the computation of the integral in \eqref{eq:inverse Ito} might require further numerical techniques. One way to avoid that is to consider different types of tree-like corrections than $L(y,\tilde{y})$ -- note that the only requirement of a tree-like correction is to be a path that connects a point $y$ to a point $\tilde{y}$ and then retraces itself back to $y$. Thus, by choosing the tree-like correction carefully, we can simplify the computation.

In particular, we consider a correction that corrects each dimension sequentially, instead of all together. In spirit, this is the same idea underlying splitting algorithms, used in the simulation of complicated Stochastic Differential Equations (SDEs), for example in \cite{Foster} and \cite{splitting}. We demonstrate how this works on the following model:

\begin{align}
dY_t^{(1)}=&-aY_t^{(1)}dt+\sqrt{Y_t^{(1)}Y_t^{(2)}}dX_t^{(1)} \label{eq:sqrt1}\\
dY_t^{(2)}=&Y_t^{(2)}dX_t^{(2)} \label{eq:sqrt2}
\end{align}
This is a 3-dimensional system of the form \eqref{eq:main}, where the first dimension, denoted by $Y^{(0)}_t$, is the deterministic time component$$
\left(\begin{array}{ccc}
dY_t^{(0)} \\
dY_t^{(1)} \\
dY_t^{(2)}
\end{array}\right)=\left(\begin{array}{ccc}
1 & 0 & 0 \\
-aY_t^{(1)} & \sqrt{Y_t^{(1)}Y_t^{(2)}} & 0 \\
0 & 0 & Y_t^{(2)}
\end{array}\right) \left(\begin{array}{ccc}
dt \\
dX_t^{(1)} \\
dX_t^{(2)}
\end{array}\right)
$$
such that the inverse of $f$ we are required to integrate over is given by $$g(Y_t)=\left(\begin{array}{ccc}
1 & 0 & 0 \\
\frac{aY_t^{(1)}}{\sqrt{Y_t^{(1)}Y_t^{(2)}}} & \frac{1}{\sqrt{Y_t^{(1)}Y_t^{(2)}}} & 0 \\
0 & 0 & \frac{1}{Y_t^{(2)}}
\end{array}\right).$$ In its present form, the algorithm therefore requires the following integral evaluated at each interval: 
$$\bigintss_0^1 \frac{1}{\sqrt{\left(\tilde{Y}(n)_{k\delta}^{(1)}+u(Y^{(1)}_{k \delta}-\tilde{Y}(n)_{k\delta}^{(1)})\right)\left(\tilde{Y}(n)_{k\delta}^{(2)}+u(Y^{(2)}_{k \delta}-\tilde{Y}(n)_{k\delta}^{(2)})\right)}}d u . $$
This has a relatively complicated solution that is cumbersome to evaluate, slowing down the algorithm. Now, suppose that when constructing the tree-like corrections, instead of using a linear segment from $\tilde{Y}(n)_{k\delta}$ to $y_{k\delta}$, we instead correct each dimension one by one. i.e. we consider the tree-like correction that connects $\tilde{Y}(n)_{k\delta}$ to $Y_{k\delta}$ via $\tilde{Y}(n)_{k\delta}^{(1)}$, then $\tilde{Y}(n)_{k\delta}^{(2)}$ and so on up to the $d$\textsuperscript{th} dimension (see figure \ref{fig:enter-label}).
\begin{figure}[h!]
    \centering
    \includegraphics[width=0.7\textwidth]{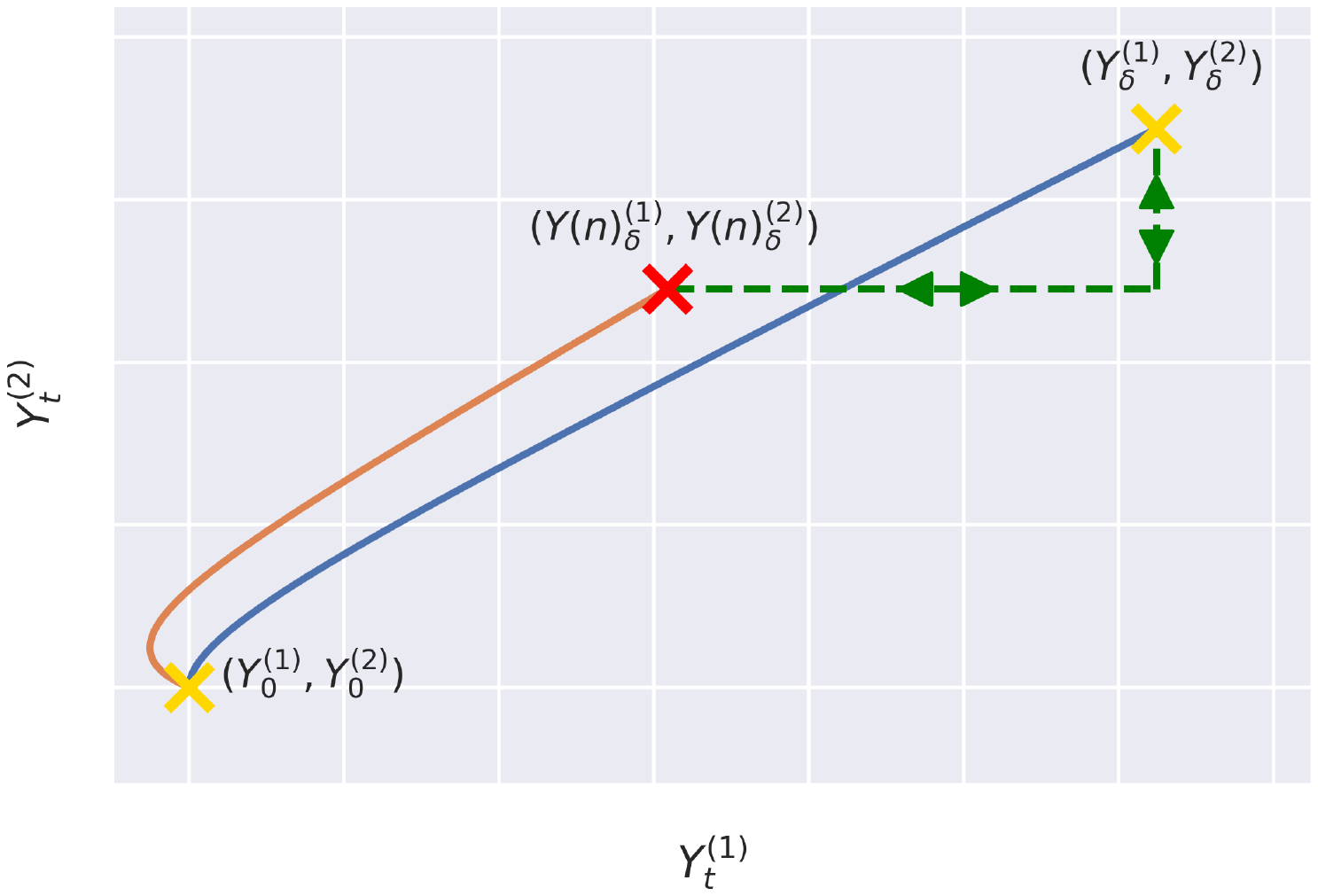}
    \caption{Visualisation of the correction.}
    \label{fig:enter-label}
\end{figure}

This new tree-like correction simplifies the update to $\tilde{c}(n)$. The proof of \ref{prop: alg} is unchanged, except the linear path between $\tilde{Y}(n+1)_{(k-1)\delta}$ and $y_{(k-1)\delta}$ is now the 'split' linear path, that corrects each dimension sequentially. Defining this path as $SL(\tilde{Y}(n+1)_{(k-1)\delta},y_{(k-1)\delta})$, the final line of the proof now looks like the following: 
\begin{eqnarray*}
    \tilde{X}^\delta(n+1)^1_{(k-1)\delta,k\delta} &=& X(n)^1_{(k-1)\delta,k\delta} - I^{-1}\left(SL(\tilde{Y}(n+1)_{(k-1)\delta},y_{(k-1)\delta})\right)^1 \\
    &&+ I^{-1}\left( SL(\tilde{Y}(n+1)_{k\delta},y_{k\delta})\right)^1.
\end{eqnarray*}

As before, we have:
\begin{eqnarray*}
     I^{-1}\left(SL(\tilde{Y}(n+1)_{(k-1)\delta},y_{(k-1)\delta})\right)= I^{-1}\left(\bigotimes_{i=1}^d L_i(\tilde{Y}(n+1)_{(k-1)\delta},y_{(k-1)\delta})\right),
\end{eqnarray*} 
where  $L_i(y_1,y_2)$ is the i'th segment of the split linear path. Concretely, the j'th coordinate of $L_i(y_1,y_2)$ is given by \begin{equation*}L_i(y_1,y_2)_{u}^{(j)}= 
\left\{
    \begin{array}{lr}
        y_2^{(j)}, & \text{if } j <i\\
        y_1^{(j)}+u(y_2^{(j)}-y_1^{(j)}), & \text{if } j=i\\
        y_1^{(j)}, & \text{if } j>i
    \end{array}\right.
\end{equation*} 
Therefore
\begin{eqnarray*}
    \tilde{X}^\delta(n+1)^1_{(k-1)\delta,k\delta} &=& X(n)^1_{(k-1)\delta,k\delta} - \sum_{i=1}^d I^{-1}\left(L_i(\tilde{Y}(n+1)_{(k-1)\delta},y_{(k-1)\delta})\right) \\
    &&+ \sum_{i=1}^d I^{-1}\left( L_i(\tilde{Y}(n+1)_{k\delta},y_{k\delta})\right).
\end{eqnarray*}

Now, computing the inverse It\^o map is much simpler, as all dimensions except one are constant: we must solve a one-dimensional integral of the form $\int g(a+bu)du$ where $g$, $a$ and $b$ are known. A simple Taylor expansion of $g$ around $u$ will give us a polynomial expansion of the integral, if we do not know the exact anti-derivative. Note that this method not only allows us to simplify \eqref{eq:inverse Ito}, which requires that $m=d$ but it can be used to construct inverse It\^o maps under weaker assumptions.

\subsection{Reducibility}
\label{sec: reducibility}

Another way of simplifying the computation of \eqref{eq:inverse Ito} is by using the reducibility property. Suppose that $g(y)$ has an antiderivative $\tilde{g}$. Then, we can write
\[
I^{-1}\left( Y(n) \right)_{s,t} = \int_s^t g(Y(n)_u)dY(n)_u = \int_s^t d(\tilde{g}(Y(n)_u)) = \tilde{g}(Y(n)_t) - \tilde{g}(Y(n)_s).
\]
In this case, the algorithm converges within one step, as the solution to the discrete inverse problem will only depend on ${\bf y}_{{\mathcal D}_\delta}$. Note that, in our notation, this would be equivalent to ${\bf I}^b_{1R}$ mapping any ${\bf P}_R(Z)$ to $0$ and convergence follows directly from \eqref{eq: correction term 2}. For the anti-derivative to exist, the system needs to satisfy a commutativity condition (or to be `reducible'; see \cite{Foster} or \cite{ait_sahalia} for the diffusive case). 

While this is usually not the case, it is often the case that the system is reducible for some of the dimensions. Thus, to increase the speed of the algorithm, we only need to iterate over the irreducible dimensions (the ones not satisfying the commutativity condition). When paired with appropriately chosen tree-like corrections, this can substantially speed up the algorithm.

\section{Numerical Examples}
\label{sec: examples}

We compare the classical approach of subsection \ref{sec: classic} to the algorithm presented in subsection \ref{sec: signature}. We consider three different model settings, which cover a range of modeling paradigms. Initially, we consider the model mentioned in section \ref{sec: extensions_1} specified via equations \eqref{eq:sqrt1} and \eqref{eq:sqrt2}, under two separate configurations.  We show convergence over this 2-dimensional irreducible process, where the piecewise linear driving path is drawn from the increment process of a fractional Brownian motion (fBm) to highlight the potential in rougher, irreducible regimes. We further extend this to a case where the path of $Y_t$ also depends on the history of $Y_t$. Secondly, we detail a model from opinion dynamics, where a number of dimensions are unobserved, and affect the process through the drift coefficient. Finally, we use an exact simulation algorithm from financial modeling to illustrate how the solutions of the discrete inverse problem converge to those of the continuous inverse problem, highlighting the power of the method in approximating complex, multivariate stochastic models.      

It is rare that we can exactly simulate from the stochastic model, so in each of the scenarios below (aside from the final one), we simulate the solution to the ODE \eqref{eq: general ODE} at a very fine time-step $\delta_{sim}<<\delta_{sample}$, with fixed initial condition, and then sample this ODE at a coarser time-step $\delta_{sample}=\delta=\frac{T}{N}$. Here, we use the increments of an fBm to build this path, however the technique would work with other paths. We thus assume that the increments of the driving path, $X^\mathcal{D}_t$, are distributed according to fractional Gaussian noise, i.e. with covariance structure
$$\mathbb{E}(X^{\mathcal{D}}_{i\delta}X^{\mathcal{D}}_{j\delta})=\frac{\delta^{2h}}{2}\left(i^{2h}+j^{2h}-|i-j|^{2h}\right),$$
where $i\delta$ and $j\delta$ are elements of the partition $\mathcal{D}$. The Hurst parameter $h\in (0,1)$ determines the roughness of the noise and coincides with a Brownian motion at $h=0.5$. For this analysis, we assume that an $n$-dimensional fBm is just a set of $n$ pairwise independent fBm's, ensuring that all dimension-based correlation is incorporated by the model $f(Y_t)$, although we see no reason why the method could not work by a structure of dependence across different noise paths.  

Crucially, to mimic the discretisation of the true rough differential equation, the piecewise noise must converge in the p-variation metric to the continuous one, which only occurs in the case of fBm for $h>0.25$ (see \cite{TerryBook}), so we restrict the analysis to this setting.  Recall that we are still using Algorithm~\ref{alg: sig} to infer a driving path $X^{\mathcal{D}}_t$ that is piecewise linear over the \textit{sampling grid} ${\mathcal D} = {\mathcal D}_\delta$, rather than the simulation grid ${\mathcal D}_{\delta_{sim}}$.


\subsection{Irreducible 2-dimensional process}
In this section, we show how the method introduced works on a simple, irreducible example, under two settings. Firstly, where the drift parameter is homogeneous; it does not change in time. Secondly, to illustrate the power of the method on more complex paths, we show how the method can work in a setting where the drift parameter alters throughout the path's trajectory. 

\subsubsection{Fixed drift}
\label{sec: sqrt}

\begin{figure}
\centering
    \subfloat[\centering Varying sampling rate (averaged over roughness)]{{\includegraphics[width=7.4cm]{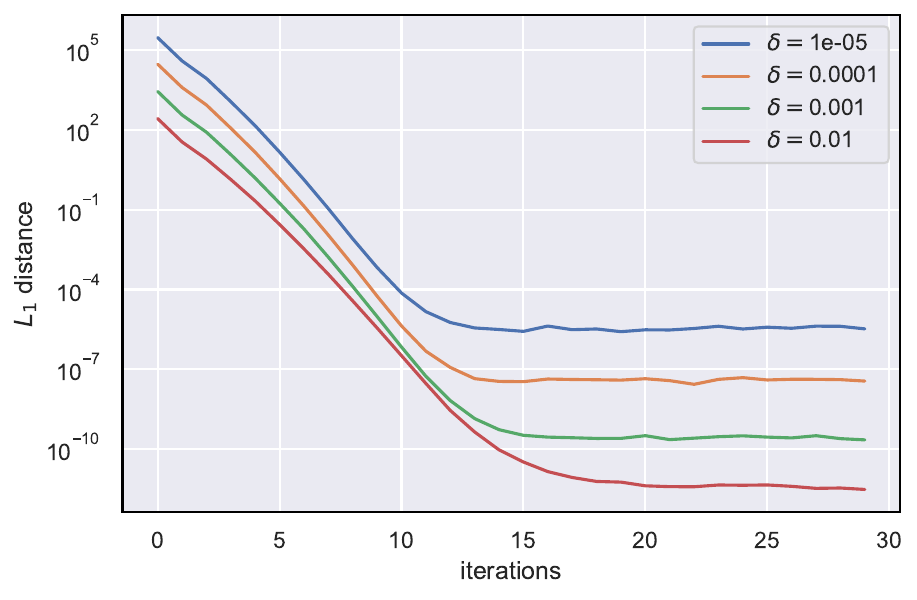}}}
    \subfloat[\centering Varying roughness (averaged over sampling rate)]{{\includegraphics[width=7.4cm]{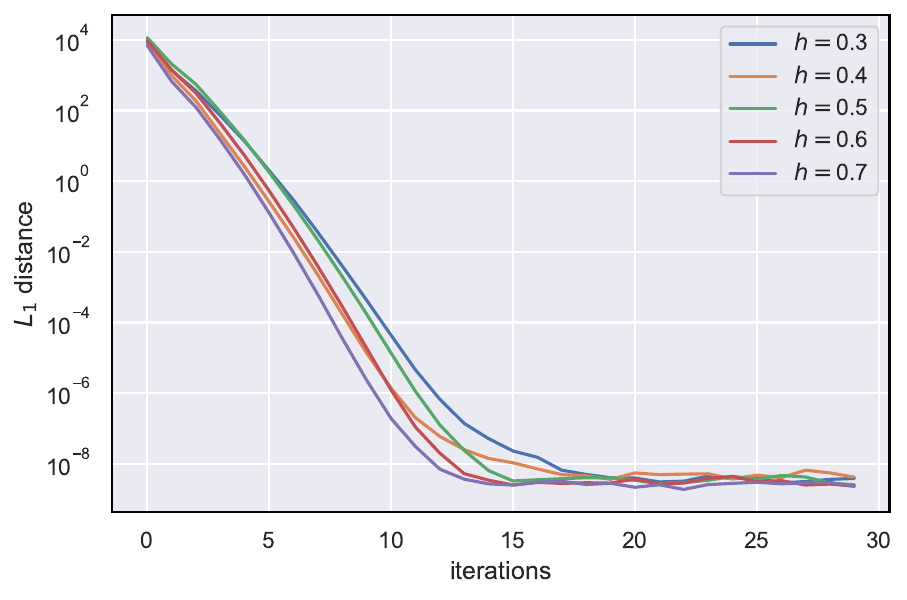}}}
    \caption{Convergence in $L_1$ distance averaged over different iterations, while varying roughness and sampling rate.}
    \label{fig:fixed_drift}
    
\end{figure}

We first consider the algorithm on the toy model specified by equations \eqref{eq:sqrt1} and \eqref{eq:sqrt2} of Section \ref{sec: extensions}, which does not admit a transform to unit variance. Here, it is not possible to write the process $X_t$ in terms of the increments of the response $Y_t$, as it depends on the iterated integrals of $Y_t$ which we do not observe. Therefore, we initialize $X_t$ by setting $Y_t$ to be linear between the observations $Y^\mathcal{D}$, as mentioned in Algorithm~\ref{alg: sig}. Further, as the problem is immediately solved over the second dimension (which is reducible), we are only required to update the first dimension, reducing computational time. To compare the sample points to the process driven by $X^\delta(c^\delta(n))$, we use the $L_1$ distance over the sample points (taking the log for comparability -- see remark \ref{rem: L1 error bound} for justification), against computation time (in seconds)
$$\log\left(\sum_{t_i\in\mathcal{D}}|Y_{t_i}-I(X^\delta({c}^\delta(n)))_{t_i}|\right).$$In Figure \ref{fig:fixed_drift} we show the results. Firstly, in Figure \ref{fig:fixed_drift}(A) we show how convergence is affected by the sampling rate. Here we have averaged the results across paths driven by fBMs with different levels of roughness.  This plot shows that the algorithm converges uniformly across $\delta$, with all paths converging to the correct sample points at around 12 iterations. The gap we see between each $\delta$ is representative of the increase in computer error occurring in each interval, as $N$ increases by a factor of 10 between successive interval estimations.

In Figure \ref{fig:fixed_drift}(B), we visualise the same set of results, but average across the sampling rate $\delta$, to examine how the driver roughness affects the convergence of the scheme. As expected, the Algorithm~\ref{alg: sig} takes longer to converge for rougher signals, but it still converges for an fBm with $h=0.3$ within 20 iterations of the algorithm. 

\subsubsection{Path-dependent drift}
\label{sec:path_dependent_drift}
We now consider a more challenging setting for the inverse algorithm by letting the drift $a_t$ depend on the past trajectory of $Y_t$. We let $a_t$ switch between positive (attraction to home state) and negative (repellent), depending on whether the trajectory of $Y_t$ has hit an upper barrier $u$ and a lower barrier $l$, respectively. We are therefore looking to find the set of $N\times 2$ gradients ${c}^\delta_x = \left( c^{\delta,k}_{x,i} \right)_{i=1,k=1}^{N,2}$ that solve the equations
\begin{align}
dY_t^{(1)}=&-a_tY_t^{(1)}dt+\sqrt{Y_t^{(1)}Y_t^{(2)}}c_{x,i}^{\delta,1}dt \label{eq:sqrt_unobs1},\\
dY_t^{(2)}=&Y_t^{(2)}c_{x,i}^{\delta,2}dt, \label{eq:sqrt_unobs2}
\end{align}
in each interval $t \in [i\delta,(i+1)\delta]$, with initial condition $Y_{i\delta}=y_{i\delta}$ and terminal condition $Y_{(i+1)\delta}=y_{(i+1)\delta}$. The drift $a_t$ is initialized at $0$. If the path hits the upper barrier $u$, $a_t$ becomes $1$, while if it hits the lower barrier $l$, $a_t$ becomes $-1$. Such an example could represent an animal switching from home location attraction to repellent behavior. We must therefore simulate $a_t$, along with the path $Y_t$.

This presents an issue for local algorithms such as the Newton-Raphson (NR) algorithm, as we cannot, for a given ${c}^\delta(n)$, calculate $D_c F(\delta; Y_{k\delta},c^\delta(n)_k)$ or even $F(\delta; Y_{k\delta},c^\delta(n)_k)$, as $F(\cdot)$ now depends on $Y_t$ for $t\in[0,k\delta]$ in equation (\ref{eq:NR}). We consider four approaches:
\begin{enumerate}
    \item A naive approach is to set $a=0$ in the calculation of $D_c F$ and $F$. This is fast to simulate, as we can simulate over the coarse grid, but induces a sizeable error.
    \item We can improve the accuracy of approach (1) by first simulating the full process on the simulation grid, i.e. using ${c}^\delta(n)$, we retrieve $Y_t({c}^\delta(n))$ and $a_t({c}^\delta(n))$. Now, in each coarse interval, for $t \in [i\delta,(i+1)\delta]$, we set $a_t=a_{i\delta}({c}^\delta(n))$. This is more accurate than approach (1) as we use the dynamics of $X_t({c}^\delta(n))$, but with a restriction to remove the requirement of simulating twice at the fine time-step (here we simulate at the fine level once, and the coarse level once). However, over the intervals where the path crosses the barrier, we induce an error.
    \item Repeat approach (2) but now allow $a_t$ to change over each interval when computing $D_c F$ and $F$, i.e. we simulate the system once, in order to get $a_t({c}^\delta(n))$. We then simulate $D_c F$ and $F$ inside each interval over the fine grid, with initial conditions $a_{i\delta}=a_{i\delta}({c}^\delta(n))$. This is more expensive as we need to simulate the system over the fine grid twice. Due to an application of the chain rule, we also need to calculate the derivative of the solution with respect to the initial condition $y_{i\delta}$. Let $m$ be the number of intervals in the fine grid within each interval of the coarse grid. To calculate the gradient at $(i+\frac{j+1}{m})\delta$, given the path $Y$ and gradient $G$ at $(i+\frac{j}{m})\delta$, one must iteratively apply the following formula:
    \begin{align*}
    G\left(\frac{(j+1)\delta}{m},Y_{i\delta}\right)=&G\left(\frac{j\delta}{m},Y_{i\delta}\right)\cdot D_{Y_{i\delta}}F\left(\frac{\delta}{m},Y_{(i+\frac{j}{m})\delta}(c^\delta(n))\right) \\
        &+ G\left(\frac{\delta}{m},Y_{(i+\frac{j}{m})\delta}(c^\delta(n))\right).
    \end{align*}
    \item We repeat approach (3), but consider that we cannot find the true derivative with respect to the initial condition or with respect to $c$. Thus, we are required to solve the initial value problem corresponding to \eqref{Z equation} to gather $D_c F$, which adds a layer of computation on top of the simulation. This mimics the situation where we can obtain an accurate simulator, but are required to gather a numerical or approximate derivative, which is often expensive.
\end{enumerate}

We initialize the algorithms with ${c}^\delta(0)={0}$, as we need to simulate the system forwards to get an initial estimate for the path $a_t$. We then simulate the true path five times, with $\delta_{sim}=0.00125$, and sample the data at increasing rates. Thus, we apply the above algorithms on different sample sizes, namely $N=20,50,100,200$. The ability to parallelise some of the computations has a further affect on computational efficiency. Due to the results of Proposition~\ref{prop: alg}, and the splitting scheme detailed in Section~\ref{sec: extensions}, we can perform the computation of $r(n)$ across multiple cores. This entails a maximum of $N d^2$ calculations. The NR approaches all involve a simulation step (not parallelisable), and an update using a gradient calculation. This gradient calculation requires the simulation of all dimensions, so cannot be parallelised more than across the number of points in the sample grid. We take an approach that gives us a holistic estimate for the effects of this parallelisation, by dividing the running time of the update by the number of cores (we assume 8 across the simulations).

\begin{figure}
    \centering
    \subfloat[\centering $\delta_{sample}=0.05; \ N=20$]{{\includegraphics[width=7.4cm]{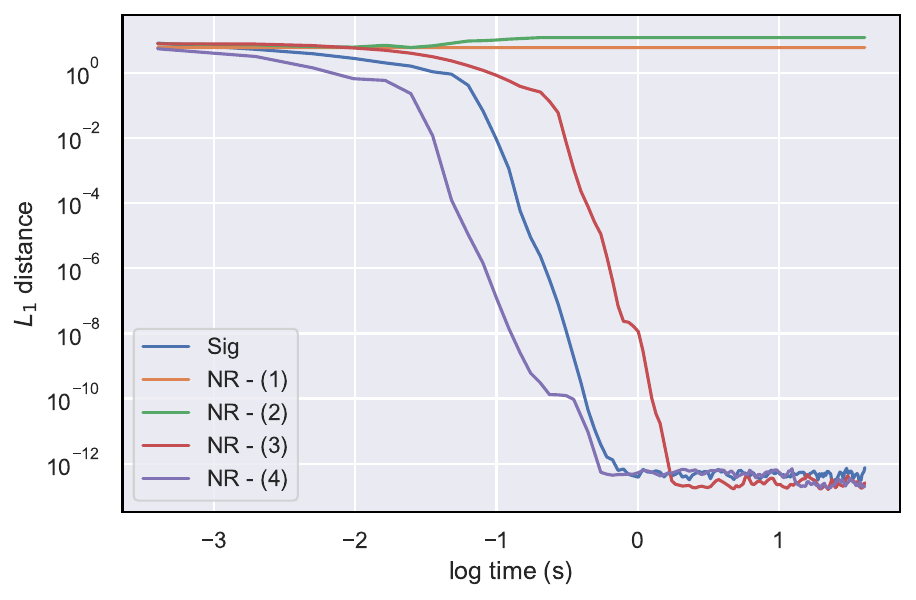}}}
    \subfloat[\centering $\delta_{sample}=0.02; \ N=50$]{{\includegraphics[width=7.4cm]{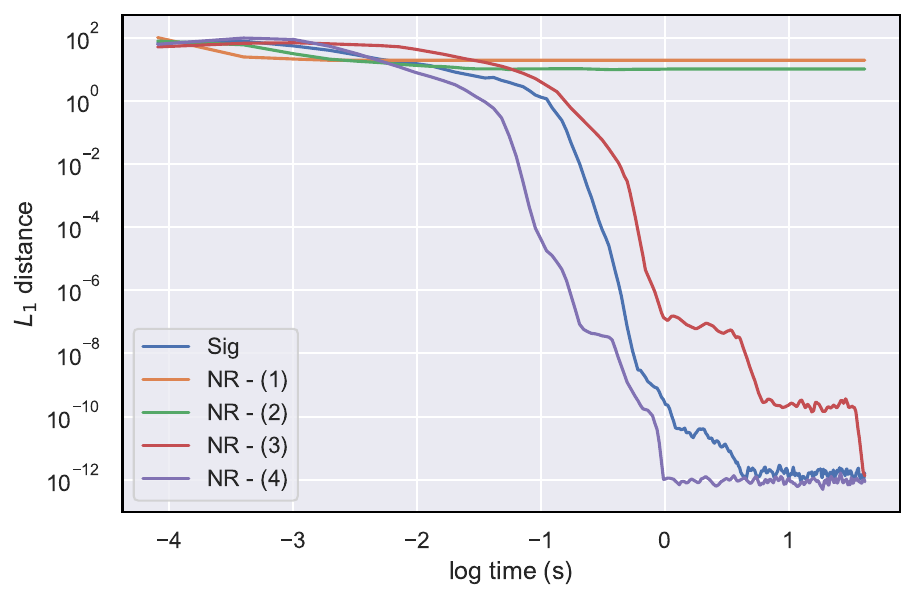}}}
    \\
    \subfloat[\centering $\delta_{sample}=0.01; \ N=100$]{{\includegraphics[width=7.4cm]{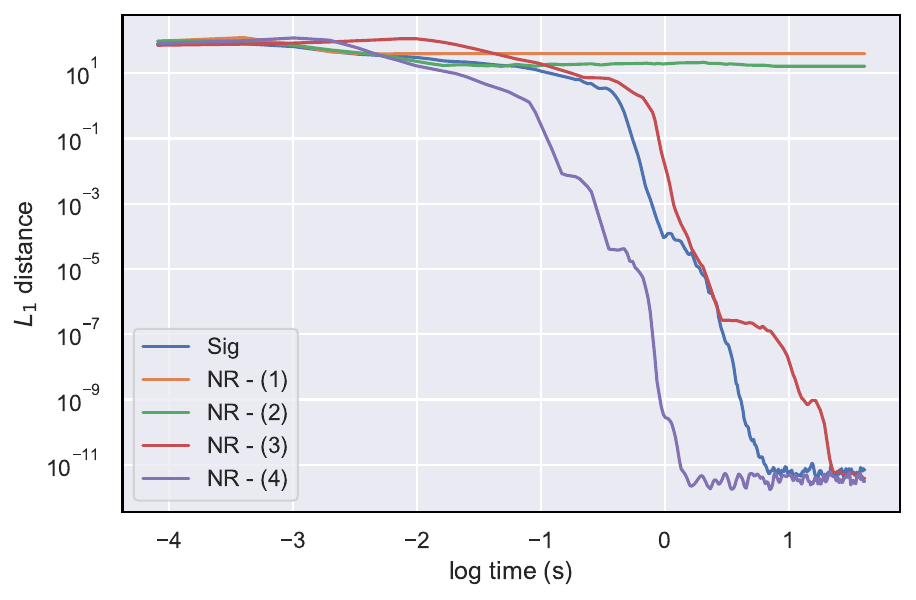}}}
    \subfloat[\centering $\delta_{sample}=0.005; \ N=200$]{{\includegraphics[width=7.4cm]{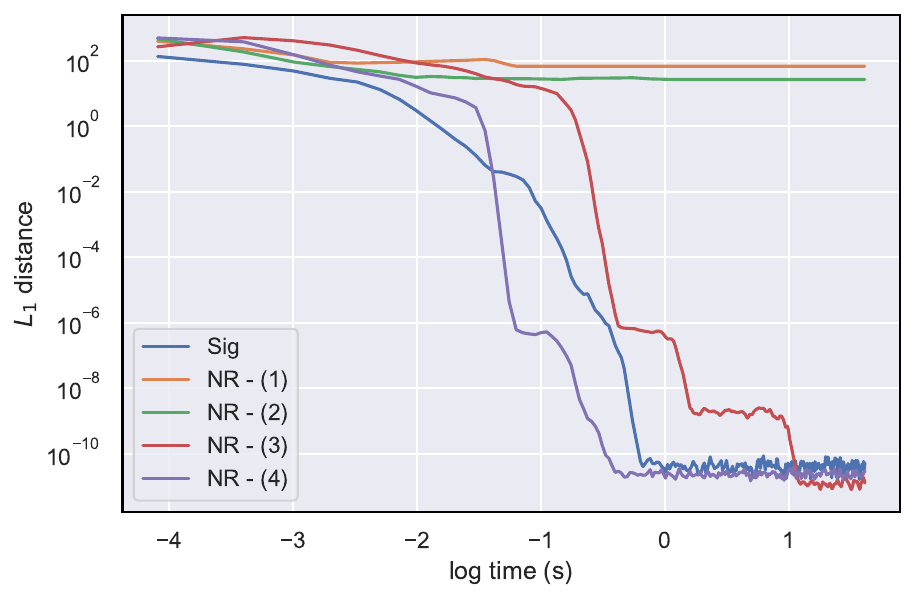}}}
    \caption{Response convergence over different sampling rates with  $\delta_{sim}=0.00125; \ h=0.3; \ l=2.8;u=3.2;a=5$, averaged across 5 simulations, assuming parallelisation across 8 cores.}
    \label{fig:sqrt_switching}
\end{figure}

First, we can see in Figure~\ref{fig:sqrt_switching} that naive approaches (1) and (2) converge fast, but do not capture the original data due to missing details on the model. The NR algorithm converges faster than the signature detailed in algorithm~\ref{alg: sig} in all simulations related to approach (3). This is expected since the NR algorithm uses more information from the model when exact derivatives are used, and therefore can target the correct value faster. The signature algorithm reaches the true path faster than NR in approach (4), where the derivative must be approximated. Furthermore, as $N$ decreases, the signature algorithm becomes more efficient compared to NR as the proportion of time that both algorithms spend on simulation compared to updating ${c}^\delta(n)$ decreases. This is a benefit for the signature algorithm, as, once simulated, it can update ${c}^\delta(n)$ faster than the NR algorithm. 

\begin{figure}
    \centering
    \includegraphics[width=9cm]{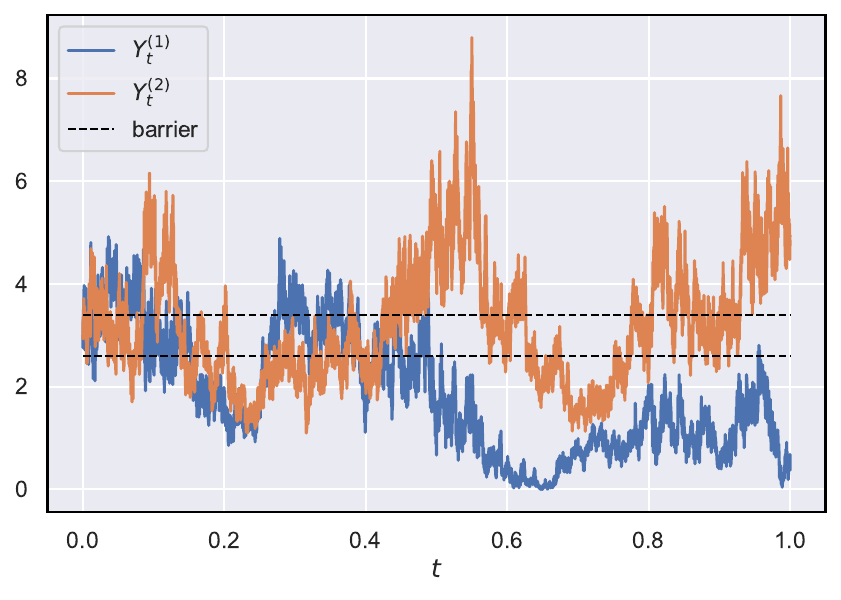}
    \caption{Example sample path of $Y^{\mathcal{D}}_t$ with parameters $h=0.3;\delta=1\times 10^{-4};l=2.6;u=3.4;a=3$}
    \label{fig:enter-label}
\end{figure}

\subsection{Opinion dynamics with unobserved trajectories}
\label{sec: opinion}

\begin{figure}
    \centering
    \includegraphics[width=9cm]{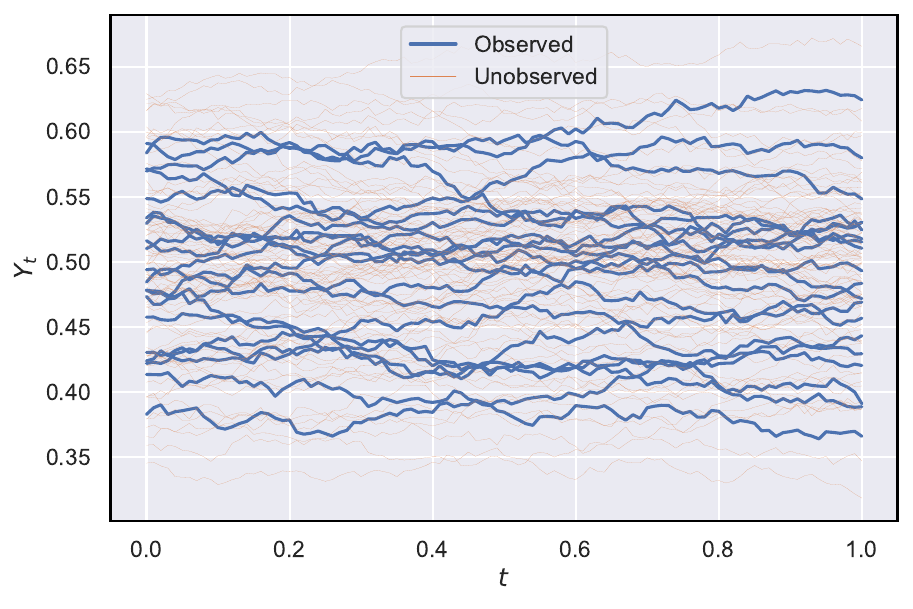}
    \caption{Example opinion dynamics for $\sigma=0.1, \ R=0.1$, with $d_{o}=20, \ d=100$}
    \label{fig:opinion_example}
\end{figure}

We now give an example which has path dependence through unobserved trajectories. Below is an adaptation of a model for opinion dynamics from \cite{opinion}, where we have enforced a controlling term in the variance, to prevent the particles from hitting their boundaries too quickly.  Particle $i$ moves with dynamics represented by the following differential equation:
\begin{equation}
\label{eq: opinion}
    dY_t^{(i)} = -\frac{1}{d}\sum_{j: \ |Y_t^{(i)}-Y_t^{(j)}|<R}(Y_t^{(i)}-Y_t^{(j)}) dt + \sigma \sqrt{Y_t^{(i)}(1-Y_t^{(i)})}dX_t^{(i)},
\end{equation}
for $i,j=1,\dots,d$. In reality, only a number of these trajectories would be observed. We assume that we observe the trajectories in the index set $d^\mathcal{O}:=\{1,2,\hdots,d_o\}$ out of $d$ trajectories, without loss of generality if the initial starting points are independent and random. In Figure~\ref{fig:opinion_example} we show an example sample path for 100 particles, where we may observe only 20 of the trajectories.   

As alluded to in Remark~\ref{rmk: uneven dimensions} (1), we now have a number of degrees of freedom in the unobserved trajectories, so in both the signature and NR algorithm, we set ${c}^{\delta,k}(0)=0$ for $k \in d^\mathcal{U}:= \{d_o+1,d_o+2,\hdots,d\}$, and let the unobserved trajectories $Y^{\mathcal{U}}$ be free to change to ensure a fit of the observed trajectories. In reality, one may use a Kalman filter-like approach to update the unobserved driving paths $c^{\mathcal{U}}$.  This is one approach for the NR algorithm; another would be to ignore the unobserved paths fully, running the algorithm with only the observed paths. As $d_o$ grows with respect to $d$, this approaches the true system. 

At this point, we consider how the algorithms compare in efficiency, and in what scenario the signature algorithm can outperform a local approach. With a given ${c}^\delta(n)$, the signature algorithm takes two steps of computation to perform an update: 
\begin{enumerate}
    \item Simulate the path $I(X^\delta({c}^\delta(n)))$ (forward step).
    \item Calculate a maximum of $d^2 N$ integrals (backward step).
\end{enumerate}
Step (1) cannot be parallelised, but step (2) can. The integrals in step (2) can be calculated simultaneously, so with enough parallel cores, one could calculate all $d^2 N$ integrals in one step. Therefore, the additional computation after the simulation of step (1) is minimal with enough cores. We estimate this effect using the same methodology as detailed in Section~\ref{sec:path_dependent_drift}.

The NR algorithm requires two steps of computation at each time-step:
\begin{enumerate}
    \item Calculate $F(\delta; Y_{k\delta},c^\delta(n)_k)$.
    \item Calculate $D_c F(\delta; Y_{k\delta},c^\delta(n)_k)$.
\end{enumerate}
If there is no path dependence or dependence on unobservable variables, then both steps can be parallelised over the $N$ steps. However, in the case of hidden paths, we do not have access to the full $Y_{k\delta}$, so cannot accurately carry out step (1) or (2). We therefore consider a few options for updating ${c}^\delta(n)$, in the case of hidden trajectories: 
\begin{enumerate}
    \item Simplify the model, by ignoring the unobservable dimensions. This essentially replaces $F$ with an approximation $\tilde{F}$ in the computations above (if $R=0$, i.e. there is no interaction between particles, then it is exact). This has three efficiency benefits: there is no need to simulate the full system to update ${c}^{\delta, \mathcal{O}}(n)$, we can parallelise the computations over the $N$ time steps, and we only need to solve an initial value problem of dimension $d_o+d_o^2$ . 
    \item   Fix ${c}^{\delta,\mathcal{U}}(n)=0, \ \forall n$, and simulate $I (X^\delta({c}^\delta(n)))$ to get $Y_{k\delta}^{\mathcal{U}}, \ \forall k\in \{1,\hdots,N-1\}$, then carry out steps (1) and (2) of the NR algorithm (parallel in $N$). This is less efficient, as we need to simulate the full system first, then calculate $d+d^2$ solutions to the initial value problem for both $F$ and $G$.
    \item  Fix ${c}^{\delta,\mathcal{U}}(n)=0, \ \forall n$, but instead of simulating the system as in case (2) right above, we use the output of the previous time-step computation as an initial condition for the next computation. i.e. we update $ c^{\delta,\mathcal{O}}(n)_{2}$ by computing $F\left(\delta;\{Y^{\mathcal{O}}_\delta,F^\mathcal{U}_1\},c^\delta(n)_1\right)$ and $D_cF\left(\delta;\{Y^{\mathcal{O}}_\delta,F^\mathcal{U}_1\},c^\delta(n)_1\right)$,where $F_1$ are the unobservable dimensions of the output to step (1) of NR calculated in the previous time-step, so in this case $F_1:=F(\delta;Y_0,c^\delta(n)_0))^{\mathcal{U}}$. This is, in spirit very close to case (2) right above, but we now cannot parallelise the update over $N$, as each step of the NR algorithm requires the update from the previous interval. However, there is also no need to simulate the system initially in this situation. Both case (2) and (3) use more information from the underlying model than case (1) at the expense of efficiency.
 \end{enumerate}

Thus, the signature approach is computationally beneficial in two scenarios, namely when simulation is fast and dimensions are large, or where we have fixed computer resources and cannot use $N$ cores to compute local (NR) updates. As we do not have access to the solution of differential equation (\ref{eq: opinion}), we use a numerical initial value problem solver to simulate the system for a given piecewise linear path $X^{\mathcal{D}}_t$. We are not necessarily concerned with sampling rate here or the intricacies of exact simulation, so we simulate the system over the same grid as we perform the sampling and we only focus on the simulations that do not touch the boundaries to examine the effectiveness of the algorithm. The signature algorithm only improves with the ability to exactly simulate.


We simulate the system 3 times over $[0,1]$, with $\delta_{sim}=0.001$, using a piecewise linear path with firstly Brownian increments, and then a slightly rougher path where we draw increments from a fBm with $h=0.4$. We further use a differing number of initial particles $d\in\{10,50,100,200\}$. Once simulated, we sample the data at increasing intervals, so that $N\in\{10,100,200,500\}$ and alter the number of observed paths (we assume that we observe a proportion of the paths, letting $d_o/d\in\{0.2,0.5,0.8\}$). We set $\sigma=0.06, \ R=0.1$. We sample the initial positions from a $\mathcal{N}(0.5,0.005^2)$ distribution. To run the signature and NR algorithms (cases (2) and (3)), we assume that we do observe the initial positions of all the variables. 
The metric of convergence we use here is the $L_2$ distance, but we standardise by the number of observed dimensions, so that we can compare across path dimensions.  

\begin{figure} 
    \centering
    \includegraphics[width=\textwidth]{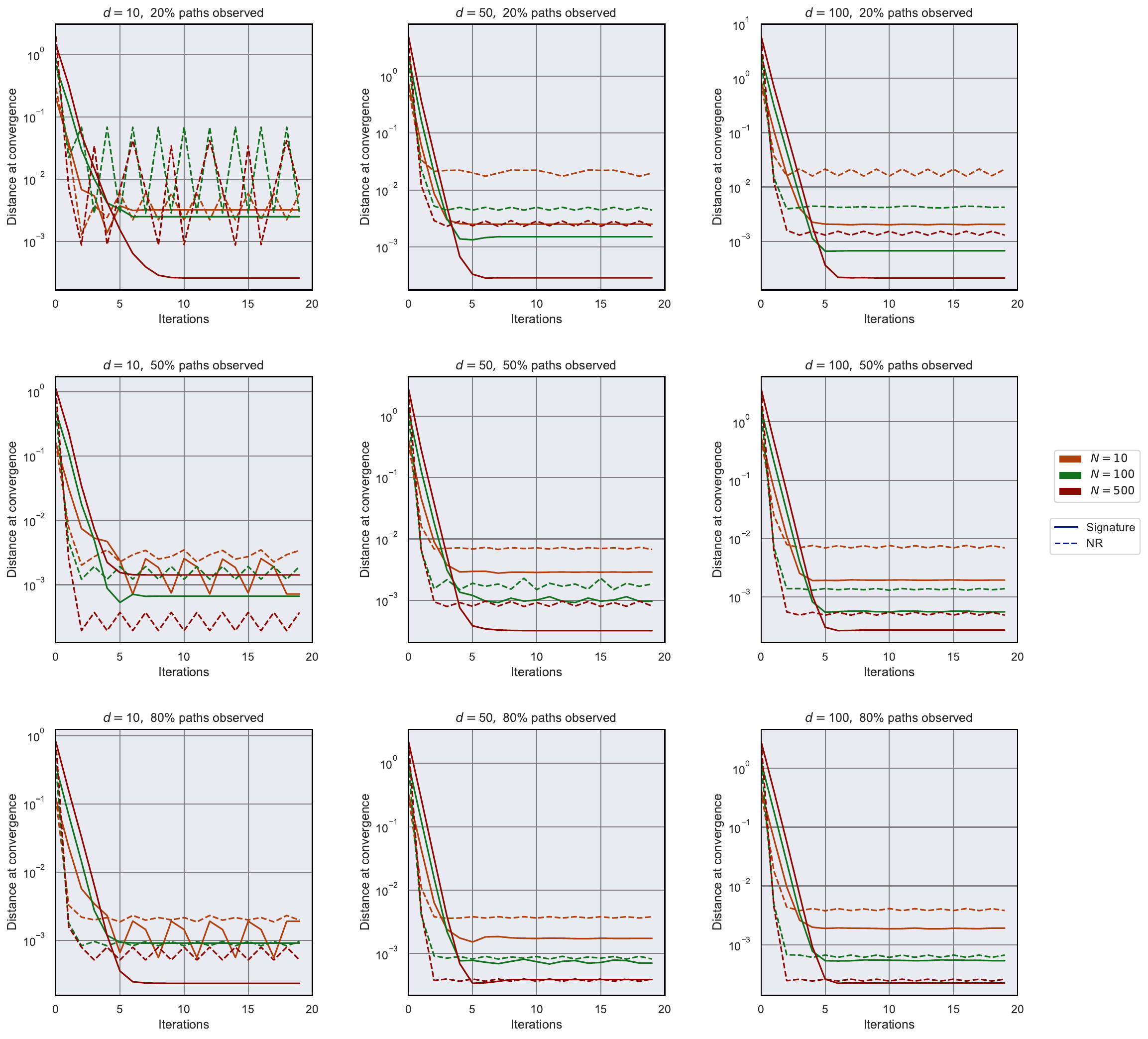}
    \caption{Comparison between signature and naive NR (case (1)) for Brownian driving path across different sampling rates}
    \label{fig:opinion_bm}
\end{figure}

In Figure~\ref{fig:opinion_bm}, we show how the paths converge to the sampled data in the Brownian case, compared to the 'naive' attempt at NR detailed above in case (1). As expected, the signature algorithm converges to a closer path than the naive attempt, in all situations. It does best in comparison in the case where we have more dimensions, and less observed paths. This is expected, as it allows for interactions with the unobserved paths, so can tailor the observed responses to the correct ones using a truer reflection of the underlying model.

\begin{figure} 
    \centering
    \includegraphics[width=\textwidth]{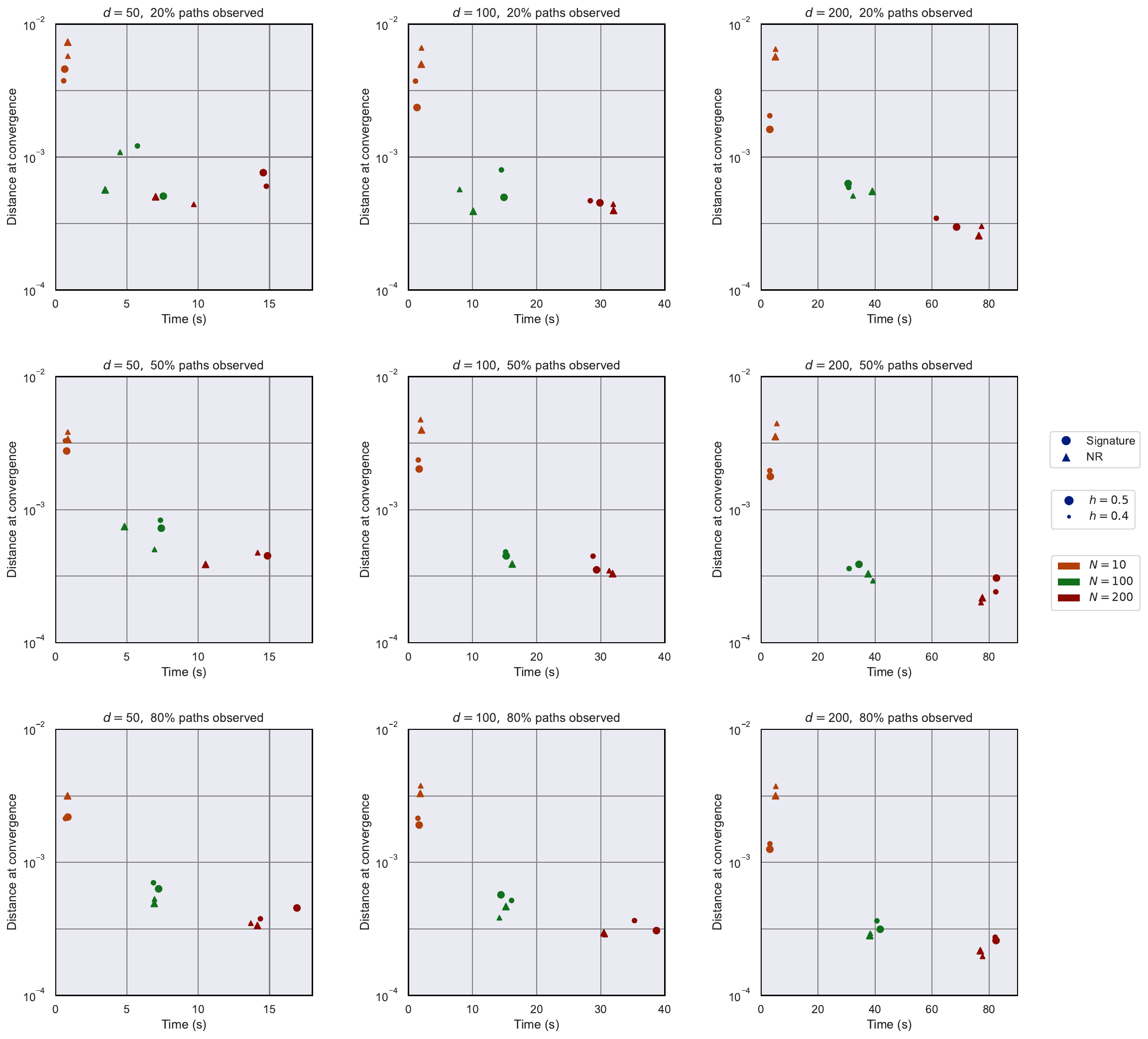}
    \caption{Convergence comparison for signature and NR (case (2)) across different model roughness and sampling rates}
    \label{fig:opinion_fbm}
\end{figure}

In Figure~\ref{fig:opinion_fbm}, we compare to case (2) of the NR algorithm (note that case (3) is not included, as we found it to behave similar to case (2), but was more inefficient). We further include the system simulated with a rougher path, i,e. where increments are drawn from the increment distribution of the fBm with $h=0.4$. We plot how quickly both algorithms converge in computation time (assuming the ability to parallelise over a maximum of 8 cores), and to what level they converge to (all of the paths converged to the data in these simulations). It is clear that the signature performs best relative to the NR in the case with more dimensions, and less data (e.g. $d=200$, $N=10$ and $d_o/d=0.2$). However, as we increase the availability of data, by increasing the observed paths and with a higher sampling rate, the NR improves relative to the signature. This illustrates how the method that relies more on the model, can attain a faster convergence when it has a more accurate realization of the model.   


\subsection{Convergence to the continuous inverse problem}
\label{sec:cts_inverse_example}

Finally, we show how signature algorithm detailed in Algorithm~\ref{alg: sig} can be used as a robust approach to approximate the solution to the continuous inverse problem. Thanks to the algorithm detailed in \cite{li_wu_exact} we can simulate \textit{exactly} the SDE system:
\begin{align}
dY_t^{(1)}=&\mu Y_t^{(1)}dt+Y_t^{(1)}Y_t^{(2)}\left(\sqrt{1-\rho^2}dW_t^{(1)}+\rho dW_t^{(2)} \right)\label{eq:sqrt_unobs1},\\
dY_t^{(2)}=&-\kappa(\theta-Y_t^{(2)})dt+\xi dW_t^{(2)}.\label{eq:sqrt_unobs2}
\end{align}

The model above is the Ornstein-Uhlenbeck Stochastic Volatility (OUSV) model for the joint movement of stock prices and volatility. Volatility is usually unobserved, but we are interested in the behavior of the solutions to the discrete inverse problem as we gain more data, so here we assume we can observe the stochastic volatility component. In reality, one could use a Bayesian algorithm, such as a particle filter, to simulate forward the volatility process.     

We first apply a Stratonovich correction term to the drift (so that the piecewise linear approximation converges) and then simulate the system at a very fine timescale, proceeding to sample the process over grids of varying size. We then calculate the p-variation distance between the exact finely sampled path and $I(X^\delta({c}^\delta(n)))$ over different $\delta$; the exact algorithm is constructed for an It\^{o} diffusion, so we require $p>2$.

\begin{figure} 
    \centering
    \subfloat[\centering p-variation distance to piecewise projection of exact sampled $y_t$ at the 20th iteration across sampling rate $\delta$]{{\includegraphics[width=7.4cm]{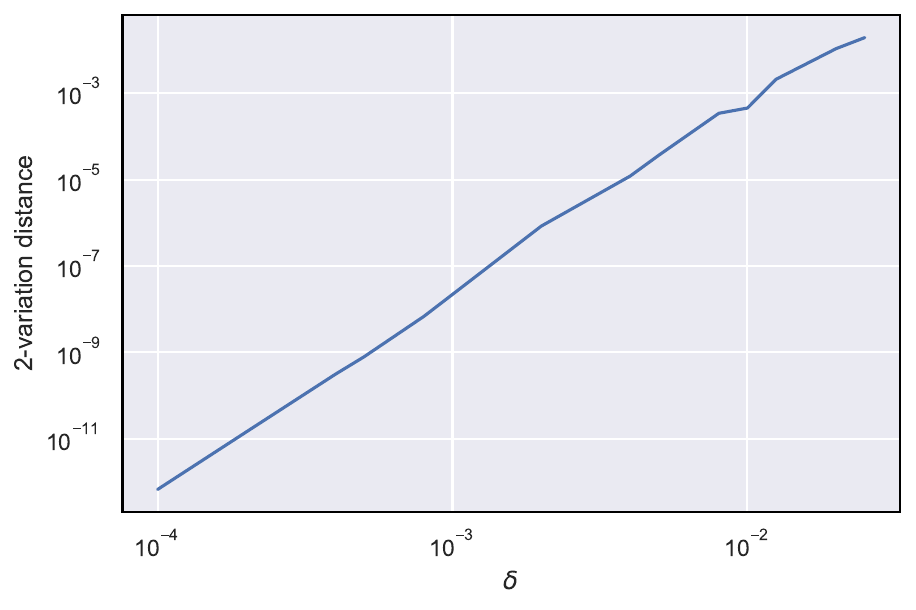}}}
    \subfloat[\centering Convergence to piecewise projection of exact sampled $y_t$ across $\delta$, as $n$ grows]{{\includegraphics[width=7.4cm]{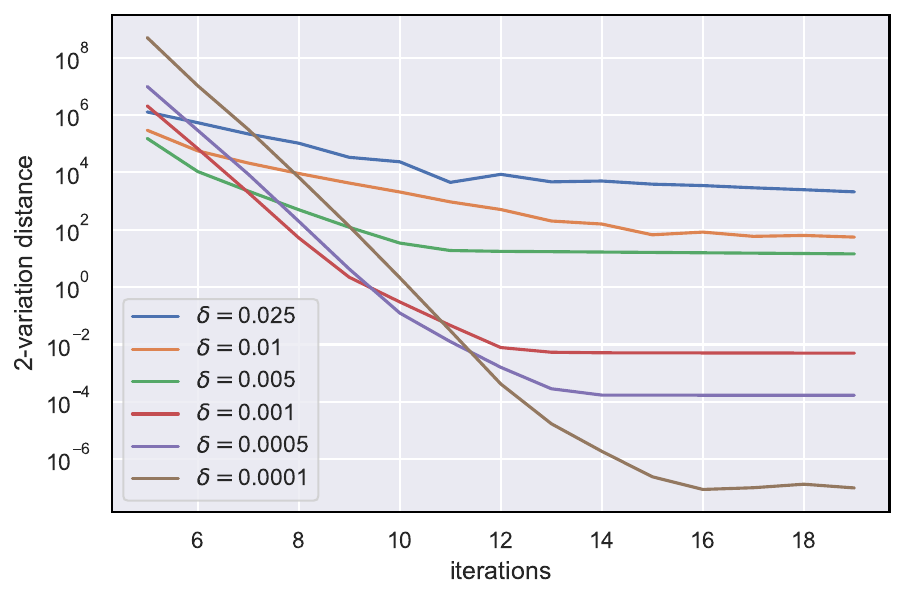}}}
    \\
    \centering \subfloat[\centering Approximate (due to sampling) p-variation convergence of $I(X^\delta({c}^\delta(20)))$ to the exact trajectory $y_t$ as $\delta\rightarrow 0$]{{\includegraphics[width=7.4cm]{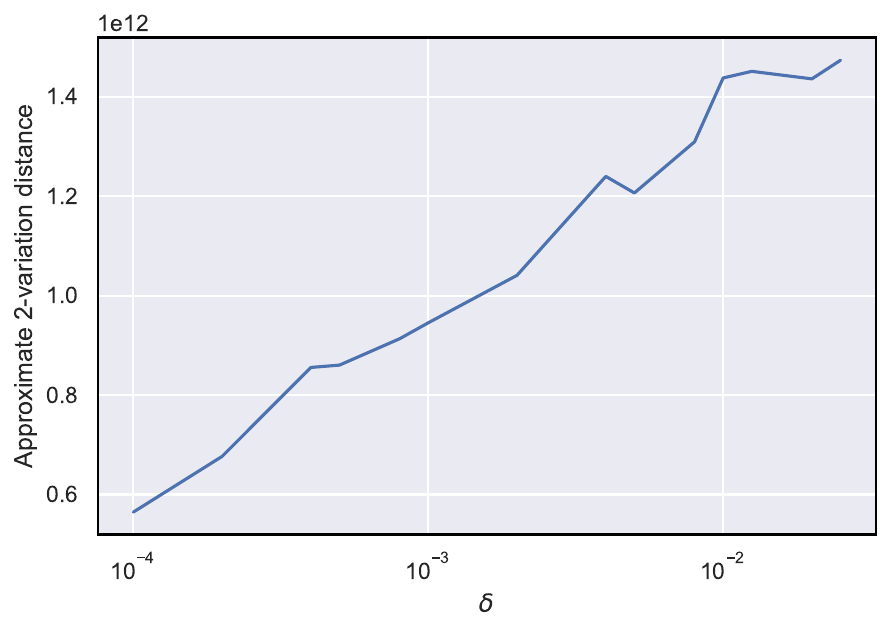}}}
    \caption{Convergence of solutions in the p-variation metric as $\delta\rightarrow 0$.}
    \label{fig:exact}
\end{figure}

The plots in Figure \ref{fig:exact} display empirically that, firstly, we have convergence of the algorithm to the projected exactly sampled piecewise trajectory of the response, $\pi^{\mathcal{D}}(Y_t)$, and further, convergence in p-variation of the continuous response to the piecewise linear $X_t^{\delta}({c}^\delta(n)$ to the actual trajectory $Y_t$ (the latter fact is impossible to verify with finite computational resources so we use an approximate fine partition). Moreover, we note that although the finer partition is initialized further away than in the coarser partition (local errors are cumulative), convergence occurs at the same rate in each delta.

\subsection{Conclusion from the Numerical Examples}

We have introduced a signature-based algorithm that computes the set of gradients of a piecewise-linear path which drives the observations $y$ via the map $I(\cdot)$. The algorithm is very straightforward to implement for a given model. When the model is known and simulation is possible, the signature method eliminates the need to compute derivatives of the solution explicitly, as required by forward methods such as Newton-Raphson. Instead, it requires the evaluation of an integral involving simulation outcomes and the vector field $f(\cdot)$. Although it offers a slower convergence on a single interval compared to methods such as the Newton Raphson procedure, it considers the entire path in each iteration, enabling global control of the error. In its current form, the algorithm requires a `split' piecewise linear path which is advantageous if we are aiming to construct a likelihood for the data that is consistent with the model and the observations. However, in principle, one could change the assumption of $\Tilde{X}(n)$ being piecewise-linear to coming from a different family of approximations to Gaussian processes (e.g. Karhunen-Lo\'eve approximation), so long as the inverse map $\frac{1}{\delta}I^{-1}\left( L(\tilde{Y}(n)_{k\delta},y_{k\delta}) \right)$ is still tractable.  

\begin{acks}[Acknowledgments]
T.~M. is supported by the Warwick Statistics Department Centre for Doctoral Training, and gratefully acknowledges funding from the University of Warwick and the UK Engineering and Physical Sciences Research Council (Grant number: EP/T51794X/1). A.~P. acknowledges the support of the UK Engineering and Physical Sciences Research Council (Grant number: EP/W00707X/1).

The authors would like to thank the anonymous referee for their comments and suggestions that significantly improved the paper.
\end{acks}


\begin{thebibliography}{99}

\bibitem{ait_sahalia}
Y. Aït-Sahalia. Closed-Form Likelihood Expansions for Multivariate Diffusions.
{\it The Annals of Statistics} 36(2): 906-937, 2008.

\bibitem{BailleulDiehl}
I. Bailleul and J. Diehl. The Inverse Problem for Rough Controlled Differential Equations. {\it SIAM J. Control Optim.} 53(5): 2762--2780, 2015.

\bibitem{Bailleul}
I. Bailleul. Regularity of the It\^o-Lyons map. {\it Confluentes Math.}, 7(1):3--11, 2015. 

\bibitem{Baudoin}
F. Baudoin. {\it Rough Path Theory}. {\tt https://fabricebaudoin.blog/wp-content/uploads} {\tt /2017/08/rough.pdf}

\bibitem{Horatio}
H. Boedihardjo, Xi Geng, T. Lyons and DanyuYang. The signature of a rough path: Uniqueness. {it Advances in Mathematics}: 293: 720-737, 2016.

\bibitem{splitting}
E. Buckwar, A. Samson, M. Tamborrino, I. Tubikanec. A splitting method for SDEs with locally Lipschitz drift: Illustration on the FitzHugh-Nagumo model.
{\it Applied Numerical Mathematics},
Volume 179,
2022, 191-220,
{\tt https://doi.org/10.1016/j.apnum.2022.04.018}

\bibitem{Chevyrev}
I. Chevyrev. Rough Path Theory. {\tt arXiv: 2402.10331}

\bibitem{Craigmile}
P. Craigmile, R. Herbei, Ge Liu and G. Schneider. Statistical inference for stochastic differential equations. {\it WIREs Computational Statistics, 15(2)}, e1585. {\tt https://doi.org/10.1002/wics.1585}, 2023.

\bibitem{CoutinLejay}
L. Coutin and A Lejay. Sensitivity of rough differential equations: An approach through the Omega lemma. {\it J. Differential Equations} 264: 3899--3917, 2018.

\bibitem{Driver}
B. Driver. {\it Rough Path Analysis}
{\tt https://mathweb.ucsd.edu/~bdriver/247A-Winter2012/} {\tt Lecture\%20Notes/rps.pdf}

\bibitem{Flint}
G. Flint and T. Lyons. Pathwise approximation of SDEs by coupling piecewise abelian rough paths. {\tt arXiv: 1505.01298}

\bibitem{FrizVictoir}
P. Friz and N. Victoir. {\it Multidimensional stochastic processes as rough paths: theory and applications}, vol.120. CUP, 2010.

\bibitem{Foster}
J. Foster, G. dos Reis and C. Strange. High order splitting methods for SDEs satisfying a commutativity condition. {\it SIAM J. Num. An.} , 62(1): 500-532, 2024. {\tt arXiv:2210.17543}

\bibitem{Gluckstad}
M. Gl\"uckstad, N. Muca Cirone and J. Teichmann. Signature Reconstruction form Randomiased Signatures {\tt arXiv:2502.03163}



\bibitem{HuNualart}
Y. Hu and D. Nualart. Parameter estimation for fractional Ornstein-Uhlenbeck processes, {\it Statistics \& Probability Letters}, 80(11-12): 1030--1038, 2010. 

\bibitem{kelley2010}
W. G. Kelley and A. C. Peterson. {\it The Theory of Differential Equations}, (Second Edition), Springer, 2010.

\bibitem{Kidger}
P. Kidger, J. Foster, Xuechen Li and T Lyons. Efficient and Accurate Gradients for Neural SDEs. {\it NeurIPS}, 2021b. {\tt arXiv:2105.13493}

\bibitem{KubiliusSkorniakov}
K. Kubilius and V. Skorniakov. On some estimators of the Hurst index of the solution of SDE driven by a fractional Brownian motion, {\it Statistics \& Probability Letters}, 109:159--167, 2016. 

\bibitem{Kutoyants}
Y. A. Kutoyants. {\it Statistical inference for ergodic diffusion processes}. Springer Science \& Business Media, 2013.

\bibitem{li_wu_exact}
C. Li and L. Wu. Exact simulation of the Ornstein–Uhlenbeck driven stochastic volatility model {\it European Journal of Operational Research} 275(2): 768-779, 2019.

\bibitem{loeve_book}
M. Loève. {\it Probability Theory I}. Springer New York, 1977.

\bibitem{TerryBook}
T. Lyons and Z. Qian. {\it System Control and Rough Paths} OUP, 2002.




\bibitem{Darrick}
P. Semnani, V. Guan, E. Robeva and D. Lee. Path-Dependent SDEs: Solutions and Parameter Estimation. {\tt arXiv:2505.22646}

\bibitem{Unterberger}
J. Unterberger. A rough path over multidimensional fractional Brownian motion with arbitrary Hurst index by Fourier normal ordering. {\it Stochastic Processes and their Applications} 120: 1444--1472, 2010.

\bibitem{opinion}
C. Wang, Q. Li, W. E and B. Chazelle. Noisy Hegselmann-Krause Systems: Phase transition and the 2R-conjecture. {\it Journal of Statistical Physics} 166(5): 1209–1225, 2017.

\bibitem{ZhangXiao}
P. Zhang, W.-L. Xiao, X.-L. Zhang, P.-Q. Niu. Parameter identification for fractional Ornstein-Uhlenbeck processes based on discrete observation. {\it Economic Modelling} 36: 198--203, 2014.


\end{thebibliography}
\end{document}